\documentclass[bj,numbers]{imsart}

\usepackage{amsmath,amsfonts,amssymb,mathtools}
\usepackage{graphicx,subcaption}

\startlocaldefs
\numberwithin{equation}{section}

\theoremstyle{plain}
\newtheorem{theorem}{Theorem}[section]
\newtheorem{lemma}[theorem]{Lemma}

\theoremstyle{definition}
\newtheorem{assumption}{Assumption}[section]
\newtheorem{definition}[theorem]{Definition}

\newtheorem*{remark}{Remark}

\newcommand{\R}{\mathbb{R}}
\newcommand{\Rd}{\mathbb{R}^d}
\newcommand{\norm}{\|\cdot\|}
\newcommand{\Z}{\mathcal{Z}}
\newcommand{\X}{\mathcal{X}}
\newcommand{\E}{\mathbb{E}} 

\newcommand{\ltwo}{L^2}

\newcommand{\nprior}{\Pi_N}
\newcommand{\nlan}{\Pi_N^{LAN}}
\newcommand{\dnlan}{\Pi_N^{DLAN}}
\newcommand{\data}{\mathcal{D}^N}
\newcommand{\ptheta}{P_{\theta_0}^N}

\newcommand{\Hs}{\widetilde{H}}
\newcommand{\opinfo}{\mathbb{I}_{\theta_0}}
\newcommand{\opcomb}{\mathbb{K}_{\theta_0,\Lambda}}
\newcommand{\barepsu}{\overline{\epsilon}_{N,0}}
\newcommand{\barepsc}{\overline{\epsilon}_{N,1}}
\newcommand{\barpsi}{\overline{\psi}_N}
\newcommand{\bigpsi}{\widehat{\Psi}_N}
\newcommand{\lpsi}{\psi_\lambda}
\newcommand{\bartheta}{\overline{\theta}_N}
\newcommand{\thetat}{\theta^{(t)}}
\newcommand{\loss}{\ell_N}
\newcommand{\llan}{\ell_N^{LAN}}
\newcommand{\dllan}{\ell_N^{DLAN}}
\newcommand{\stoch}{\mathcal{F}_N}
\newcommand{\etal}{\underline{\eta}}
\newcommand{\etau}{\overline{\eta}}
\newcommand{\kappal}{\underline{\kappa}}
\newcommand{\kappau}{\overline{\kappa}}

\endlocaldefs

\begin{document}

\begin{frontmatter}
\title{Frequentist Coverage of Bayes Posteriors in Nonlinear Inverse Problems with Gaussian Priors}

\runtitle{Frequentist Coverage of Bayes Posteriors in Inverse Problems}

\begin{aug}

\author[A]{\inits{F.}\fnms{Youngsoo}~\snm{Baek}\ead[label=e1]{youngsoo.baek@duke.edu}\orcid{0000-0002-0143-2839}}
\author[B]{\inits{S.}\fnms{Katerina}~\snm{Papagiannouli}\ead[label=e2]{aikaterini.papagiannouli@unipi.it}\orcid{0000-0002-4488-4843}}

\address[A]{Department of Biostatistics and Bioinformatics, Duke University, Durham, NC\printead[presep={,\ }]{e1}}

\address[B]{Department of Mathematics, University of Pisa, IT\printead[presep={,\ }]{e2}}
\end{aug}

\vspace{5mm}
\begin{center}

\emph{In the memory of Sayan Mukherjee}

\end{center}
\begin{abstract}
We study asymptotic frequentist coverage and approximately Gaussian properties of 
Bayes posterior credible sets in nonlinear inverse problems when a Gaussian prior is placed on the parameter of the PDE. 
The aim is to ensure valid frequentist coverage of Bayes credible intervals when estimating continuous linear functionals of the parameter. Our results show that Bayes credible intervals have conservative coverage 
under certain smoothness assumptions on the parameter and a compatibility condition between the likelihood and the prior, regardless of whether an efficient limit exists or Bernstein von-Mises (BvM) theorem holds. 
In the latter case, our results yield a corollary with more relaxed sufficient conditions than previous works. The theory is illustrated with 
a PDE that arises in predicting the transport of radioactive waste from underground repositories and optimizing oil recovery from subsurface fields: 
an elliptic inverse problem for Darcy flow. 
In this case, a near-$1/\sqrt{N}$ contraction rate and conservative coverage results are obtained for linear functionals that were shown not to be estimable efficiently.
\end{abstract}


\end{frontmatter}

\section{Introduction} \label{sec:intro}

We study a general class of statistical nonlinear inverse problems arising from partial differential equations (PDEs). 
Such problems occur when an unknown function $\theta \in \Theta$, defined on a spatial domain $\Z \subset \R^d$, can only be observed indirectly through a nonlinear forward operator
$G: \Theta \to L^{2}(\mathcal{X})$,
where $\Theta$ is a Hilbert space and $L^{2}(\mathcal{X})$ denotes the space of square-integrable functions on $\mathcal{X} \subset \R^d$, where observations are collected subject to measurement noise. 
Typical examples arise in medical imaging, geophysical fluid dynamics, and climate modeling, where $\theta$ encodes an unknown physical parameter of the system.

In this paper, we consider nonlinear inverse problems under the statistical model
\begin{equation} \label{eq:model}
Y_i = G(\theta_0)(X_i) + \varepsilon_i, \qquad i=1,\ldots,N,
\end{equation}
where $\theta_0 \in \Theta$ is the unknown ground-truth parameter, the design points $X_i \in \mathcal{X}$ represent sampling locations, and the noise variables $\varepsilon_i$ are independent Gaussian random variables. The sample size is denoted by $N \geq 1$.

In the class of problems studied here, the forward operator $G$ is given by the parameter-to-solution map of a PDE describing the underlying physical system. The parameter $\theta_0$ typically enters the PDE as a coefficient function, and $G(\theta_0)$ represents the observable solution. The primary inferential task is to recover $\theta_0$ from the data $(Y_i, X_i)_{i=1}^N$. 
In many applications, however, an estimate of $\theta_0$ is not sufficient. The modern paradigm of uncertainty quantification (UQ) requires statistically valid measures of uncertainty derived from the data.

A widely used approach to UQ in inverse problems is Bayesian inference. Influential works such as \cite{stuart2010inverse,dashti2013bayesian} advocated the use of Bayesian methods in PDE-based inverse problems, building on earlier developments \cite{tarantola2005inverse,kaipio2006statistical}. In practice, posterior distributions are often explored using Markov chain Monte Carlo (MCMC) algorithms \cite{cotter2013mcmc}. 
From a theoretical perspective, much of the literature has focused on posterior consistency and contraction rates, that is, whether the posterior concentrates around the true parameter $\theta_0$ as $N \to \infty$, and at what rate \cite{giordano2020consistency, monard2021consistent}.

A more delicate question concerns the \emph{frequentist coverage} of Bayesian procedures: whether posterior credible sets provide valid confidence intervals for the unknown parameter. Initial results in this direction were obtained for linear inverse problems \cite{knapik2011bayesian, castillo2013nonparametric, castillo2015bernstein}. More recently, \cite{nickl2020bernstein} and \cite{monard2021statistical} developed Bernstein--von Mises (BvM) type results for nonlinear inverse problems, including PDE models such as the Schr{\"o}dinger equation; see also \cite[Sections~2--3]{nickl2023bayesian}.

Roughly speaking, the BvM theorem asserts that, as $N \to \infty$, the posterior distribution becomes approximately Gaussian, centered around an efficient estimator and with variance given by the inverse Fisher information. As a consequence, posterior credible sets asymptotically achieve correct frequentist coverage. However, in infinite-dimensional settings, BvM may fail \cite{freedman1999wald}, so we follow \cite{nickl2020bernstein} and \cite{monard2021statistical} in studying whether these properties are met for certain functionals of $\theta$. 
In particular, it is known that when estimating a functional of form $\langle \theta, \psi \rangle$ for $\psi\in\Theta$ and when $\psi$ does not belong to the range of the adjoint of an ``information operator'' (to be rigorously defined below), efficient estimation at the $1/\sqrt{N}$ rate may be impossible, leading to a breakdown of BvM \cite{nickl2022some}. These observations raise the following fundamental question.

~

\begin{center}
\textbf{Frequentist coverage problem for Bayesian posteriors.}
Consider the nonlinear regression model \eqref{eq:model}. 
Set a target level $\gamma$ (e.g., 5\%).
Can we find a Bayesian posterior credible set $I_N$ for estimating continuous linear functionals $\langle \theta, \psi \rangle$, containing at least $(1-\gamma)$-posterior mass, such that:

\begin{enumerate}
\item \textit{asymptotically correct}  coverage,
\[
\liminf_{N\to \infty} P_{\theta_0}^{N}\big(\langle \theta_0,\psi\rangle \in I_N\big)\geq 1-\gamma;\text{ or,}
\]

\item \textit{asymptotically exact} coverage,
\[
\lim_{N\to \infty} P_{\theta_0}^{N}\big(\langle \theta_0,\psi\rangle \in I_N\big)= 1-\gamma,
\]
even in regimes where Bernstein--von Mises theorems fail?
\end{enumerate}
\end{center}

~

The goal of this paper is to study this question for a broad class of nonlinear inverse problems with Gaussian priors. A key difficulty arises from the fact that BvM-based techniques are not applicable in general. In particular, for PDE-based inverse problems such as Darcy flow, the failure of BvM is closely tied to the absence of efficient estimators. Indeed, classical results \cite{van1991differentiable} show that efficient estimation of a functional requires a compatibility condition between the functional and the score operator. When this condition fails, neither efficient estimation nor BvM can hold \cite{nickl2022some}.
This motivates the need for new tools that do not rely on BvM. Our main contribution is to develop a general framework for establishing frequentist coverage of Bayesian credible sets in such settings. The key technical ingredient is a novel use of \emph{variational inequalities}, which allow us to control the interaction between the forward operator and the prior covariance without relying on spectral diagonalization like in \cite{knapik2011bayesian}.

As a motivating example, we study the Darcy flow problem, which models subsurface flow and has been widely used in geophysical applications \cite{stuart2010inverse}. Specializing to $\R^2$ domains, the goal is to estimate the coefficient of a second-order divergence form elliptic PDE
\begin{equation} \label{eq:darcy} 
-\nabla\cdot(a(s,t) \nabla u(s,t)) = f(s,t) 
\quad \text{where} \quad 
(s,t)\in \mathcal{X}\subset\mathbb{R}^2,
\end{equation}
subject to problem-specific boundary conditions on $u$.
Here $a$ is the unknown parameter, a function with a domain $\X$ that is referred to as the ``conductivity'' or ``diffusion'' coefficient. 
We get to observe finitely discretized measurements of $u = u(a)$, contaminated with observation errors.
The description suits the observation model \eqref{eq:obs_model}, provided the uniqueness of $u(a)$ for each $a$, so that the resulting implicit mapping $G$ associates each $a$ to $u(a)$.
If uniqueness and regularity are not guaranteed for $a$ in a linear space, one can smoothly reparametrize the problem in terms of a linear space variable $\theta$.

Our results show that Bayesian credible intervals can still achieve valid frequentist coverage, even in regimes where BvM fails. Under suitable smoothness and compatibility conditions, we prove that posterior credible sets are asymptotically conservative, and in some cases exact, while achieving near-optimal contraction rates. These results demonstrate that the failure of BvM does not necessarily compromise uncertainty quantification in nonlinear inverse problems.

\subsection{Our Contribution}
Our main contribution is to develop a general framework for analyzing the frequentist coverage of Bayesian posterior credible sets in nonlinear inverse problems of the form \eqref{eq:model}, under Gaussian priors on $\Theta$. We then apply this framework to show the validity of posterior credible intervals in UQ for the Darcy flow problem.

The failure of BvM in Darcy flow has a fundamental connection to the \emph{lack of an efficient estimator}. By efficiency, we refer to the phenomenon that a statistical estimator converges to its estimand at a $1/\sqrt{N}$-rate and attains the ``best possible variance,'' which should be the inverse of Fisher information, based on the Cram{\'e}r-Rao lower bound 
(cf. \cite[Section~2]{kullback1959information} and \cite[Chapter~8]{van2000asymptotic}). \cite{van1991differentiable} studied making sense of such a statement for estimating functionals of an infinite-dimensional parameter belonging to a normed linear space. A consequence of \cite[Theorems~3.1,~4.1]{van1991differentiable} is that the ``efficient information'', whose inverse is a lower bound to the variance of all $1/\sqrt{N}$-rate estimators of a (not necessarily linear) functional $\Psi(\theta)$ is non-zero if and only if a suitably defined gradient of $\Psi$ belongs to the range of the adjoint of a certain \emph{score operator}. \cite[Section~7]{van1991differentiable} already contained examples demonstrating not all functionals satisfy such a property. \cite[Theorems~1--2]{nickl2022some} related the score operator for the likelihood defined by \eqref{eq:obs_model} to a derivative $\opinfo$ of $G$ at $\theta_0$ (defined in Section \ref{sec:setup}), before using this relation to show both the lack of estimators convergent to the true functional at $1/\sqrt{N}$-rate and the failure of BvM in Darcy flow when $\psi$ does not belong to the range of the adjoint of $\opinfo$.

In this paper, we focus on inference for \emph{linear} functionals $\Psi(\theta) = \langle \theta, \psi \rangle$, and establish conditions under which posterior credible intervals provide valid frequentist uncertainty quantification, even in regimes where Bernstein--von Mises (BvM) theorems fail.
More precisely, we show that:

\begin{enumerate}
\item Under suitable regularity conditions on the forward operator $G$, the prior covariance, and the functional $\psi$, posterior credible intervals achieve \emph{asymptotically correct} frequentist coverage.

\item Even in the absence of efficient estimators (and hence failure of BvM), posterior distributions of $\Psi(\theta)$ contract around $\Psi(\theta_0)$ at rates that depend on the smoothness of $\theta_0$ and increasingly approach $1/\sqrt{N}$.
\end{enumerate}

A key technical challenge is that classical approaches rely on spectral analysis of the Fisher information operator or diagonalization arguments, which are generally unavailable in nonlinear PDE-based inverse problems, unlike in linear inverse problems where the forward operator admits a suitable spectral decomposition \cite{knapik2011bayesian}. 
To overcome this difficulty, we introduce a novel use of \emph{variational inequalities}, which allow us to control the interaction between the forward operator and the prior covariance without requiring explicit spectral decompositions.

\begin{definition}[Variational inequality]
\label{def:vi}
Let $X, Y$ be Hilbert spaces and let $A : X \to Y$ be an injective bounded linear operator with non-closed range. 
Let $f : [0,\infty)\to[0,\infty)$ be continuous and strictly increasing, with $f(0)=0$. 
We say that $x \in X$ satisfies a \emph{variational inequality} with constant $b > 0$ and modifier $f$ if, for all $u \in X$,
\[
b \|x - u\|_X^2 
\leq 
\|x\|_X^2 - \|u\|_X^2 + f(\|Ax - Au\|_Y).
\]
\end{definition}

In our setting, this inequality is applied to operators of the form 
\[
\Lambda^{-1/2}\opinfo,
\]
where $\opinfo$ denotes the linearization of the forward map $G$ at $\theta_0$, and $\Lambda$ is the inverse of the Laplace operator which encodes the prior covariance structure. The variational inequality provides indirect control over the spectrum of the operator $\Lambda^{-1/2}\opinfo^*\opinfo\Lambda^{-1/2}$, which governs the fluctuations of the posterior distribution. 
We can draw on several consequences established in the literature, which we summarize in Appendix \ref{appendix:vi_lemmas}. 

As we show in Section \ref{sec:darcy_flow}, previous results like \cite[Theorem~5]{nickl2022some} can be used to verify this inequality when $\theta_0$ is smooth enough. The novelty of our application is that we use an inequality in Definition \ref{def:vi} to study the frequentist coverage of Bayes posteriors, whereas previous authors used it to derive rates of convergence for estimators in deterministic inverse problems. Under reasonable conditions stated in terms of smoothness of $\theta_0$ and $\psi$, and a local injectivity estimate on the Cameron-Martin space of the prior, we prove both an asymptotically conservative, sometimes exact, coverage of Bayes posterior quantile-based interval estimates and an upper bound on posterior convergence rate to $\theta_0$, closer to $1/\sqrt{N}$ for a smoother truth $\theta_0$. 

This approach replaces classical spectral arguments and yields the key estimates required to establish frequentist coverage. In particular, it allows us to treat nonlinear inverse problems where diagonalization techniques are not available. As a consequence, we show that Bayesian credible intervals provide valid frequentist confidence sets even in settings where BvM fails due to the lack of efficient estimators.

\section{Problem setup} \label{sec:setup}

\subsection{Notation} \label{sec:notation}

For two sequences $a_N$ and $b_N$, $a_N = O(b_N)$ (also $a_N \lesssim b_N$) implies there exists a constant $C$ such that $|a_N| \leq Cb_N$ for all sufficiently large $N$. 
$a_N \asymp b_N$ when $a_N \lesssim b_N$ and $a_N \gtrsim b_N$. 
$a_N = o(b_N)$ (also $a_N \ll b_N$)when for sufficiently large $N$, $|a_N| \leq \varepsilon b_N$ for any constant $\varepsilon > 0$.
$Law(Y)$ is a shorthand for the probability measure for a random variable $Y$.
For two sequences of random variables $Y_N:\Omega\to\R$ and $Z_N:\Omega\to\R$, $Y_N = O_P(Z_N)$ if $Y_N/Z_N$ is a stochastically bounded sequence, and $Y_N = o_P(Z_N)$ if $Y_N/Z_N$ converges in probability to zero. 
$Y_N \stackrel{P}{\to} 0$ implies $Y_N = o_P(1)$.

Let $(X,\norm_X)$ be a normed linear space with topological dual space $X^*$. 
A ball in $X$ centered at zero with radius $r$ is denoted $B_X(r)$. 
$\mathcal{N}(S,\|\cdot\|_X,\epsilon)$ is the number of balls needed to cover a subset $S$ using $\|\cdot\|_X$-induced metric balls of radius $\epsilon$ centered around elements of $S$. 
For two normed linear spaces $X$ and $Y$ and a continuous linear operator $A:X\to Y$, $\|A\|_{\rm op}$ is its operator norm and $A^*$ is its adjoint. 

Let $\Z\subset\Rd$ be an open bounded connected set (i.e., domain) with a smooth boundary $\partial\Z$. 
For each $p\in[1,\infty]$, $L^p(\Z)$ is the space of Lebesgue-integrable real-valued functions.
$C(\Z)$ is the Banach space of continuous functions on $\Z$ equipped with a uniform norm $\norm_{\infty}$.
$C^r(\Z)$ for each $r > 0$ denotes the space of $\lfloor r\rfloor$-times continuously differentiable functions with partial derivatives of order $\lfloor r\rfloor$ that are H{\"o}lder continuous with an exponent $r - \lfloor r\rfloor$.
$C_0^\infty(\Z)$ is the space of all infinitely differentiable functions with compact support contained in $\Z$.
$H^m(\Z)$ for each $m > 0$ is the Sobolev space of $\R$-valued functions on $\Z$ for which all $m$-th order weak derivatives exist and belong to $\ltwo(\Z)$.
Sobolev space $H^s(\Z)$ with a fractional order $s\in\R$ is defined based on \cite[Definition~3.22,~3.37]{haroske2007distributions}.
$H_0^s(\Z)$ for each $s \geq 0$ is the closure of $C_0^\infty(\Z)$ in $H^s(\Z)$, with a dual space $H^{-s}(\Z) = (H^s_0(\Z))^*$.

\subsection{Observation model}\label{sec:model}

Let $\mathcal{Z}, \mathcal{X}\subset \mathbb{R}^d$ be bounded domains with smooth boundaries. 
We assume the data arise from the nonparametric regression model
\begin{equation}\label{eq:obs_model}
    Y_i = G(\theta)(X_i) + \varepsilon_i, \quad i=1,\ldots,N,
\end{equation}
where the observation errors $\varepsilon_i$ are i.i.d.\ standard normal random variables.  
The unknown parameter $\theta = \theta_0$ is assumed to belong to a separable linear subspace $\Theta$ of $C(\mathcal{Z})$, the space of continuous functions on $\mathcal{Z}$.  
The joint law of the sample $(Y_i,X_i)_{i=1}^N$ is denoted by $P_{\theta_0}$.
We adopt a random design in which the covariates $X_i$ are i.i.d. uniformly distributed over $\mathcal{X}$, a bounded domain with a smooth boundary $\partial\mathcal{X}$.  
Random design, first proposed for numerical analysis by \cite{diaconis1988bayesian}, provides a probabilistically robust alternative to fixed design and has found wide application.  
Although our proofs are carried out under random design, the main results will remain valid for suitably chosen fixed designs with good interpolation properties, based on arguments used by \cite{yan2024bayesian}.

A natural observation space is $L^2(\mathcal{X})$, the space of functions on $\mathcal{X}$ that are Lebesgue square-integrable.
The operator $G\colon \Theta \to L^2(\mathcal{X})$ encodes the forward mathematical model.  
In this work we focus on nonlinear forward maps that arise, for example, when the coefficient of an elliptic PDE is mapped to the solution of that PDE. Gaussian process priors are a natural modeling choice for $\Theta$ provided their smoothness is compatible with the Sobolev or H\"older scales that describe the regularity of $\theta_0$.  
Throughout, we assume that $G$ is injective, ensuring that $\theta_0$ can in principle be recovered from noiseless observations $G(\theta_0)$.

Assuming the observation model \eqref{eq:obs_model}, we consider general rescaled, $N$-dependent \emph{prior distributions} induced by construction $\theta = \tau_N\theta',\; \theta'\sim \Pi_0$,
following \cite{giordano2020consistency,monard2021statistical}.
The choice of the ``base measure'' $\Pi_0$ and the precise conditions on the scaling sequence $\tau_N$ are given in Section \ref{sec:prior}.
The Bayes posterior distribution corresponding to an $N$-dependent prior $\nprior$ is given by the formula
\begin{equation} \label{eq:posterior}
    \nprior(\theta\in A|X_1,\ldots,X_N,Y_1,\ldots,Y_N) := \frac{\int_A e^{-\loss(\theta)}~d\nprior(\theta)}{\int_\Theta e^{-\loss(\theta)}~d\nprior(\theta)},
\end{equation}
for every $A\subset\Theta$ that is measurable with respect to the canonical Borel $\sigma$-algebra of a separable linear subspace $\Theta\subset\ltwo(\Z)$.
$\loss$ defines the ``empirical loss'' term
\[
\loss(\theta) = \sum_{i = 1}^{N} \ell_i(\theta),\;
\ell_i(\theta) = \frac{1}{2}|Y_{i} - G(\theta)(X_i)|^{2}.
\]
From hereon, we will use the abbreviation $\data = (Y_1,\ldots,Y_N,X_1,\ldots,X_N)$. 

Throughout our study, it will be assumed that the posterior \emph{contracts around $\theta_0$} under $\mathcal{D}^\infty\sim P_{\theta_0}^{\infty}$ (cf. Assumption \ref{assmp:posterior_contraction} for a precise definition).
This is an intuitive consistency demand: all the mass of the posterior distribution should eventually collapse (with high probability) to the true value of $\theta_0$.
It is thus natural to also study a linearization of the posterior around $\theta_0$.
Simplifying \cite[Definition~3.1.2]{nickl2023bayesian}, we assume that $G$ is Fr{\'e}chet differentiable as a map from $\Hs^p\subset C(\Z)$, viewed as a subset with inherited norm $\|\cdot\|_\infty$, to $\ltwo(\X)$ (cf. 
Assumption \ref{assmp:derivative} for a precise definition). Denote by $\mathbb{I}_\theta$ its derivative at $\theta\in\Hs^p$.
Assuming $\theta_0\in\Hs^p(\Z)$, we can define a local approximation of $\ell_N$ as
\begin{equation} \label{eq:lan_lik}
    \llan(\theta - \theta_0) := -\sum_{i=1}^{N}\epsilon_i\opinfo(\theta-\theta_0)(X_i) + \frac{N\|\opinfo(\theta-\theta_0)\|_{\ltwo(\X)}^2}{2},
\end{equation}
and the ``remainder process'' as 
\begin{equation} \label{eq:rem_lik}
    \stoch(\theta,\theta_0) := \loss(\theta) - \loss(\theta_0) - \ell^{LAN}_N(\theta-\theta_0).
\end{equation}
The ``linearized posterior'' replaces the exact likelihood in \eqref{eq:posterior} with its local approximation $\ell^{LAN}_N$:
\begin{equation} \label{eq:posterior_lan}
    \Pi_N^{LAN}(\theta\in A|\data) := \frac{\int_A e^{-\ell^{LAN}_N(\theta)}~d\Pi_N(\theta)}{\int e^{-\ell^{LAN}_N(\theta)}~d\Pi_N(\theta)}.
\end{equation}

\subsection{Prior Distribution} \label{sec:prior}

We restrict the scope of priors considered to Gaussian measures $\Pi_0$ supported on $\ltwo(\Z)$. A good reference for Gaussian measures on Banach spaces is \cite[Chapter~2]{bogachev1998gaussian}. A more condensed account is given in \cite[Setion~6]{stuart2010inverse}.
Let, therefore, $\Pi_0$ be a centered Gaussian distribution on $\ltwo(\Z)$. By definition, a random variable $\theta'\sim\Pi_0$ if $\langle \theta', h\rangle_{\ltwo(\Z)}$ has a Gaussian distribution on $\R$.
By ``centered,'' we mean the mean of every such linear continuous functional is zero.
Centered Gaussian distribution on a separable Hilbert space correspond bijectively to covariance operators \cite[p.2095]{gugushvili2020bayesian}.
The covariance operator $\mathcal{C}$ of $\Pi_0$ can be characterized by the variance relation: $\E\langle \theta',h\rangle_{\ltwo(\Z)}^2 = \langle\mathcal{C}h,h\rangle_{\ltwo(\Z)}$ for every $h\in\ltwo(\Z)$. 
This operator is self-adjoint, non-negative, and trace-class; in particular, it has an explicit eigendecomposition with discrete spectrum.
We assume $\mathcal{C}$ is injective with an inverse $\Lambda$, an $\ltwo(\Z)$-valued unbounded linear operator whose domain $D(\Lambda)$ is dense in $\ltwo(\Z)$. We also assume there exists some $c>0$ such that for every $x\in D(\Lambda)$, we have the numerical range condition
$\|\Lambda^{1/2} x\|_{\ltwo(\Z)} \geq c\|x\|_{\ltwo(\Z)}$.
Then, using the results from \cite[Section~8.4]{engl1996regularization}, we can show that
    $\Hs^\infty := \bigcap_{p\in\R} D(\Lambda^{p})$
is dense in $\ltwo(\Z)$. 
We define the associated \emph{Hilbert scales} $(\Hs^s)_{s\in\R}$ as the completion of $\Hs^\infty$ with respect to the norm $\norm_{\Hs^s}$ induced by the following (well-defined) inner product
$\langle f, g \rangle_{\Hs^s} := \langle \Lambda^{s/2}f, \Lambda^{s/2}g \rangle_{\ltwo(\Z)}$.
Some well-known results about the Hilbert scales are summarized in Appendix \ref{appendix:scales}.

As stated, $\Hs^0 = \ltwo(\Z)$ is assumed to be the support of $\Pi_0$.
In many practical examples, there exists an $s\in (0,1)$ such that $\Pi_0$ has full measure on $\Hs^s$, a Hilbert space dense in $\ltwo(\Z)$. We call any such $\Hs^s$ a Hilbert support of $\Pi_0$, following \cite[Chapter~3]{bogachev1998gaussian}.
This provides a natural choice of a space of possible parameters for our problem; in our previous notation, $\Theta = \Hs^s$.
The \emph{Cameron-Martin space} of $\Pi_0$ for a centered Gaussian measure on $\ltwo(\Z)$ has a simple characterization of being $\Hs^1 = \Lambda^{-1/2}(\ltwo(\Z))$, which by Theorem \ref{thm:hilbert_scales} is equal to $ \Lambda^{-1}(\Hs^{-1}(\Z))$;
it is a separable Hilbert space equipped with the inner product $\langle\cdot,\cdot\rangle_{\Hs^1}$.
Definitions in \cite[Section~2.3--4]{bogachev1998gaussian} show that such a characterization is equivalent to the usual definitions of the Cameron-Martin space.

\begin{remark}
    We make two remarks on our choice of notation.
    First, the Cameron-Martin space of a Gaussian measure is often referred to as the ``reproducing kernel Hilbert space (RKHS)'' in nonparametric statistics literature \citep{ghosal2017fundamentals,gine2021mathematical}. 
    We use a different terminology to maintain consistency with \cite{bogachev1998gaussian}.
    Second, as mentioned by \cite{gugushvili2020bayesian}, it seems natural to force the covariance operator $\mathcal{C}$ correspond to $L^{-2\alpha}$ for some unbounded operator $L$, where $\alpha>0$ is our desired ``smoothness parameter.'' They also use the notation ``$(\Hs^s)_{s\in\R}$'' to indicate Hilbert scales generated by $L$, so that the Cameron-Martin space of the resulting prior must be written as ``$\Hs^\alpha$.''
    Such a notation suits certain covariance choices, such as spectral powers of the negative Dirichlet Laplacian, 
    but is not satisfactory when one wants to decouple the smoothness of $\theta_0$ (measured in the Sobolev sense) from the desired boundary conditions.
    For example, in Section \ref{sec:darcy_flow}, we will consider the prior of \cite{giordano2020consistency} and \cite{monard2021statistical} obtained by smoothly truncating sample paths of a centered stationary Whittle-Mat{\'e}rn Gaussian process to a compact subset of a bounded domain $\Z$. 
    The construction ensures that the Cameron-Martin space of the resulting prior is continuously embedded in $H_0^\alpha(\Z)$, to which $\theta_0$ is assumed to belong.
    However, the boundary conditions are ensured through multiplying sample paths with a smooth cutoff function, so it is awkward to write the resulting covariance operator, implicitly defined, as an $\alpha$-th spectral power of some other operator.
\end{remark}

\section{General Results for Gaussian Priors}
\label{sec:main}

Let $\Psi:\ltwo(\Z)\to\R$ be a non-zero continuous linear functional with Riesz representer $\psi\in\Hs^q,\; q \geq 0$. 
The case $q = 0$ is when $\psi\in\ltwo(\Z)$. 
The goal of this section is to derive results on the contraction and frequentist coverage properties of the induced posterior on the linear functional $\Psi(\theta) = \langle\psi,\theta\rangle_{\ltwo(\Z)}$.
We say that the posterior \emph{contracts around $\Psi(\theta_0)$ at a rate $s_N$} if there exists $b>0$ such that for any sequence $M_N\to\infty$, we have $\ptheta(\Pi_N(|\Psi(\theta-\theta_0)|>M_Ns_N|\data)\geq e^{-bN\epsilon_N^2})\to 0$.
We define a Bayesian credible interval corresponding to the posterior distribution \eqref{eq:posterior} as
    \begin{equation} \label{eq:cred_set}
        I_N := \{r : |r-\langle\psi,\bartheta\rangle_{\ltwo}| \leq Q_{N,\gamma}\},
    \end{equation}
where $Q_{N,\gamma}$ is a $(1-\gamma)$-th quantile satisfying  $\nprior\circ\Psi^{-1}(I_N \mid \data) = 1-\gamma$ and
$\bartheta$ is the mean of the random variable $\theta$ with respect to the posterior distribution \eqref{eq:posterior}, which makes sense as a Bochner integral.
An often small, user-chosen level $\gamma$ specifies the nominal \emph{credible level} of this interval; the hope is that $I_N$ has asymptotically the same \emph{coverage} as $1-\gamma$. Precisely, coverage is defined as $\ptheta(\theta_0\in I_N)$.
We say that the coverage is \emph{asymptotically correct} if
\[\liminf_{N\to\infty} \ptheta(\theta_0\in I_N) \geq 1-\gamma,\] 
and \emph{asymptotically exact} if 
\[\lim_{N\to\infty} \ptheta(\theta_0\in I_N) = 1-\gamma.\]

\subsection{Main Results} 
\label{sec:thm_general}

This section contains the main result of this work on the asymptotic frequentist coverage of Bayesian credible intervals.
We make the following assumptions about the prior, operator $G$, and the posterior.

\begin{assumption}\label{assmp:truth_prior}
    For some $p\in(0,1)$, the Gaussian prior distribution $\Pi_0$ defined by a sequence $\tau_N$ and a centered Gaussian measure $\Pi_0$ places full measure on $\Hs^p$ which is continuously embedded into $C(\Z)$.
    The true parameter $\theta_0$ (i.e., $(Y_i,X_i) \stackrel{iid}{\sim} P_{\theta_0}$) belongs to $B_{\Hs^\beta}(M)$ for some $M>0$ and $\beta\geq 1$.
\end{assumption}

\begin{assumption}
    \label{assmp:unif_bd_continuity}
    For every $\theta\in B_{\Hs^p}(M)$, a Hilbert support from Assumption \ref{assmp:truth_prior}, 
    there exists $U = U(M) > 0$ such that $\|G(\theta)\|_\infty \leq U$. Furthermore, for every $\theta_1,\theta_2\in B_{\Hs^p}(M)$,
    there exists $L = L(M) > 0$ such that 
    \[\|G(\theta_1) - G(\theta_2)\|_{\ltwo(\X)} \leq L\|\theta_1 - \theta_2\|_{\ltwo(\Z)}.\]
\end{assumption}

\begin{assumption} \label{assmp:derivative}
For every $\theta$ belonging to a Hilbert support $\Hs^p$ of $\Pi_0$ from Assumption \ref{assmp:truth_prior},
there is a continuous linear operator $\mathbb{I}_\theta : \Hs^p \to \ltwo(\X)$ such that the remainder term $R_\theta(h) := \|G(\theta + h) - G(\theta) - \mathbb{I}_\theta(h)\|_{\ltwo(\X)}$ defined for each $h\in\Hs^p$ satisfies $\|R_\theta(h)\|_{\ltwo(\X)} = O(\|h\|_{\infty})$ as $\|h\|_{\infty}\to 0$.
Furthermore, $\opinfo$ is in fact a continuous linear operator from $\Hs^p$ to $C(\X)$.
\end{assumption}

\begin{assumption}
    \label{assmp:prior_mass}
    We have $\frac{1}{\sqrt{N}} \ll \tau_N \lesssim 1$ and 
    $\nprior(\theta : \|G(\theta) - G(\theta_0)\|_{\ltwo(\X)} \leq \epsilon_N) \geq e^{-AN \epsilon_N^2}$
    for some $A> 0$ and a sequence $\epsilon_N\to 0$ satisfying $\sqrt\frac{\log(N)}{N} \ll \epsilon_N$. 
\end{assumption}

\begin{assumption} \label{assmp:posterior_contraction}
    There exists a sequence of measurable sets $\Theta_N$
    such that 
    $\sup_{\theta\in\Theta_N}|\Psi(\theta-\theta_0)| = o_N(1)$ and the posterior distribution \eqref{eq:posterior} satisfies, for some $b > A + 2$ for 
    $A$ from Assumption \ref{assmp:prior_mass},
$\ptheta(\nprior(\Theta_N|\data) \geq e^{-bN\epsilon_N^2}) \to 0$.
\end{assumption}

\begin{assumption}
\label{assmp:deriv_stability}
    For every $M > 0$ and every $s\geq 1$, there exists a $\delta = \delta(s) \in (0,1]$ and $c=c(s,M) > 0$ so that we have
    \[\|h\|_{\ltwo(\Z)} \leq c(s,M)\|\opinfo(h)\|_{\ltwo(\X)}^{\delta(s)}\] for every $h\in B_{\Hs^s}(M)$.
\end{assumption}

\begin{remark}
\cite[Conditions~3.1--2,~4--5]{monard2021statistical} are essentially comparable with our Assumptions \ref{assmp:truth_prior}--\ref{assmp:posterior_contraction},
though we simplified \cite[Condition~3.2]{monard2021statistical} to state a (slightly stronger) Assumption \ref{assmp:derivative} for an easier exposition.
Assumption \ref{assmp:deriv_stability} is a novel ``injectivity condition,'' mentioned in the introduction as a tool for verifying a variational inequality. 
It has a similar flavor \cite[Theorem~5]{nickl2022some}, but our assumption only considers local estimates for uniformly bounded $h\in\Hs^1$. 
\end{remark}

Before proceeding any further, we first define several non-stochastic quantities that will be useful in the statement of our main results. 
Their motivation is a basic linearization of the operator $G$ around $\theta_0$: were $G$ truly linear, $\Lambda$ can be suitably chosen such that results of \cite{knapik2011bayesian} imply the frequentist coverage results we want to show. 
Since this is not the case, several quantities will arise from linearizing $G$ whose orders in $N$ need to be estimated.

The bounded linear operator $\opcomb : \ltwo(\Z) \to \ltwo(\X)$ is defined by 
\begin{equation} \label{eq:opcomb}
    \opcomb := \opinfo\Lambda^{-1/2}.
\end{equation}
$\opcomb^*\opcomb$ is a positive, self-adjoint, trace-class linear operator on $\ltwo(\X)$, its trace bounded above by the operator norm of $\opinfo^*\opinfo$ times the trace of $\Lambda^{-1}$.
Next, we define $\psi_{\lambda} := \Lambda^{-1/2}(\psi)\in D(\Lambda^{1/2})$,
based on which we define a sequence of positive real numbers $s_N$, depending on $\psi$:
\begin{equation} \label{eq:asymp_scale}
    s_N = s_N(\tau_N,\psi,\opcomb) = \sqrt{
        \langle \psi, 
        \Lambda^{-1/2}(N\opcomb^*\opcomb + \tau_N^{-2})^{-1}\lpsi 
        \rangle_{\ltwo(\Z)}}.
\end{equation}
This sequence is deterministic. The proofs will clarify that this sequence may be thought of as the asymptotic variance of the posterior distribution of $\Psi(\theta)$; importantly, it may have a different order in $N$ depending on $\psi$.
Furthermore, the proofs will rely on a construction of ``perturbative'' sequence
\begin{equation} \label{eq:perturbation}
    \barpsi = \overline{\psi}_N(\tau_N,\opcomb,\psi) = 
    \Lambda^{-1/2}(N\opcomb^*\opcomb + \tau_N^{-2})^{-1}\lpsi.
\end{equation}

The proof of the following claim is postponed to the Supplement \cite{supplement}.
\begin{lemma} \label{lem:perturbation_lemma}
    For every finite $N$, the operator $N\opcomb^*\opcomb + \tau_N^{-2}$ has an injective bounded linear inverse on $\ltwo(\Z)$
    that maps $D(\Lambda^{s})$ to an image contained in $D(\Lambda^{s})$ for every $s\geq 0$.
    In particular, $\barpsi$ belongs to $\Hs^2=D(\Lambda)$.
    Furthermore,
    $u = (N\opinfo^*\opinfo + \tau_N^{-2}\Lambda)\Lambda^{-1/2}(N\opcomb^*\opcomb+\tau_N^{-2})^{-1}\Lambda^{-1/2}u$ for every $u\in\ltwo(\Z)$.
\end{lemma}
This lemma sheds some light on the importance of the sequence of vectors $\barpsi$. 
The proof of BvM from \cite{monard2021statistical} relies on constructing a sequence of perturbations to a posterior Laplace transform, which we also follow (see the proof of Lemma 5.1 \cite{supplement}). 
Under \cite[Condition~3.3]{monard2021statistical}, such a sequence is obtained by solving an ``information equation'' $\opinfo^*\opinfo x = \psi$ in $x\in \ltwo(\Z)$.
Furthermore, \cite[Condition~3.6]{monard2021statistical} assumes that the solution belongs to $\Hs^1$ (written as ``$\mathcal{H}_N$'' therein).
These conditions not only stipulate that a solution exists, but also it has sufficient regularity \cite[p.81]{nickl2023bayesian}.
The equation is not solvable if $\psi\notin R(\opinfo^*)$, a fact related to the failure of BvM and lack of efficient estimators.
Our sequence $\barpsi$, while playing an analogous role in the proof, may be viewed as solutions of $N$-dependent equations $(N\opinfo^*\opinfo + \tau_N^{-2}\Lambda)x = \psi$, which, based on the above lemma, is always well-defined and automatically has $\Hs^2$-regularity.

We now state the main result of this work: under some, for the moment opaque, conditions, the frequentist coverage of $I_N$ \eqref{eq:cred_set} is indeed asymptotically correct.
As a side result, it also yields an upper bound for the rate of contraction for the posterior pushforward under the map $\Psi$ around $\Psi(\theta_0)$.

\begin{theorem} \label{thm:freq_coverage}
    Suppose Assumptions \ref{assmp:truth_prior}-\ref{assmp:posterior_contraction} hold and $\psi\neq 0$.
    Suppose, furthermore, that for $\mathcal{F}_N$ defined in \eqref{eq:rem_lik} we have
    \begin{equation} \label{eq:stoch_unif_control}
        \sup_{\theta\in\Theta_N}|\mathcal{F}_N(\theta,\theta_0)| = o_{\ptheta}(N\epsilon_N^2)
    \end{equation}
    and that the ``linearized posterior'' from \eqref{eq:posterior_lan} satisfies, for each $t\in\R$ and for large enough $N=N(t)$,
    \begin{equation} \label{eq:lan_ratio_control}
        \ptheta\left(
        \left| \frac{\nlan\{ \Theta_N - ts_N^{-1}\barpsi| \data \}}{\nlan\{\Theta_N | \data \}} - 1\right| \geq e^{-cN\epsilon_N^2}\right) \to 0,
    \end{equation}
    for some $c>0$. If, additionally,
    \begin{equation} \label{eq:bias}
    b_N := b_N(\tau_N,\theta_0,\Lambda,\opcomb,\psi) = \tau_N^{-2}\langle \theta_0, \Lambda\barpsi \rangle_{\ltwo(\Z)}
    \end{equation}
    satisfies $|b_N|/s_N\ll 1$, then:
    \begin{enumerate}[label=(\roman*)]
        \item ${\liminf}_{N\to\infty} \ptheta
        (\Psi(\theta_0) \in I_N) \geq 1-\gamma$ if $\psi\notin R(\opinfo^*)$;
        \item ${\lim}_{N\to\infty} \ptheta
        (\Psi(\theta_0) \in I_N ) = 1-\gamma$ 
        if $\psi\in R(\opinfo^*)$.
    \end{enumerate}
    Moreover, for an arbitrarily slowly growing sequence $M_N\to \infty$ and any $\epsilon > 0$, we have
    \begin{equation}
        \label{eq:fnl_contraction}
        \nprior\{ \Psi(\theta - \theta_0) \geq M_N s_N | \data \} = o_{\ptheta}(1).
    \end{equation}
\end{theorem}

\begin{remark}
\cite[Theorem~3.8]{monard2021statistical} is a previous BvM for PDE parameter identification problems.
\cite[Condition~3.6]{monard2021statistical} was imposed to prove a uniform $o_{\ptheta}(1)$ bound on a similar ``remainder empirical process'' as $\mathcal{F}_N(\theta,\theta_0)$.
Condition \eqref{eq:stoch_unif_control} is more relaxed and allows a non-vanishing, even slowly divergent, empirical process. 
The additional cost we pay is condition \eqref{eq:lan_ratio_control}, which relates to the concentration properties of a ``linearized posterior'' \eqref{eq:posterior_lan} and the size of shift sequence $s_N^{-1}\barpsi$ with respect to the high posterior probability sets $\Theta_N$.
The proof of Theorem \ref{thm:darcy_coverage} shows that, as long as $\nlan$ also concentrates on $\Theta_N$ with the same rate sequence $\epsilon_N$ as in Assumption \ref{assmp:posterior_contraction},
condition \eqref{eq:lan_ratio_control} is a reasonable one.
\end{remark}

The theorem reduces the study of the relative rate of decay between two deterministic sequences $b_N$ and $s_N$, which depend on both $\psi$ and $\theta_0$. 
\eqref{eq:fnl_contraction} asserts that the sequence $s_N$ yields a posterior rate of contraction around $\Psi(\theta_0)$.
An interpretation of $b_N$ is that it is the $\theta_0$-dependent ``asymptotic bias'' of the posterior, which may be non-vanishing when $\psi\notin R(\opinfo^*)$. 
If $|b_N|/s_N\ll 1$ as stipulated, then our previous comment that $s_N^2$ may be interpreted as the asymptotic posterior variance is justifiable.

The following results provide a way of estimating these sequences, but the sufficient conditions are quite technical.
Importantly, we now use Assumption \ref{assmp:deriv_stability}, and our results show that one can ``trade off'' smoothness $\beta$ of $\theta_0$ with the smoothness of $\psi$ to ensure $|b_N|/s_N \ll 1$.

\begin{theorem} \label{thm:order_estimates}
    Under Assumptions \ref{assmp:truth_prior}-\ref{assmp:posterior_contraction},
    let $\beta\geq 1$, 
    $\psi \in \Hs^q,\; q\geq 0$ with $\psi\notin {\rm ker}(\Lambda^{-1/2})$.
    $\frac{|b_N|}{s_N} \ll 1$ if
    \begin{enumerate} [label=(\roman*)]
        \item  \label{condition1:tau1}
        Assumption \ref{assmp:deriv_stability} holds, $\tau_N\asymp 1$ and
        $\delta(2 + q)(\frac{\beta}{2}\wedge 1) > \frac{1}{2}$;
        or,
        \item \label{condition2:tau_decay}
        Assumption \ref{assmp:deriv_stability} holds,
        $\tau_N \gg N^{-\frac{2\tilde{\delta} - 1}{4( 1 - \tilde{\delta})}}$ for
        $\tilde{\delta} = \delta(2+q)(\frac{\beta}{2}\wedge 1)$;
        or,
        \item \label{condition3:bvm_rate}
        $\psi\in R(\opinfo^*\opinfo\Lambda^{-\frac{1}{2}}),\; \beta = 1$ and $\tau_N \gg N^{-\frac{3}{8}}$.
    \end{enumerate}
    We also have:
    \begin{enumerate}[label=(\roman*)]
        \item $s_N\lesssim \tau_N^{1-\delta}N^{-\delta/2}$, for every $\delta\in (0,\delta(1))$, if $\psi\notin R(\opinfo^*)$ and Assumption \ref{assmp:deriv_stability} holds; 
        \item $s_N\asymp N^{1/2}$ and $\sqrt{N}s_N\to\|(\opinfo^*\opinfo)^{-1/2}\psi\|_{\ltwo(\Z)}$, if $\psi\in R(\opinfo^*)$. 
    \end{enumerate}
\end{theorem}

\begin{remark}
    That we split the upper bounds for $s_N$ depending on whether $\psi\in R(\opinfo^*)$ is related to our discussion in the introduction regarding the necessary condition for the existence of efficient estimators.
    In the context of PDE-based inverse problems and model \eqref{eq:obs_model}, the condition boils down to $\psi\in R(\opinfo^*)$, as shown by \cite{nickl2022some}.
    Inspection of the proof of Theorems \ref{thm:freq_coverage} and Theorem \ref{thm:order_estimates} shows that with \eqref{eq:stoch_unif_control}, \eqref{eq:lan_ratio_control} and $|b_N|/s_N\ll 1$ being true, 
    the condition $\psi\in R(\opinfo^*)$ yields a BvM, in the sense that a re-centered Bayes posterior of a linear functional $\Psi(\theta)$, re-scaled by $\sqrt{N}$, weakly converges in probability to $Normal(0,
    \|(\opinfo^*\opinfo)^{-\frac{1}{2}}\psi\|_{\ltwo(\Z)}^2)$ (cf. Lemma \ref{lem:gaussian_limit} for a more precise statement).
\end{remark}

\subsection{Discussion}

Theorem \ref{thm:order_estimates} uses Assumption \ref{assmp:deriv_stability} and further conditions on $\beta$ to derive an asymptotically correct coverage.
If $\tau_N\asymp 1$, then since $\delta(s)\leq 1$ for any $s\geq 1$ by assumption, condition \ref{condition1:tau1} is met only if $\beta > 1$.
To see that $\beta > 1$ is a reasonable requirement, let us consider an ``$\alpha$-smoothing prior'' that has been studied multiple times in the literature.
Let $\Lambda = L^{2\alpha}$ for $L$ the negative Dirichlet Laplacian on $\Z\subset\R$ (so $d=1$ and $\tau_N\asymp 1$),
so that $\Hs^1 = D(\Lambda^{1/2})$ can be continuously embedded into $H^\alpha_0(\Z)$, with $\|\cdot\|_{\Hs^1}$ equivalent to $\|\cdot\|_{H^\alpha_0}$ on $\Hs^1$.
If we now additionally assume that $G$ is linear and $\theta_0\in D(L^{\alpha + \frac{1}{2}})$, then \cite[Theorem~6.1]{gugushvili2020bayesian} shows Assumption \ref{assmp:posterior_contraction} is met with $\epsilon_N$ minimax optimal over a ball in Sobolev space $H^{\alpha + \frac{1}{2}}(\Z)$.
In our notation, then, $\theta_0\in \Hs^{1+\frac{1}{2\alpha}}$, so $\beta > 1$.
Whether condition \ref{condition1:tau1} is now met depends on the smoothness of $\psi$ and the obtainable dependence of $\delta(s)$ on $s\geq 1$ in Assumption \ref{assmp:deriv_stability}.

Conditions \ref{condition2:tau_decay} and \ref{condition3:bvm_rate} provide alternative conditions under which we can alleviate the demand $\beta > 1$. 
The second condition, though technical, will be useful when we analyze a Gaussian prior for Darcy flow with $\tau_N\to 0$ in Section \ref{sec:darcy_flow}.
The third condition states that $\psi\in R(\opinfo^*\opinfo\Lambda^{-\frac{1}{2}})$;
as previously mentioned,     \cite[Conditions~3.3,~3.6]{monard2021statistical} stipulate this exact condition, before proving BvM (Theorems 3.7--8).
Let us compare, then, how our other sufficient conditions compare with those of \cite{monard2021statistical}.
\cite[Conditions~3.1--2,~4--5]{monard2021statistical} are comparable with our Assumptions \ref{assmp:truth_prior}--\ref{assmp:posterior_contraction},
and Assumption \ref{assmp:deriv_stability} is not necessary when $\psi\in R(\opinfo^*\opinfo\Lambda^{-1/2})$.
Furthermore, their re-scaling sequence $\tau_N$, given in display (10) but with a different notation ``$\alpha$'' for $\beta$, the known smoothness of true $\theta_0$, does satisfy $\tau_N\gg N^{-\frac{3}{8}}$. 
The remaining conditions of \cite[Condition~3.6]{monard2021statistical} yield their Lemma 4.3, which derives a uniform $o_{\ptheta}(1)$ bound on a comparable empirical process as \eqref{eq:rem_lik}. 
On the other hand, we assumed a more relaxed upper bound \eqref{eq:stoch_unif_control} to be the case. To derive such a bound, we go through similar entropy number calculations as the \cite[Condition~3.6]{monard2021statistical} in the proof of Theorem \ref{thm:darcy_coverage} \cite{supplement}, but since we only demand \eqref{eq:stoch_unif_control}, we can make less strong assumptions about $\beta$.
Finally, while \eqref{eq:lan_ratio_control} was introduced by us, it is proven to be satisfied for Darcy flow in \cite{supplement}, and by following similar arguments, one can do so for the Schr{\"o}dinger equation as well.
Furthermore, since only $\psi\in R(\opinfo^*)$ is needed to ensure $\sqrt{N}s_N\to\|(\opinfo^*\opinfo)^{-1/2}\psi\|_{\ltwo(\Z)}$, Theorem \ref{thm:order_estimates} \ref{condition1:tau1}--\ref{condition2:tau_decay} can be used with Theorem \ref{thm:freq_coverage} to still obtain a BvM when $\psi\in R(\opinfo^*)$ but not $R(\opinfo^*\opinfo\Lambda^{-1/2})$. 
Therefore, our results provide both an improved version of an existing BvM and conditions on the smoothness of $\theta_0$ to obtain BvM when $\psi$ has less regularity.

\section{Application: Darcy Flow Problem} 
\label{sec:darcy_flow}

This section summarizes an application of Theorems \ref{thm:freq_coverage} and \ref{thm:order_estimates} to the case of identifying the unknown parameter in PDE \eqref{eq:darcy}.
As mentioned in the introduction, 
\cite{nickl2022some} showed for this problem how $\psi\in C_0^\infty(\Z)$ need not belong to $R(\opinfo^*)$, even in the simplest Darcy flow model where a solution to the PDE can be derived by hand, which implied the impossibility of BvM.
This is disappointing, as one is left without any knowledge about the asymptotic coverage of posterior credible intervals.
We now show that under natural smoothness conditions on $\theta_0$ and a target linear functional, one can recover correct (albeit not exact) asymptotic coverage of Bayesian credible intervals.

\subsection{Observation Model}
\label{sec:forward_model_darcy}

Consider a bounded domain $\Z=\X \subset \R^d$ and an elliptic PDE with inhomogeneous boundary condition:
\begin{equation}\label{eq:diffusion_eq}
    \begin{cases}
        \nabla\cdot(a\nabla u) = f &\text{on }\Z,\\
        u = u_D &\text{on }\partial\Z.
    \end{cases}
\end{equation}
To highlight only the relevant statistical aspect of the inverse problem, we will assume the source term $f$ and the boundary condition $u_D$ belong to $C^\infty(\overline{\Z})$.
For a fixed, unknown conductivity field $a_0$, we denote the corresponding PDE solution $u = u_{a_0}$.
Suppose $\alpha > 1 + d/2$ and define $n=n(x)$ as the outward-pointing normal vector at $x\in\partial\Z$, with the boundary derivative $\frac{\partial a}{\partial n}$ understood in the trace sense if necessary.
If $a_0\in \Gamma$, defined by:
\begin{equation}
    \label{eq:constraint_set}
    \Gamma_{\alpha, K_{min}} :=
     \bigg\{a \in H^\alpha(\Z) : 
    a > K_{min}\text{ on }\Z,\;a=1\text{ on }\partial\Z,\;
   \frac{\partial^j a}{\partial n^j} = 0 \text{ on }\partial\Z\text{ for } j=1,\ldots,\alpha-1 \bigg\},
\end{equation}
then the differential operator $u\mapsto \nabla\cdot(a\nabla u)$ is strongly elliptic on $\Z$. We assume a fixed lower bound $K_{min}\in (0,1)$ is also known.
Given $a_0 \in H^\alpha(\Z)$, \cite[Theorems~8.3--13]{gilbarg1977elliptic} imply that the PDE \eqref{eq:diffusion_eq} has a unique strong solution $u = u_{a_0}$ belonging to $H^{\alpha+1}(\Z)$.
The parameter space $\Gamma$, however, is not a linear space, so we consider a bijective reparameterization of $a$ into $\theta$ belonging to some linear space.

\begin{definition}
    \label{def:link_fn}
    \cite[Definition~7]{nickl2020convergence}
    We say $\phi$ is a \emph{link function} if it is a smooth, strictly increasing bijection from $\R$ to $(K_{min},\infty)$, 
    satisfying $\phi(0) = 1$ and $\phi' > 0$.
    We say a link function $\phi$ is \emph{regular} if $
    \sup_{x\in\R}|\phi^{(k)}(x)| < \infty$ for every integer $k\geq 1$.
\end{definition}

\cite{nickl2020convergence} showed that a link function $\phi$ is a smooth bijection between the linear space $H^\alpha_0(\Z)$ and $\Gamma$.
Regularity of $\phi$ was used in the authors' proof of convergence properties of a penalized M-estimator.
An easy-to-compute link function considered by \cite{stuart2010inverse} of the form $\phi(x) = e^x + K_{min}$, on the other hand, is not regular.
Whether $\phi$ is regular or not, we can smoothly parameterize $a = \phi(\theta)$ and define the forward map as $G(\theta) = u_{\phi\circ\theta}$,
based on the uniqueness of the solution.

The goal of the inverse problem is to estimate the reparameterization $\theta_0$ of $a_0$ based on noisy observations given by the model \eqref{eq:obs_model}.
The following assumption was introduced by \cite{richter1981inverse} and has been used by previous analyses of the problem from a Bayesian angle by \cite{giordano2020consistency}.
\cite{richter1981inverse} already highlighted the fact that without such a condition, estimating $\theta_0$ becomes a very ill-conditioned problem irrespective of $\alpha$. 
Readers interested in this ``low-regularity'' problem may refer to the more recent results \cite{bonito2017diffusion}.
\begin{assumption}
    \label{assmp:darcy_flow}
    $\theta_0$ is compactly supported on some $K\subsetneq \Z$. 
    There exists a partition of $\Z$ into subsets $\Z_1$ and $\Z_2$, each possibly disconnected, such that
    $G(\theta_0)$ satisfies
$\inf_{x\in \Z_1}\|\nabla G(\theta_0)(x)\|_{\R^d}\geq k_1,\;
    \inf_{x\in\Z_2}\Delta G(\theta_0)(x) \geq k_2$ for positive constants $k_1,\; k_2$.
\end{assumption}
Note that the above assumption is implied by a positive source term $f$ strictly bounded away from 0.
For us, the above assumption has a further benefit that Assumption \ref{assmp:deriv_stability} now becomes verifiable for Darcy flow; see the proof of Lemma D.1 \cite{supplement}.

\subsection{Gaussian Process Prior Distribution}
\label{sec:prior_giordano_nickl}

We briefly review the properties of a previously proposed prior construction from \cite[Example~25]{giordano2020consistency}.
In our notation, a sequence of $N$-dependent priors $\Pi_N$ is constructed by choosing a decaying sequence $\tau_N = (\sqrt{N}\epsilon_N)^{-1},\; \epsilon_N = N^{-\frac{\alpha+1}{2\alpha+2+d}}$,
and a Gaussian ``base measure'' $\Pi_0$ induced by the following procedures.
First, take a centered stationary Whittle--Mat{\'e}rn Gaussian process $\{Z(s)\}_{s\in\mathcal Z}$ on a smooth, bounded domain $\mathcal Z\subset\mathbb R^d$.  
Its covariance kernel is given by
\begin{equation}
    \label{eq:whittle_matern}
    \mathcal{K}(t-s)
   = \int_{\mathbb R^d} e^{-\sqrt{-1}\langle s-t,\xi\rangle}\,\bar{\mu}(d\xi),
   \qquad t,s\in\mathcal Z,
\end{equation}
where the spectral density is
$\bar{\mu}(d\xi) = (1+\|\xi\|_{\mathbb R^d}^2)^{-\alpha}\,d\xi$ for a user-chosen $\alpha > d/2$.
We state the facts given in \cite[Example~25]{giordano2020consistency} about the prior measure induced, denoted by $\tilde{\Pi}_0$:
namely, the Gaussian process sample paths have versions belonging to $C^s$ for any choice of $s<\alpha-d/2$, and the Cameron-Martin space equals the set of restrictions of all functions in $H^\alpha(\Rd)$ to $\Z$, which is a Hilbert space with norm equivalent to $H^\alpha(\Z)$-norm.
Next, to avoid dealing with delicate boundary issues, we assume that $\theta_0$ is supported on a compact subset $K\subsetneq \mathcal Z$ and define $\Pi_0'$ as the pushforward measure $\tilde{\Pi}_0\circ M_\chi^{-1}$ with $M_\chi f = \chi f$,
$\chi$ being a smooth bump function supported on a compact set $K'\supsetneq K$, evaluates to 1 on $K$, and vanishes only at $\partial K'$.
The Cameron-Martin space $\Hs^1$ of $\Pi_0'$ is continuously embedded into $H^\alpha_0(\Z)$ and contains the closed subspace $H_K^\alpha(\mathcal Z)$ of functions in $H^\alpha_0(\Z)$ supported on $K$ \cite[Section~2.2.1]{giordano2020consistency}. 

The implicit definition of a covariance operator $\mathcal{C} : \ltwo(\Z)\to\ltwo(\Z)$ is given by the integral equation $(\mathcal{C}g)(s) = \int_\Z\mathcal{K}(t-s)g(t)~dt$ \cite[Lemma~I.15]{ghosal2017fundamentals}.
Since $\chi$ vanishes outside $K'\subsetneq\Z$, $\mathcal{C}$ has a non-trivial kernel; however, since $\theta\sim \Pi_0'$ is supported on $K'$ with probability 1 anyway, we can view the prior as being defined on a smaller Hilbert space $\ltwo(K')$ and only consider linear functionals in $\ltwo(K')$.

\begin{lemma} \label{lem:higher_cm_embedding}
The inverse $\Lambda$ of the covariance operator $\mathcal{C}$ for $\Pi_0' = \tilde{\Pi}_0\circ M_\chi^{-1}$ is defined on a domain that is dense in $\ltwo(K')$.
It satisfies $\|\Lambda^{1/2}x\|_{\ltwo(\Z)}\geq\|x\|_{\ltwo(\Z)}$ for all $x\in D(\Lambda^{1/2})$.
Furthermore, for every $s>1$, $\Hs^s$ is continuously embedded into $H^{s\alpha}(\Z)$.
\end{lemma}

The suggestive notation for the sequence $\epsilon_N$ was used to indicate that it is a rate of contraction for the posterior distribution around $\theta_0$ (cf. Assumption \ref{assmp:posterior_contraction}).
More precisely, consider $\alpha > 1 + d$. \cite[Theorems~4--5]{giordano2020consistency} prove that for a sequence of sets
\begin{equation}
    \label{eq:lsupset}
    \Theta_N := 
    \{\theta : \|G(\theta) - G(\theta_0)\|_{\ltwo(\Z)}\leq M\epsilon_N\}\cap
    \{\theta : \|\theta\|_{H^p(\Z)} \leq M\},
\end{equation}
where $M>0$ is large enough and the exponent $p$ is arbitrarily chosen in the interval $(1+\frac{d}{2},\alpha-\frac{d}{2})$, 
the posterior satisfies $\ptheta(\nprior(\Theta_N|\data) \geq e^{-bN\epsilon_N^2}) \to 0$ for some $b>0$.
Known ``conditional stability'' results can then be used to derive convergence towards $\theta_0$ in more natural norms and thus verify Assumption \ref{assmp:posterior_contraction}, as we do in the proof of Theorem \ref{thm:darcy_coverage} below.

\subsection{Frequentist Coverage Results}
\label{sec:coverage_darcy}

\begin{figure}[!ht]
    \centering
    \includegraphics[width=0.6\linewidth]{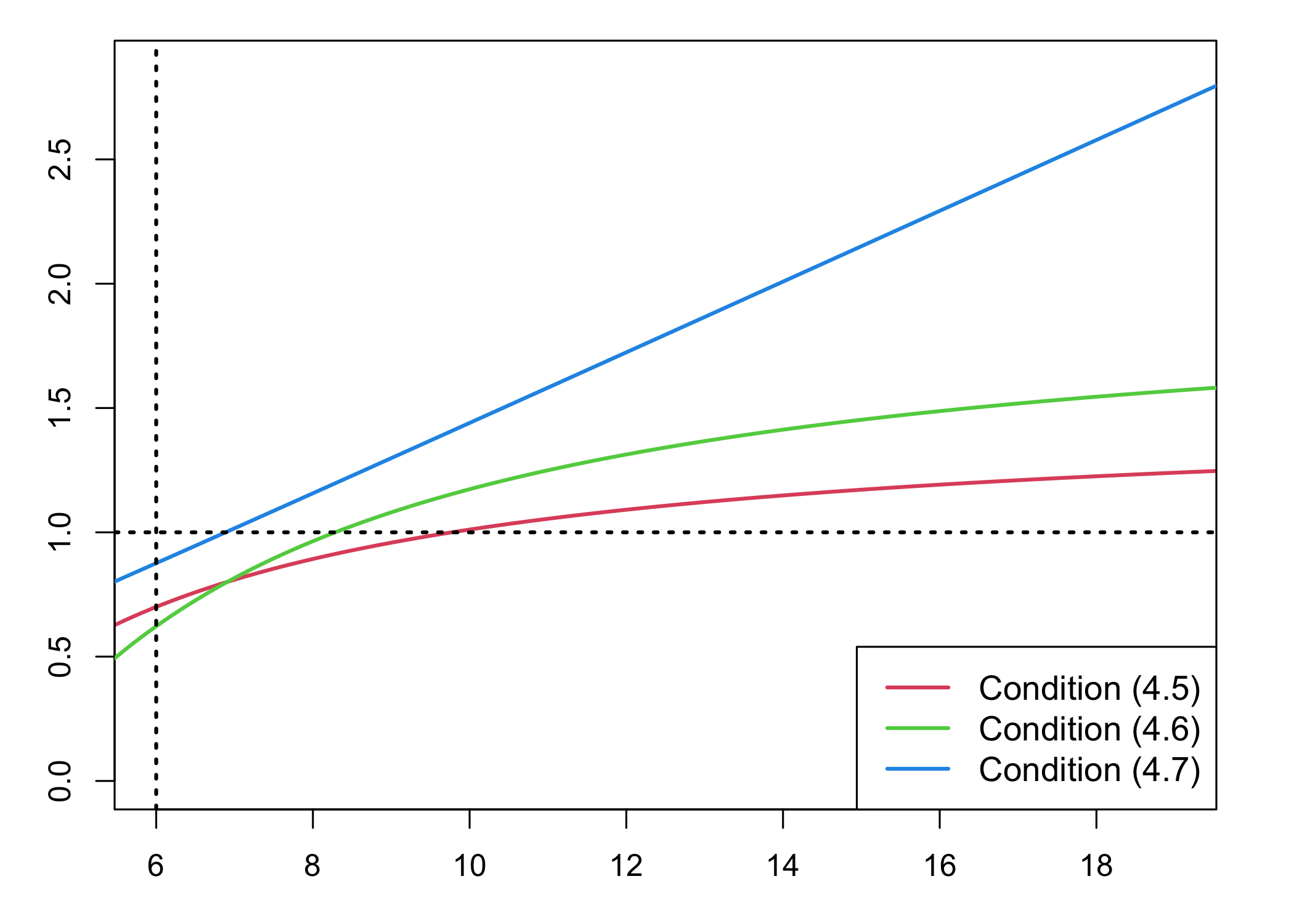}
    \caption{Left-hand sides of conditions \eqref{eq:complex_condn1}, \eqref{eq:complex_condn2} and \eqref{eq:complex_condn3} plotted as functions of $\alpha$ when $d = 2$. 
    The vertical line is at our stipulated lower bound $=3 + \frac{3}{2}d$ on $\alpha$. The horizontal line is at 1.}
    \label{fig:conditions_illustrated}
\end{figure}

Let $\Psi : \theta\mapsto \langle\theta,\psi\rangle_{\ltwo(\Z)}$ be a continuous linear functional on $\ltwo(\Z)$.
Define a quantile-based posterior credible band $I_N$ for $\Psi(\theta)$ as in \eqref{eq:cred_set}.
We now state a frequentist coverage result for this Bayes posterior, using Theorems \ref{thm:freq_coverage} and \ref{thm:order_estimates}. Its proof is given in the Supplement \cite{supplement}.

\begin{theorem}
    \label{thm:darcy_coverage}
    Suppose $\alpha > 1 + d$. 
    Let $\nprior = Law\left(\tau_N\theta'\right),\; \theta'\sim \Pi_0,\;
    \tau_N = N^{-\frac{d}{4\alpha+4+2d}}$ 
    for $\Pi_0$ constructed in Section \ref{sec:prior_giordano_nickl}.
    Assume $\theta_0 \in \Hs^\beta,\; \beta > \alpha$, 
    $(Y_i,X_i) \stackrel{iid}{\sim} P_{\theta_0}$ and the linear functional of interest satisfies $\psi\in \Hs^{q}$ for some $q \geq 0$ and $\psi\neq 0$ inside $K'$.
    Under Assumption \ref{assmp:darcy_flow},
    Assumptions \ref{assmp:truth_prior}--\ref{assmp:deriv_stability}, 
    conditions \eqref{eq:stoch_unif_control} and \eqref{eq:lan_ratio_control}, and $|b_N|/s_N = o_N(1)$ are all satisfied 
    if $\alpha > 3 + \frac{3}{2}d$ and
        \begin{eqnarray}
        \frac{1}{2}\frac{\alpha-1-d}{\alpha+1-d/2}\left(
        2 + \frac{\alpha-1-d/2}{\alpha-d/2}\right) & > & 1,
        \label{eq:complex_condn1}\\
        \frac{(\alpha+1)(\alpha-1-d/2)(2\alpha-2-3d)}{(\alpha+1-d/2)^2(\alpha-d/2)} & > & 
        1,
        \label{eq:complex_condn2}\\
        \left\{
        \frac{2\alpha+2}{2\alpha+2+d} +
        \frac{2d}{\alpha+1-d/2} - \frac{(\alpha-1-d/2)(\alpha-1-d)}{(\alpha+1-d/2)(\alpha-d/2)}
        \right\}^{-1} & > & 1,
        \label{eq:complex_condn3}\\
        \frac{d}{4\alpha+4+2d} < \frac{2\tilde{\delta}-1}{4(1-\tilde{\delta})},\; \tilde{\delta} 
        = \frac{(2+q)\alpha-1}{(2+q)\alpha+1}
        & \times & \left(\frac{\beta}{2\alpha}\wedge 1\right).
        \label{eq:nobias_condn_darcy}
        \end{eqnarray}
\end{theorem}

\begin{remark}
    Conditions \eqref{eq:complex_condn1}-\eqref{eq:complex_condn3} control the ``remainder process'' appearing in  \eqref{eq:stoch_unif_control} and ensure \eqref{eq:lan_ratio_control}.
    These conditions are eventually true for large enough $\alpha$, because the left-hand sides of the conditions are increasing functions of $\alpha$ bounded from below by 1 outside some neighborhood of the origin.
    When $d = 2$, all conditions are met and $\alpha > 3 + \frac{3}{2}d$ for integer order $\alpha\geq 10$;
    when $d = 3$, integer order $\alpha\geq 12$ suffices. See Figure 
    \ref{fig:conditions_illustrated}
    for a graphic illustration of these conditions.
\end{remark}

We have also conducted numerical experiments verifying coverage of Bayes posterior credible intervals for finite $N$,
results of which are included in \cite{supplement}.
We took a toy example from \cite{nickl2022some}, which shows that BvM is not satisfied for a functional with the representer $\psi\in C^\infty_0(\Z)$ and $\theta_0 = 0$.
Note, however, that $\theta_0 = 0$ implies the sequence $b_N$ \eqref{eq:bias} is also 0, so the implication of Theorem \ref{thm:order_estimates} that $|b_N|/s_N\ll 1$ is trivially met. 
The fact hints that, despite the failure of BvM, good coverage should be easily met by the posterior.
Results of our simulated study demonstrate not only that $\theta_0$ is contained in the credible intervals, but also that the posterior of the functional $\Psi(\theta)$ has a Gaussian-like unimodal shape around the origin, providing evidence for the empirical efficacy of our results.

\section{Proofs for Section 3}
\label{sec:proof3}

We first define some notations that are standard in studying the spectra of bounded linear operators between Hilbert spaces.
For a compact linear operator $K:X\to Y$, we define the spectral projector 
\begin{equation}
\label{eq:spectral_projector}
E_tx := \sum_{n=1,2,\ldots:\sigma^2_n<t}\langle x, v_n\rangle_X~v_n + P(x), 
\end{equation}
where $(\sigma_n^2,v_n)_{n=1}^{\infty}$ are eigenvalues and eigenfunctions of $K^*K$, and $P$ is the orthogonal projector onto the nullspace of $K^*K$.
For any piecewise continuous function $f:\R\to\R$, we define 
\begin{equation}
\label{eq:spectral_integral1}
\int_\R f(t)~d\|E_tx\|^2_X := \sum_{n=1}^{\infty} f(\sigma_n^2)|\langle x, v_n\rangle_X|^2,
\end{equation}
for which Lebesgue dominated convergence theorem and integration by parts hold \cite[Section~2.3]{engl1996regularization}.
By extension, we also define $f(K^*K) := \sum_{n=1}^{\infty}f(\sigma_n^2)\langle \cdot, v_n\rangle_X~v_n$ and
\begin{equation}
\label{eq:spectral_integral2}
\|f(K^*K)x\|_{X} = \sqrt{\int_0^{\|K^*K\|_{\rm op}^2}f^2(t)~d\|E_tx\|_X^2}.
\end{equation}
We refer the reader to \cite[Chapter~2.3]{engl1989convergence} for more details on the spectral projection theory.

\subsection{Proof of Theorem \ref{thm:freq_coverage}}

To begin the proofs of Theorems \ref{thm:freq_coverage} and \ref{thm:order_estimates}, we state two lemmas on the asymptotic weak convergence of re-scaled, re-centered posterior measures to a Gaussian limit.
Their proofs follow closely the arguments of \cite{monard2021statistical} (for Lemma \ref{lem:gaussian_limit}, see also \cite{castillo2015bernstein} for similar ideas). The proofs are included in the Supplement \cite{supplement}.

\begin{lemma} \label{lem:gaussian_limit}
    Suppose Assumptions \ref{assmp:truth_prior}-\ref{assmp:posterior_contraction} hold, $\psi\notin{\rm ker}(\Lambda^{-1/2})$,
    and that \eqref{eq:stoch_unif_control} and \eqref{eq:lan_ratio_control} of Theorem \ref{thm:freq_coverage} hold true. Then we have
    \begin{equation}
        d_{weak}(Law(s_N^{-1}\{ 
        \Psi(\theta-\theta_0) - \bigpsi + b_N
        \} \mid \data), Normal(0,1)) \stackrel{P_{\theta_0}^N}{\to} 0,
    \end{equation}
    where
    \begin{equation} \label{eq:seq_lf}
        \bigpsi := \sum_{i=1}^{N} \epsilon_i\opinfo(\barpsi)(X_i).
    \end{equation}
\end{lemma}

\begin{lemma} \label{lem:posterior_mean_conv}
    Under the conditions of Lemma \ref{lem:gaussian_limit},
    if $|b_N|/s_N = o_N(1)$, then we have
    \begin{equation} \label{eq:pm_centering}
        d_{weak}(Law(s_N^{-1}\Psi(\theta - \bartheta),\;
        Normal(0,1)) \stackrel{\ptheta}{\to} 0,
    \end{equation}
    and also
    \begin{equation}
    d_{weak}(
    Law(t_N^{-1}\Psi(\bartheta - \theta_0 ) ),\;
    \mathrm{Normal}(0, 1)
    ) \to 0, \label{eq:limit_mean}
    \end{equation}
    where 
    \begin{equation} \label{eq:t_se}
        t_N := \sqrt{N}\|\opinfo(\barpsi)\|_{\ltwo(\X)}.
    \end{equation}
\end{lemma}

We now recall that, based on the arguments in the proof of \cite[Corollary~7.3.22]{gine2021mathematical}, 
whenever a sequence of probability measures $\mu_k$ weakly converges to a non-degenerate Gaussian measure on $\R$, the collection of rays $\{(-\infty,r]\}_{r\in\R}$ forms a ``uniformity class,'' i.e., we have $\sup_{r\in\R}|\mu_k((-\infty,r]) - \Phi(r)|\stackrel{k\to\infty}{\to}0$,
for $\Phi$ the standard normal cumulative distribution function.
Thus, under the weak convergence conclusions from Lemma \ref{lem:gaussian_limit} and 
Lemma \ref{lem:posterior_mean_conv}, we have
\begin{equation}
    \label{eq:weak_uniformized}
    \sup_{r}
    \left\{|
    \nprior( s_N^{-1}|\Psi(\theta-\bartheta)| \leq r 
    \mid \data) - \Phi(r)
    |\vee|\ptheta(t_N^{-1}
    |\Psi(\bartheta - \theta_0)| \leq r) - \Phi(r)| \right\} \stackrel{\ptheta}{\to} 0.
\end{equation}
Set $z_{1-\gamma/2}$
to be the number such that 
$\Phi(z_{1-\gamma/2}) - \Phi(-z_{1-\gamma/2}) = 1-\gamma$.
Set $Q_{N,\gamma}$ to be a number such that $\nprior\{|s_N^{-1}\Psi(\theta - \bartheta)|\leq Q_{N,\gamma}\mid\data\} = 1-\gamma$.
We then obtain the convergence
$s_N^{-1}Q_{N,\gamma}\stackrel{\ptheta}{\to} z_{1-\gamma/2}$ by the continuous mapping theorem.
The conclusion follows by noting
\begin{eqnarray*}
    \ptheta(| \Psi(\bartheta - \theta_0) |\leq Q_{N,\gamma} \mid \data)
     & = & \ptheta\left(t_N^{-1}|\Psi(\bartheta - \theta_0)| \leq \frac{s_N}{t_N}
    s_N^{-1}Q_{N,\gamma}\right)\\
    & \geq & \ptheta\left(t_N^{-1}|\Psi(\bartheta - \theta_0)| \leq
    s_N^{-1}Q_{N,\gamma}\right),
\end{eqnarray*}
since, by Lemma \ref{lem:st_order_equivalence} below, $s_N/t_N$ is uniformly bounded from below by 1.
Hence, the limit infimum of the left hand side is bounded from below by $\Phi(z_{1-\gamma}/2)$ as $N\to\infty$. 
On the other hand, if $\psi\in R(\opinfo^*)$, Lemma \ref{lem:st_order_equivalence} below shows that $s_N/t_N\to 1$, so we can replace the inequality above with the statement 
$\ptheta(|\Psi(\bartheta-\theta_0)|\leq Q_{N,\gamma}|\data)\to \Phi(z_{1-\gamma/2})$.

The last assertion about posterior contraction rate follows because we may write
\begin{equation*}
    \nprior\{|\Psi(\theta - \theta_0)| \leq M_Ns_N\} = 
    \nprior\{s_N^{-1}|\{\Psi(\theta - \theta_0) - \bigpsi\}|
    \leq M_N - s_N^{-1}\bigpsi\},
\end{equation*}
where $s_N^{-1}\bigpsi = O_{\ptheta}(1)$,
since $t_N/s_N\leq 1$ for every $N$ (Lemma \ref{lem:st_order_equivalence}) and $t_N^{-1}\bigpsi$ weakly converges to a standard normal distribution,
a claim which is proven by the triangular central limit theorem (cf. proof of Lemma \ref{lem:posterior_mean_conv} \cite{supplement}).
Thus, for any sequence $M_N\to\infty$, we may write
\begin{equation*}
    \nprior\{|\Psi(\theta-\theta_0)| \leq M_Ns_N\} \geq
    \nprior\{s_N^{-1}|\{\Psi(\theta - \theta_0) - \bigpsi| \leq M_N/2\}
\end{equation*}
and conclude the lower bound converges in probability to 1 by Lemma \ref{lem:gaussian_limit} and its corollary \eqref{eq:weak_uniformized}. 

\begin{lemma} \label{lem:st_order_equivalence}
    $t_N/s_N$ is a positive sequence bounded above by 1.
    If $\psi\in R(\opinfo^*)$, then $t_N/s_N\to 1$.
\end{lemma}

\begin{proof}
    We will express $s_N$ and $t_N$ in a spectral notation involving the spectrum of $\opcomb^*\opcomb$.
    Define a sequence of functions $f_N:(0,\infty)\to(0,\infty)$ by $f_N(t) = \frac{1}{\sqrt{1+N\tau_N^2t}}$ and $g_N(t) = \sqrt{N}\tau_N\sqrt{t}f_N^2(t)$.
    Based on the discrete spectrum of a trace-class operator $\opcomb^*\opcomb$, we can re-write the definitions \eqref{eq:asymp_scale} and \eqref{eq:t_se} using the notation from \eqref{eq:spectral_integral2}:
    \begin{eqnarray}
    \label{eq:eq_sn_spec}
    s_N^2 & = & \tau_N^2 \|f_N(\opcomb^*\opcomb) \lpsi\|_{\ltwo(\Z)}^2,\\
    t_N^2 & = & \tau_N^2\|g_N(\opcomb^*\opcomb)\lpsi\|_{\ltwo(\Z)}^2.
    \end{eqnarray}
    Since $f_N$ strictly dominates $g_N$ on $(0,\infty)$ for each $N$, we deduce that $t_N/s_N$ is a positive sequence bounded from above by 1.    
    Next, we prove that $t_N/s_N\to 1$ if $\psi\in R(\opinfo^*)$.
    We note that $s_N/t_N \to 1$ if and only if, for every $c > 0$,
    we have
    $\int_0^{\frac{c}{N\tau_N^2}} \frac{1}{1 + N\tau_N^2t}~dF^2_{\lpsi}(t) \ll \tau_N^{-2}s_N^2.$
    The arguments for the proof of this equivalence can be found in the proof of \cite[Theorem~5.4]{knapik2011bayesian} (see display (7.15) therein and below). 
    Therefore, it is sufficient to prove the above display.
    But $\psi \in R(\opinfo^*)$ if and only if $\lpsi \in R(\opcomb^*)$, since $\opinfo$ is injective on the image of $\Lambda^{-1/2}$ by Assumption \ref{assmp:deriv_stability}, so the left-hand side is $o(\frac{1}{N\tau_N^2})$ by Lemma \ref{lem:source_spectral}. 
    On the other hand, Lemma \ref{lem:order_s} below will show that whenever $\psi\in R(\opinfo^*)$, we have $s_N^2 \asymp \frac{1}{N\tau_N^2}$.
    Hence, we conclude that if $\psi\in R(\opinfo^*)$, then $s_N/t_N\to 1$.
\end{proof}

\subsection{Proof of Theorem \ref{thm:order_estimates}}

We want to show that $|b_N|/s_N \ll 1$ under conditions \ref{condition1:tau1}--\ref{condition3:bvm_rate} and derive the stated upper bound estimates for $s_N$ depending on whether $\psi\in R(\opinfo^*)$.
Once we have precise enough estimates of the sequences $s_N$ and $b_N$, the proof of Theorem \ref{thm:order_estimates} should be a simple verification.
As emphasized in the introduction, the following results draw on proving an inequality of the form given in Definition \ref{def:vi} and using its consequences.
Known consequences of such an inequality are summarized in Appendix \ref{appendix:vi_lemmas}.

Concretely, our first goal is to verify the following inequality for some exponent $\eta\in (0,1]$ and positive constants $b = b(\Lambda,\theta_0,\psi,\eta)$, $C = C(\Lambda,\theta_0,\psi,\eta)$:
\begin{equation}
     \label{eq:vi}
     b\|\lpsi - u\|_{\ltwo(\Z)}^2 \leq
     \|\lpsi\|_{\ltwo(\Z)}^2 - \|u\|_{\ltwo(\Z)}^2 + 
     C\|\opcomb ( \lpsi - u )\|_{\ltwo(\X)}^\eta.
\end{equation}
This is the inequality given in Definition \ref{def:vi} specialized to our context.
The importance of Assumption \ref{assmp:deriv_stability} lies in its verification.
For the rest of this section, we reserve notation $\delta(s)$ for the exponent yielding the statement in Assumption \ref{assmp:deriv_stability} for each $\Hs^s,\; s \geq 0$.

\begin{lemma} \label{lem:vi_sufficient}
    Suppose $\psi \in \ltwo(\Z)$.
    Under Assumption \ref{assmp:deriv_stability}, 
    we have \eqref{eq:vi} for some positive constants $b = b(\Lambda,\theta_0,\psi)\leq 1,\; 
    C = C(\Lambda,\theta_0,\psi) > 0$ and 
    $\eta = \frac{2\delta(1)}{1 + \delta(1)} \in (0,1]$.
\end{lemma}

\begin{proof}
    By self-adjointness of $\Lambda$,
    $|\langle \Lambda^{-1/2}\psi,u\rangle_{\ltwo(\Z)}| 
        \leq \|\Lambda^{-1/2}u\|_{\ltwo(\Z)}
        \|\psi\|_{\ltwo(\Z)}$
    for every $u\in B_{\ltwo(\Z)}(1)$.
    For every such $u$, $\Lambda^{-1/2}u\in B_{\Hs^{1}}(M)$ for some $M=M(s)$, 
    so the first norm is bounded above by 
    $\|\opinfo\Lambda^{-1/2}u\|_{\ltwo(\X)}^{\delta(1)}$ up to a constant uniform in $u\in B_{\ltwo(\Z)}(1)$ by Assumption \ref{assmp:deriv_stability}.
    The conclusion then follows by Lemma \ref{lem:vi_characterization}.
\end{proof}

We will now prove consequences of \eqref{eq:vi} that allow us to estimate $s_N$ \eqref{eq:asymp_scale} and $b_N$ \eqref{eq:bias}.

\subsubsection{Estimating 
\texorpdfstring{$s_N$}{sN}}

There are two possible regimes regarding the choice of functional: either $\psi\in R(\opinfo^*)$, if and only if $\lpsi\in R(\opcomb^*)$;
or, $\lpsi\notin R(\opcomb^*)$. 
The first scenario is the ``efficient regime'' in light of \cite[Theorem~3.1]{van1991differentiable}. 
We first focus on the second ``non-efficient regime,'' which can be further split into two scenarios depending on the quantity
\begin{equation}
\label{eq:eta_star}
    \eta^* := \sup
    \{\eta\in(0,1] : \text{There exist }C\geq0,\;b\leq 1
    \text{ for which ineq. (5.6) holds}\}.
\end{equation} 
Since we have Lemma \ref{lem:vi_sufficient}, we know that there is at least some $\eta>0$ for which \eqref{eq:vi} is met, and thus the above definition makes sense, under Assumption \ref{assmp:deriv_stability}.
Here is the importance of $\eta^*$:
Lemma \ref{lem:vi_source_equiv} tells us that when $\psi \notin R(\opinfo^*)$, \eqref{eq:vi} with $\eta=1$ cannot be satisfied for any constants.
However, it is still possible that \eqref{eq:vi} is met for every $\eta < 1$. Lemmas \ref{lem:dist_upper_bound} and \ref{lem:dist_source_equiv} imply that in this case $\psi\in R((\opcomb^*\opcomb)^{\mu/2})$ for every $\mu < 1$, so $\psi$ is ``almost efficiently estimable.''

We now derive estimates for $s_N$ depending on whether $\eta^* < 1$, $\eta^* = 1$ but $\psi \notin R(\opinfo^*)$, or $\psi \in R(\opinfo^*)$.

\begin{lemma} \label{lem:order_s}
    The following holds true if we assume Assumption \ref{assmp:deriv_stability} for cases (i) and (ii) below.
    \begin{enumerate}[label=(\roman*)]
        \item Suppose $\psi\notin R(\opinfo^*)$ and $\eta^* < 1$. 
        Set $\kappa^* = \frac{\eta^*}{2(2-\eta^*)} < \frac{1}{2}$.
        Then, for any $\kappal \in (0,\kappa^*)$ and
        $\kappau \in (\kappa^*,\frac{1}{2})$,
        \begin{equation}
        \label{eq:order_s1}
        \tau_N(\sqrt{N}\tau_N)^{-2\kappau} \lesssim 
            s_N
            \lesssim
        \tau_N(\sqrt{N}\tau_N)^{-2\kappal}.
        \end{equation}
        \item Suppose $\psi \notin R(\opinfo^*)$ and $\eta^* = 1$. Then $\sqrt{N}s_N \asymp S(N)$
        for some sequence $S(N)$ that grows slower than polynomially in $N$, i.e., $N^{-\nu}S(N)\to 0$ for every $\nu > 0$.
        \item Suppose $\psi\in R(\opinfo^*)$. Then $\sqrt{N}s_N\to 
        \|(\opinfo^*\opinfo)^{-1/2}\psi\|_{\ltwo(\Z)}$.
    \end{enumerate}
\end{lemma}
The proof is based on a supporting lemma. To state the lemma, we define the following quantities (cf. Appendix \ref{appendix:vi_lemmas}): 
    the distance function,
    $d_{\lpsi}(R) = \inf_{v\in B_{\ltwo(\Z)}(R)}
    \left\| \lpsi - \sqrt{\opcomb^*\opcomb}v \right\|_{\ltwo(Z)}$,
    and the spectral distribution function, defined by a standard spectral projection:
    $F_{\lpsi}(t) := \|E_t \lpsi\|$, using the notation from \eqref{eq:spectral_projector}.
    Viewed as a function of $R\in[0,\infty)$, $d_x$ is decreasing and continuous \cite[p.7]{flemming2011sharp}.
For each finite $R$, a minimizer $v$ in \eqref{eq:dist_fn} exists, since the functional is weakly lower semicontinuous and $B_X(R)$ is weakly compact on a reflexive space $X$.
It is clear that there exists a finite $R$ and a minimizer $v$ such that a zero infimum is attained if $u \in R(\sqrt{H})$.
More properties of these quantities and equivalences are summarized in Appendix \ref{appendix:vi_lemmas}, which we will be using below.

\begin{lemma} \label{lem:dist_order}
    Let $\psi\notin R(\opinfo^*)$ and assume $\eta^*>0$. For any $\underline{\eta} \in (0,\eta^*)$ and $\etau\in (\eta^*, 1)$ (with $\etau=1$ if $\eta^*=1$), we have
    \begin{equation}
    \label{eq:dist_order}
        R^{-\frac{\etau}{2 - 2\etau}}
        \lesssim d_{\lpsi}(R) \lesssim R^{-\frac{\etal}{2-2\etal}}\quad\text{as }
        R\to\infty.
    \end{equation}
    Furthermore, $\lpsi \in R((\opcomb^*\opcomb)^{\mu})$ for every $\mu \in (0, \frac{\etal}{2(2-\etal)})$ and
    \begin{equation}
    \label{eq:spec_order}
        \epsilon^{\kappal} \lesssim 
        F_{\lpsi}(\epsilon) \lesssim \epsilon^{\kappau}
        \text{ as }\epsilon\to 0,
    \end{equation}
    for $\kappal = \frac{\etau}{2(2-\etau)}$ and $\kappau = \frac{\etal}{2(2-\etal)}$.
\end{lemma}

\begin{proof}[Proof of Lemma \ref{lem:order_s}]
    Recall that $s_N$ has a spectral representation \eqref{eq:eq_sn_spec}:
    for a sequence of functions $f_N(t) = \frac{1}{\sqrt{1+N\tau_N^2t^2}}$, we have $s_N^2/\tau_N^2 = \|f_N^2(\opcomb^*\opcomb)\lpsi\|^2$.
    \begin{enumerate} [label=(\roman*)]
        \item  First, we prove \eqref{eq:order_s1}.
        For each $N$, $f_N$ in \eqref{eq:eq_sn_spec} is continuous, bounded above by 1, and strictly decreasing.
        By strict monotonicity, whenever $0 \leq t \leq \frac{1}{N\tau_N^2}$, we have $f^2_N(t)\geq 1/2$. 
        As a result, we have
        \begin{equation} \label{eq:lower_bound}
        s_N^2 \geq \tau_N^2
        \int_{0}^{\frac{1}{N\tau_N^2}} f_N^2(t)~dF^2_{\lpsi}(t)
        \geq \frac{\tau_N^2}{2}F^2_{\lpsi}\left(\frac{1}{N\tau_N^2}\right).
        \end{equation}
    Combining the above with the lower bound of \eqref{eq:spec_order} yields the lower bound of \eqref{eq:order_s1}. 
    To derive the upper bound, we proceed similarly as the proof of \cite[Theorem~5.5]{hofmann2007analysis} and use the definition of function $\Phi$ given in Appendix \ref{appendix:vi_lemmas} \eqref{eq:fixed_point}:
    $\Phi_{\lpsi}(R) := \frac{d_{\lpsi}(R)}{R}.$
    By the continuity and monotonicity of $\Phi_{\lpsi}:(0,\infty)\to(0,\infty)$ mentioned therein, 
    which hold here since $\lpsi\notin R(\opcomb^*)$,
    for every $\epsilon > 0$, the equation $\Phi_{\lpsi}(R) = \epsilon$ has a unique solution $\widehat{R} = \widehat{R}(\epsilon)$.
    We choose $\epsilon_N = \frac{1}{N\tau_N^2}$ and set $\widehat{R}_N $ the solution for each $\epsilon_N$. 
    We denote by $v^\dagger_N$ any minimizer appearing in the definition of \eqref{eq:dist_fn} for $R=\widehat{R}$.
    By triangle inequality, we obtain:
    \begin{eqnarray*}
        \tau_N^{-1}s_N & = & 
        \|f_N(\opcomb^*\opcomb)\lpsi\|_{\ltwo(\Z)} \\
        & \leq & 
        \|f_N(\opcomb^*\opcomb)(\lpsi - \opcomb^* v^\dagger_N)\|_{\ltwo(\Z)} + 
        \|f_N(\opcomb^*\opcomb)\opcomb^* v^\dagger_N\|_{\ltwo(\Z)} \\
        & \leq & 
        d_{\lpsi}(\widehat{R}_N) + \frac{1}{\sqrt{N}\tau_N}\widehat{R}_N
        \leq 2 d_{\lpsi}\left(\Phi_{\lpsi}^{-1}\left(\frac{1}{\sqrt{N}\tau_N}\right)\right).
    \end{eqnarray*}
    In the third inequality, we used the uniform boundedness of $f_N$ (cf. Lemma \ref{lem:perturbation_lemma}) and that the map $t\mapsto tf_N(t)$ is bounded by $(\sqrt{N}\tau_N)^{-1}$ for each $N$. 
    Since $F_{\lpsi}(\epsilon)$ decays polynomially in $\epsilon$ over a neighborhood of zero from Lemma \ref{lem:dist_order} \eqref{eq:spec_order}, 
    Lemma \ref{lem:spec_fp_bounds} yields
    \begin{equation*}
    \tau_N^{-1}s_N \lesssim d_{\lpsi}\left(\Phi_{\lpsi}^{-1}\left(\frac{1}{\sqrt{N}\tau_N}\right)\right) \lesssim F_{\lpsi}\left(\frac{1}{N\tau_N^2}\right),
    \end{equation*}
    with the hidden constants independent of $\etal,\etau$.
    The asserted bound follows from \eqref{eq:spec_order}.
    \item Now suppose $\psi \notin R(\opinfo^*)$ but $\eta^* = 1$. 
    By Lemma \ref{lem:dist_order}, we have $\lpsi \in R((\opcomb^*\opcomb)^{\nu/2})$ for every $\nu < 1$.
    We can now show that for every $\nu < 1$,
    \begin{equation*}
        N^\nu s_N^2 = \tau_N^{2-2\nu}\int_0^\infty
        \frac{N^\nu\tau_N^{2\nu}t^\nu}{1+N\tau_N^2t}\frac{dF_{\lpsi}^2(t)}{t^\nu} = o_N(1).
    \end{equation*}
    The first equality is a simple consequence of the spectral representation \eqref{eq:eq_sn_spec}; 
    the second inequality then follows from the dominated convergence theorem, applicable because the sequence of integrand functions pointwise converges to 0 and is dominated by an envelope function $t\mapsto \frac{1}{t^\nu}$, integrable by the assumption $\lpsi\in R((\opcomb^*\opcomb)^{\nu/2})$.
    On the other hand, by monotonicity of the integrand below, we also have
    \begin{equation*}
        Ns_N^2 \geq \int_{\frac{1}{N\tau_N}^2}^\infty \frac{N\tau_N^2t}{1 + N\tau_N^2t}~\frac{dF_{\lpsi}^2(t)}{t}
        \geq \frac{N\tau_N^2K^2}{1+ N\tau_N^2K^2}
        \int_{\frac{1}{N\tau_N^2}}^{\infty}
        \frac{dF_{\lpsi}^2(t)}{t^2},
    \end{equation*}
    for $K$ the operator norm of $\opcomb^*\opcomb$.
    As $N\to\infty$, the lower bound diverges, since we know $\lpsi\notin R(\opinfo^*)$. Therefore, we must have $Ns_N =: S(N)$ growing slower than polynomially in $N$. 
    \item Finally, suppose $\psi \in R(\opinfo^*)$ so that $\lpsi \in R(\opcomb^*)$. Then
    \begin{equation*}
    Ns_N^2 = N\tau_N^2\int_0^\infty
    \frac{1}{1 + N\tau_N^2t}~dF^2_{\lpsi}(t)
    = \int_0^\infty \frac{N\tau_N^2t}{1 + N\tau_N^2t}
    ~\frac{dF^2_{\lpsi}(t)}{t}.
    \end{equation*}
    Applying the dominated convergence theorem with $t\mapsto t^{-1}$ as an integrable envelope, we conclude the above converges to $\|(\opcomb^*\opcomb)^{-1/2}\lpsi\|_{\ltwo(\Z)}^2= \|(\opinfo^*\opinfo)^{-1/2}\psi\|_{\ltwo(\Z)}^2$.
    \end{enumerate}
\end{proof}

\begin{proof}[Proof of Lemma \ref{lem:dist_order}]
    We first prove the assertion \eqref{eq:dist_order}.
    Since the variational inequality is met for any $\etal < \eta^*$, the upper bound follows from the implication of Lemma \ref{lem:dist_upper_bound} and continuity of $d_x$.
    The lower bound is trivial if $\eta^*=1$, so assume $\eta^*<1$.
    By definition \eqref{eq:eta_star}, an inequality \eqref{eq:vi} with $\eta=\etau$ cannot be met with any positive constants $b,C$.
    By the contrapositive of Lemma \ref{lem:dist_lower_bound},
    no matter how large $\overline{R}>0$ and $c'>0$ are, there is some $R>\overline{R}$ such that $d_{\lpsi}(R) > c'R^{-\frac{\etau}{2-2\etau}}$.
    But since $\lpsi\notin R(\opcomb^*)$ when $\psi\notin R(\opinfo^*)$, Lemma \ref{lem:dist_upper_bound} implies that $d_{\lpsi}$ does not attain zero for any finite $R$ and is strictly decreasing outside $(0,\overline{R}')$ for some $\overline{R}'>0$.
    By continuity, we conclude $d_x(R)\gtrsim R^{-\frac{\etau}{2-2\etau}}$ as $R\to\infty$.
    
    The second assertion \eqref{eq:spec_order} follows from Lemma \ref{lem:spec_dist_equiv}, along with the fact that $F_{\lpsi}(\epsilon)$ is increasing in $\epsilon$ and $d_{\lpsi}(R)$ is decreasing in $R$.    
    The assertion on the range inclusion of $\lpsi$ for every $\mu<\frac{\etal}{2(2-\etal)}$ follows from Lemma \ref{lem:dist_source_equiv}, setting $\kappa$ in the display therein to equal $\frac{\etal}{2(2-\etal)}$.
\end{proof}

\subsubsection{Estimating 
\texorpdfstring{$b_N$}{bN}}

Next, we estimate the ratio of $|b_N|$ to $s_N$, again depending on whether $\psi\in R(\opinfo^*)$ or $\eta^*$ \eqref{eq:eta_star} is 1.

\begin{lemma}
    \label{lem:order_sb}
    Under Assumption \ref{assmp:deriv_stability},
    let $\psi \in \Hs^q$ and $\theta_0 \in B_{\Hs^\beta}(M)$ for some $q\geq 0,\; \beta\geq 1$ and $M > 0$.
    \begin{enumerate}[label=(\roman*)]
        \item If $\psi\notin R(\opinfo^*)$ and $\eta^* < 1$, then
        \begin{equation}
        \label{eq:order_b1}
        \frac{|b_N|}{s_N} \lesssim
        \|\theta_0\|_{\Hs^\beta}\tau_N^{-1}
        \left( \frac{1}{\sqrt{N}\tau_N} \right)^{(1 + 2\kappau)\delta(2 + q)(\frac{\beta}{2}\wedge 1) - 2\kappau}
        \end{equation}
    for every $\kappau \in (\kappa^*, \frac{1}{2})$, with $\kappa^* := \frac{\eta^*}{2(2-\eta^*)}$.
    \item If $\psi\notin R(\opinfo^*)$ and $\eta^* = 1$, 
    then for some $S(N)\to\infty$ slower than any polynomial in $N$,
    \begin{equation}
        \label{eq:order_b2}
        \frac{|b_N|}{s_N} \lesssim
        \|\theta_0\|_{\Hs^\beta}\tau_N^{-1}
        \left( \frac{S(N)}{\sqrt{N}\tau_N} \right)^{2\delta(2 + q)(\frac{\beta}{2}\wedge 1) - 1}.
    \end{equation}
    \item If $\psi\in R(\opinfo^*)$, then we can replace $2\kappau$ with 1 in \eqref{eq:order_b1}.
    \end{enumerate}
    If $\delta(2 + q)(\frac{\beta}{2}\wedge 1) > \frac{1}{2}$, then the exponents of $\frac{1}{\sqrt{N}\tau_N}$ in all three cases can be chosen to be positive.
\end{lemma}

\begin{proof}
    We use an interpolation argument to handle all three cases.
    Recall the definition of $b_N$ \eqref{eq:bias}.
    Since $\theta_0\in \Hs^1$ (Assumption \ref{assmp:truth_prior}),
    by Cauchy-Schwarz inequality and self-adjointness of $\Lambda$,
    \begin{equation}
    \label{eq:cauchy_schwarz}
        |b_N| \leq \|\theta_0\|_{\Hs^\beta}
        \|f_N^2(\opcomb^*\opcomb)\lpsi\|_{\Hs^{1 - \beta}}
    \end{equation}
    for every $\beta\geq 0$ (with $\|x\|_{\Hs^\beta} = \infty$ if $x\notin\Hs^\beta$), and $f_N(t) = \frac{t}{\sqrt{1+N\tau_N^2t}},\; t\geq 0$, as defined in \eqref{eq:eq_sn_spec}.
    Since $\lpsi\in \Hs^{1+q}$ for $\psi\in\Hs^q$, Lemma \ref{lem:perturbation_lemma} implies $f_N^2(\lpsi) \in \Hs^{1+q}$.
    We now obtain the following using interpolation when $\beta\in[0,2]$ (Theorem \ref{thm:hilbert_scales}):
    \begin{equation}
    \label{eq:bias_ineq1}
        \|f_N^2(\opcomb^*\opcomb)\lpsi\|_{\Hs^{1 - \beta}} \leq
        \|f_N^2(\opcomb^*\opcomb)\lpsi\|_{\Hs^{-1}}^{\frac{\beta}{2}\wedge1}
        \|f_N^2(\opcomb^*\opcomb)\lpsi\|_{\Hs^{1}}^{(1-\frac{\beta}{2})\vee0}.
    \end{equation}
    In fact, the above is also true when $\beta \geq 2$, as then $\norm_{\Hs^{1-\beta}}\leq\norm_{\Hs^{-1}}$.
    Furthermore, the $\Hs^{1}$-norm appearing on the right-hand side is bounded above uniformly in $N$ by the uniform bound on the sequence $f_N$.
    On the other hand, $\Lambda^{-1/2}f^2(\opcomb^*\opcomb)\lpsi\in \Hs^{2+q}$, so using Assumption \ref{assmp:deriv_stability},
    \begin{equation}
    \label{eq:bias_ineq2}
    \|f_N^2(\opcomb^*\opcomb)\lpsi\|_{\Hs^{-1}} =
        \|\Lambda^{-1/2} f^2(\opcomb^*\opcomb) \lpsi\|_{\ltwo(\Z)} \lesssim 
        \|\opcomb f^2(\opcomb^*\opcomb)\lpsi\|_{\ltwo(\X)}^{\delta(2 + q)}.
    \end{equation}
    Define (again) a sequence of functions $g_N:(0,\infty)\to(0,\infty)$ by $g_N(t) = \frac{\sqrt{N}\tau_N t}{1+N\tau_N\sqrt{t}}$, so that we can re-write
    $\|\opcomb f_N^2(\opcomb^*\opcomb)\lpsi\|_{\ltwo(\Z)} = \frac{1}{\sqrt{N}\tau_N}
    \| g_N(\opcomb^*\opcomb)\lpsi \|_{\ltwo(\X)}
    $.
    Then $g_N$ is dominated by $f_N$ \eqref{eq:eq_sn_spec} for each $N$, but then $\|f_N(\opcomb^*\opcomb)\lpsi\|_{\ltwo(\X)} = s_N/\tau_N$ \eqref{eq:eq_sn_spec}. 
    Therefore, combining \eqref{eq:cauchy_schwarz} with \eqref{eq:bias_ineq1} and \eqref{eq:bias_ineq2}, we conclude
    \begin{equation*}
        \frac{|b_N|}{s_N} \lesssim
        \|\theta_0\|_{\Hs^{\beta}}
        \left( \frac{1}{\sqrt{N}\tau_N} \right)^{\delta(2 + q)
        (\frac{\beta}{2}\wedge 1)}
        \left( \frac{s_N}{\tau_N} \right)^{\delta(2 + q)
        (\frac{\beta}{2}\wedge 1)-1}\tau_N^{-1}.
    \end{equation*}
    Since we always have $\delta(2+q)(\beta/2\wedge 1)\leq 1$, we can plug in the lower bounds for $s_N$ from Lemma \ref{lem:order_s} to obtain our three bounds by cases.
    In particular, for the case when $\eta^* < 1$, the exponent of $1/\sqrt{N}\tau_N$ in \eqref{eq:order_b1} is positive when $\delta(2 + q)(\frac{\beta}{2}\wedge 1) > \frac{2\kappau}{1 + 2\kappau}$. This is guaranteed if $\delta(2 + q)(\frac{\beta}{2}\wedge 1)> \frac{1}{2}$, as $2\kappau \leq 1$, in which case the exponent in \eqref{eq:order_b2} is also positive.
\end{proof}

\begin{lemma}
    \label{lem:order_sb_stronger}
    If $\psi \in 
    R(\opinfo^*\opinfo\Lambda^{-1/2})$ and $ \beta=1$, then 
$|b_N|\lesssim 
    \frac{\|\theta_0\|_{\Hs^1}}{N^2\tau_N^4}.$
\end{lemma}

\begin{proof}
    Since $\psi \in R(\opinfo^*)$, Lemma \ref{lem:order_s} (iii) yields $s_N \asymp 1/\sqrt{N}$.
    Inequality \eqref{eq:cauchy_schwarz} remains true for $\beta=1$, and we now have $|b_N| \leq \|\theta\|_{\Hs^1}\|f^2(\opcomb^*\opcomb)\lpsi\|_{\ltwo(\Z)}$.
    Since $\lpsi\in R(\opcomb^*\opcomb)$, we have
    \begin{equation*}
        N^2\tau_N^4\|f^2(\opcomb^*\opcomb)\lpsi\|_{\ltwo(\Z)}
        =\int_0^\infty
        \frac{N^2\tau_N^4t^2}{(1+N\tau_N^2t)^2}~
        \frac{dF_{\lpsi}^2(t)}{t}.
    \end{equation*}
    The above display is convergent by dominated convergence.
\end{proof}

Theorem \ref{thm:order_estimates} is now an easy consequence of Lemmas \ref{lem:order_s}, \ref{lem:order_sb} and \ref{lem:st_order_equivalence}.

\begin{proof}[Proof of Theorem \ref{thm:order_estimates}]
First, we verify that the conditions \ref{condition1:tau1}--\ref{condition3:bvm_rate} imply that $|b_N|/s_N\ll 1$.
Some calculations reveal that if $\psi\notin R(\opinfo^*)$, then under either condition \ref{condition1:tau1} or condition \ref{condition2:tau_decay},
all three upper bounds of $|b_N|/s_N$ in Lemma \ref{lem:order_sb} are $o_N(1)$. 
Since $2\kappau\leq 1$, the case $\delta(2+q)(\beta/2\wedge 1) = 1$ is trivial.
Under the third condition \ref{condition3:bvm_rate},
Lemma \ref{lem:order_s} yields $s_N\asymp 1/\sqrt{N}$, since $\psi\in R(\opinfo^*)$, and Lemma \ref{lem:order_sb_stronger} yields $|b_N|\lesssim \frac{1}{N^2\tau_N^4}\ll 1/\sqrt{N}$.
Second, again if $\psi\notin R(\opinfo^*)$, under Assumption \ref{assmp:deriv_stability} we have Lemma \ref{lem:vi_sufficient} and $\eta^*\geq\frac{2\delta(1)}{1+\delta(1)}$. Combining it with Lemma \ref{lem:order_s} (i)--(ii) and choosing $2\kappau \leq \delta(1)$,
we obtain $s_N\lesssim \tau_N^{1-\delta}N^{-\delta/2}$ for any $\delta\in (0,\delta(1))$.
If $\psi\in R(\opinfo^*)$, the assertions on $s_N$ follow from Lemma \ref{lem:order_s} (iii).
\end{proof}

\begin{appendix}
\section{Hilbert Scales Result} 
\label{appendix:scales}

Let $A$ be a Hilbert space and $\Lambda^{1/2}$ an adjoint, symmetric, unbounded $A$-valued linear operator, whose domain is dense in $A$, and which satisfies $\|\Lambda^{1/2} x\|\geq c\|x\|^2$ for some $c>0$ and every $x\in A$.
\cite[Proposition~8.19]{engl1996regularization} collects the important properties of the Hilbert scales $(\Hs^p)_{p\in\R}$ generated by $\Lambda$.

\begin{theorem} \label{thm:hilbert_scales}
The Hilbert scales $(\Hs^p)_{p\in\R}$ induced by the operator $\Lambda$ has the following properties.
\begin{enumerate}[label=(\roman*)]
    \item For $-\infty < p < q < \infty$, $\Hs^q$ is dense and continuously embedded in $\Hs^p$.
    \item Let $p,q\in\R$. The operator $\Lambda^{p-q}$ defined on $\ltwo(\Z)$ has a unique extension to $\Hs^{p}$ that is an isomorphism from $\Hs^{p}$ to $\Hs^{q}$. If $p > q$, then the extension restricted to $\Hs^p$ in $\Hs^q$ is self-adjoint and strictly positive. For appropriate extensions, we have $\Lambda^{p-q} = \Lambda^p\Lambda^{-q}$ and $(\Lambda^{-q})^{-1} = \Lambda^{q}$.
    \item If $p\geq 0$, then $\Hs^p = D(\Lambda^{p/2})$ and $\Hs^{-p} = (\Hs^p)^*$.
    \item If $-\infty < p < q < r < \infty$, then we have the interpolation inequality
    \begin{equation*}
        \|f\|_{\Hs^q} \leq \|f\|_{\Hs^p}^{\frac{r-q}{r-p}} \|f\|_{\Hs^r}^{\frac{q-p}{r-p}}.
    \end{equation*}
\end{enumerate}    
\end{theorem}

\section{Consequences of Variational Inequality}
\label{appendix:vi_lemmas}

Our proofs in Section \ref{sec:main} rely on tracking the asymptotics of three quantities $s_N$, $b_N$ and $t_N$, 
whose definitions are given in \eqref{eq:asymp_scale}, \eqref{eq:bias} and \eqref{eq:t_se}.
Our proofs in Section \ref{sec:darcy_flow} in \cite{supplement} involve certain calculations involving the uniform norm of $\barpsi$ \eqref{eq:perturbation} that are needed to verify condition \eqref{eq:stoch_unif_control}.
For both of these purposes, we have drawn on results from previous literature on the use of variational inequality (Definition \ref{def:vi}).
We summarize here the known results in this literature that we have used in our proofs.

Let $X,Y$ be two Hilbert spaces and $A : X\to Y$ an injective, bounded linear operator possessing a non-closed range. Then $H := A^*A$ is a positive self-adjoint operator on $X$.
For each $x\in X$, per \eqref{eq:spectral_projector}-- \eqref{eq:spectral_integral2}, $F_{x}(\epsilon)$ is defined using the spectral projector as
$F_{x}(t) := \|E_t x\|_X$, and we also define $\|f(H)x\|_X := \sqrt{\int_0^{\|H\|_{\rm op}} f^2(t)~d\|E_t x\|_X^2}$ for any piecewise continuous function on $(0,\|H\|_{\rm op})$.

The following result is a special case of \cite[Lemma~3]{flemming2011sharp}.

\begin{lemma} \label{lem:source_spectral}
    If $x\in R(H^\nu)$ then $F_x(\epsilon) = o(\epsilon^\nu)$ as $\epsilon\to 0$.
\end{lemma}

Next, we define the distance function
\begin{equation} \label{eq:dist_fn}
    d_x(R) := 
    \inf_{v \in B_{X}(R)} \| u - \sqrt{H}v \|_{X}.
\end{equation}

The following result is from \cite[Lemma~2.5]{hofmann2006approximate}.

\begin{lemma} \label{lem:dist_prop}
    When $x \notin R(\sqrt{H})$, the function \eqref{eq:dist_fn} is positive, continuous, strictly decreasing, and convex. Furthermore, for each $R$, there exists a unique vector $v^\dagger \in B_{X}(R)$ that is orthogonal to the nullspace of $\sqrt{H}$ and satisfies
    $d_x(R) = \| x - \sqrt{H} v^\dagger \|_{X}$.
\end{lemma}

The following result \cite[Theorem~2]{flemming2011sharp} shows that $F_x(\epsilon)$ and $d_x(R)$ contain the same information.

\begin{theorem} \label{lem:spec_dist_equiv}
    Let $\kappa\in (0,1/2)$. When $\lpsi\notin R(\sqrt{H})$,
    $F_{x}(\epsilon) \asymp \epsilon^{\kappa}$ as $\epsilon \to 0$ if and only if
    $d_{x}(R) \asymp R^{-\frac{2\kappa}{1-2\kappa}}$ as $R\to\infty$.
    The equivalence remains valid if both occurrences of $\asymp$ are replaced by either $\lesssim$ or $\ll$.
\end{theorem}

The following result yields an upper bound for the distance function using the variational inequality. It is a consequence of \cite[Theorem~4.5]{flemming2012solution}, as stated in Remark 4.6 below it.

\begin{theorem} \label{lem:dist_upper_bound}
    Suppose an inequality in Definition \ref{def:vi} is met for a vector $x$ and an operator $A$ with a modifier $f(t) = Ct^\eta$ for some $\eta\in (0,1)$, $C > 0$ and $b\leq 1$. Then, we have
    \begin{equation} 
    \label{eq:decay_rate_benchmark}
        d_x(R) = O(R^{-\frac{\eta}{2-2\eta}}) \text{ as }R\to\infty.
    \end{equation}
\end{theorem}

In fact, variational inequality for a modifier $f(t) = Ct^\eta,\; C>0$ is necessary for an estimate \eqref{eq:decay_rate_benchmark} to hold true, in the following sense. 
\cite{flemming2011sharp} mentions this result as a variation of \cite[Theorem~5.2]{bot2010extension}.

\begin{theorem} \label{lem:dist_lower_bound}
    If a distance function $d_x(R) > 0$ for all $R\geq 0$, then inequality in Definition \ref{def:vi} is met for $x$ with a modifier function $f(t) = CtF^{-1}(t)$ for some $C > 0$, 
    where $F(R) := R^{-1}d_x^2(R)$. 
    In particular, if \eqref{eq:decay_rate_benchmark} holds true, then Definition \ref{def:vi} is satisfied for a modifier function $t\mapsto Ct^{\eta}$.
\end{theorem}

The following result, specializing \cite[Corollary~3.3]{duvelmeyer2007range}, shows that the exponent $\nu\in(0,1)$ for which $x\in R(H^{\nu/2})$ bijectively corresponds to the best exponent for the decay of the distance function \eqref{eq:dist_fn}.

\begin{theorem} \label{lem:dist_source_equiv}
    Suppose $x\notin R(\sqrt{H})$. Suppose that there exists some $\kappa\in(0,1/2)$ such that
    $d_{x}(R)  = O(R^{-\frac{2\kappa}{1-2\kappa}})\text{ as }R\to\infty$.
    Then $\lpsi\in R(H^{\mu})$ for all $\mu\in (0,\kappa)$. Moreover, defining
    \begin{equation*}
        \kappa_{sup} := \sup
        \left\{\kappa > 0 : d_{x}(R) = O(R^{-\frac{2\kappa}{1-2\kappa}})
        \text{ as }R\to\infty\right\}\in (0,1/2],
    \end{equation*}
    we have $\kappa_{sup} = \sup\{\mu > 0 : x \in R(H^{\mu})\}$ whenever $A$ is compact.
\end{theorem}

In the above, it is possible that $\kappa_{sup} = 1/2$ but $x\notin R(\sqrt{H})$. The following equivalence between Definition \ref{def:vi} and the range inclusion $x\in R(\sqrt{H})$ is well-known (e.g., \cite[Proposition~31]{flemming2012solution}).

\begin{theorem} \label{lem:vi_source_equiv}
    $x\in X$ satisfies an inequality in Definition \ref{def:vi} for a constant $b\leq 1$ and a modifier function $f(t) = Ct$ for some $C \geq 0$ if and only if there exists some $w$ with $\|w\|_X \leq C/2$ such that $x = \sqrt{H}w$.
\end{theorem}

The following characterizations of Definition \ref{def:vi} \cite[Propositions~2.9--10]{flemming2012solution} are also useful.

\begin{theorem}
\label{lem:vi_characterization}
    $x\in X$ satisfies an inequality in Definition \ref{def:vi} with constant $b\leq 1$ and a modifier function $f(t) = Ct$ for some $C > 0$ if and only if
    $|\langle x, u\rangle_{X}| \leq \frac{C}{2}\|\sqrt{H}u\|_X$
    for all $u\in B_X(1)$.
    Alternatively, $x$ satisfies a variational inequality with some $b < 1$, $C>0$ and $f(t) = Ct^\eta$, $\eta\in (0,1)$ if and only if, for all $u\in B_X(1)$,
    \begin{equation*}
        |\langle x, u\rangle_{X}| \leq 
        \frac{2-\eta}{2}\left(
        \frac{1-b}{1-\eta}\right)^{\frac{1-\eta}{2-\eta}}
        C^{\frac{1}{2-\eta}}\|\sqrt{H}u\|_{X}^{\frac{\eta}{2-\eta}}.
    \end{equation*}   
\end{theorem}

To state the next result, we consider the equation
\begin{equation} \label{eq:fixed_point}
    \Phi_{x}(R) := \frac{d_{x}(R)}{R} = \sqrt{\epsilon}.
\end{equation}
As shown by \cite{hofmann2007analysis} and \cite{flemming2011sharp}, 
    $\Phi_{\lpsi} : (0,\infty) \to (0,\infty)$ is a continuous, strictly decreasing function for $x\notin R(\sqrt{H})$. 
The following result is taken from \cite[Corollary~1]{flemming2011sharp}.

\begin{theorem} \label{lem:spec_fp_bounds}
    We have 
    $d_{x}(2\Phi_{x}^{-1}(\sqrt{\epsilon})) \leq F_{\lpsi}(\epsilon) \leq 
    2d_{x}(\Phi_{x}^{-1}(\sqrt{\epsilon}))$ 
    for each $\epsilon \in (0,\infty)$.
\end{theorem}

When $F_{x}(\epsilon)$ decays polynomially in $\epsilon$ in a neighborhood of zero, then there exists $c\in(0,1]$ such that the lower bound above is bounded from below by $cd_x(\Phi_x^{-1}(\sqrt{\epsilon}))$;
see the discussion in \cite[Remark~4]{flemming2011sharp}.
\end{appendix}

\begin{acks}[Acknowledgments]
The authors would like to acknowledge the late Sayan Mukherjee, whose intellectual generosity and vision played a central role in bringing the authors of this work together. 
He had a rare ability to connect researchers from different backgrounds and inspire collaboration on challenging problems, and this project reflects that influence. 
His guidance and encouragement were invaluable, and we are deeply grateful for his mentorship. 
\end{acks}

\begin{funding}
KP acknowledges the partial support of the project PNRR - M4C2 - Investimento 1.3, Partenariato Esteso PE00000013 - “FAIR - Future Artificial Intelligence Research” - Spoke 1 “Human-centered AI”, funded by the European Commission under the NextGeneration EU programme.
\end{funding}

\begin{supplement}
\stitle{Supplementary Materials to ``Frequentist Coverage of Bayes Posteriors in
Nonlinear Inverse Problems with Gaussian
Priors''}
\sdescription{This document contains the remaining proofs and results for numerical experiments.}
\end{supplement}

\bibliographystyle{imsart-nameyear} 
\bibliography{references}       


\clearpage

\section{Remaining Proofs for Section 3}
\label{ssec:proof3}

\begin{proof}[Proof of Lemma \ref{lem:perturbation_lemma}]
The properties of the operator are immediate from the fact that the identity operator is bounded below on every dense subspace of $\ltwo(\Z)$ and thus has a bounded inverse on $D(\Lambda^s)$ with range contained in $D(\Lambda^s)$ for every $s\geq 0$.
This fact then shows that since $\lpsi = \Lambda^{-1/2}\psi \in \Hs^{1}$, the perturbation \ref{eq:perturbation} belongs to $\Hs^{2} = D(\Lambda)$.
Now the equality follows from simple algebra:
\begin{eqnarray*}
    &(N\opinfo^*\opinfo + \tau_N^{-2}\Lambda)\Lambda^{-1/2}(\opcomb^*\opcomb + \tau_N^{-2})^{-1}\Lambda^{-1/2}u \\
    &= \Lambda^{1/2}(N\opcomb^*\opcomb + \tau_N^{-2})(N\opcomb^*\opcomb + \tau_N^{-2})^{-1}\Lambda^{-1/2}u = u.
\end{eqnarray*}
\end{proof}

\subsection{Proof of Lemma \ref{lem:gaussian_limit}}

The proof is broken down into several steps, similar to the proof of BvM in \cite{monard2021statistical}.

\subsubsection{Renormalization by Posterior Concentration}

First, we define by $\nprior^{\Theta_N}(\cdot \mid \data)$ the posterior distribution $\nprior(\cdot \mid \data)$ truncated and re-normalized to the set $\Theta_N$.  
For any Borel set $A\subset\Hs^p$ the support of the prior, we have
\begin{eqnarray*}
    &\nprior(A \mid \data) - \Pi_N^{\Theta_N}(A \mid \data) \\
    &= \nprior(\Theta_N^c\cap A \mid \data) - 
    \nprior(\Theta_N^c \mid \data)
    \nprior^{\Theta_N}
    (A \mid \data),
\end{eqnarray*}
so the total variation distance between $\nprior(\cdot\mid\data)$ its restriction to $\Theta_N$ is bounded above as
\begin{equation*}
    d_{TV}(\nprior(\cdot \mid \data),\nprior^{\Theta_N}(\cdot \mid \data)) = \sup_{A\text{ Borel }} |\nprior(A \mid \data) - \nprior^{\Theta_N}(A \mid \data)|
    \leq 2\nprior(\Theta_N^c \mid \data).
\end{equation*}
By Assumption \ref{assmp:posterior_contraction}, the right hand side of the above equation is $o_{\ptheta}(1)$. 
\begin{equation} \label{eq:tv_convergence_global}
d_{\mbox{\tiny weak}}( \nprior(\cdot \mid \data) , \nprior^{\Theta_N}(\cdot \mid  \data) ) = o_{P_{\theta_0}^N}(1).
\end{equation}
It thus suffices to show the assertion of Lemma \ref{lem:gaussian_limit} holds for $\nprior^{\Theta_N}(\cdot \mid \data)$.

\subsubsection{Setting Up: Posterior Laplace Transform}
\label{ssec:sub2}

Define $Z_N := \Psi(\theta-\theta_0) - \bigpsi + b_N$.
We will show that the Laplace transform corresponding to the probability law of the random variable $s_N^{-1}Z_N$, induced by the truncated posterior measure $\nprior^{\Theta_N}(\cdot \mid \data)$:
\begin{equation} \label{eq:laplace_transform}
    \E^{\Theta_N}[e^{ts_N^{-1}Z_N}\mid \data] :=
    \frac{\int_{\Theta_N} e^{t s_N^{-1}Z_N - \loss(\theta)}~d\nprior(\theta)}
    {\int_{\Theta_N} e^{-\loss(\theta)}~d\nprior(\theta)},
\end{equation}
converges for each $t\in\R$ to $e^{t^2/2}$ in probability.
When true, by \cite[Theorem~1.13.1]{van2023weak}, the pointwise convergence in probability of the above display is equivalent to
\begin{equation*}
    d_{weak}(Law(s_N^{-1}Z_N),Normal(0,1))\stackrel{\ptheta}{\to }0.
\end{equation*}
Combined with \eqref{eq:tv_convergence_global},
this concludes the proof.

\subsubsection{Local Asymptotic Expansion}

To begin with, we will recall the definitions from \eqref{eq:lan_lik} and \eqref{eq:rem_lik} imply:
\begin{equation*}
    \loss(\theta) - \loss(\theta_0) = \llan(\theta-\theta_0) + \stoch(\theta,\theta_0);
\end{equation*}
therefore, the Laplace transform \eqref{eq:laplace_transform} is rewritten as
\begin{equation*}
    \frac{\int_{\Theta_N} e^{ts_N^{-1}Z_N - \llan(\theta-\theta_0) - \stoch(\theta,\theta_0)}~d\nprior(\theta)}{\int_{\Theta_N} e^{-\llan(\theta-\theta_0) - \stoch(\theta,\theta_0)}~d\nprior(\theta)}.
\end{equation*}
Define a perturbed vector $\thetat := \theta - ts_N^{-1}\barpsi $ for each $\theta\in\Hs^p$ and each $t\in\R$. Simple algebra shows
\begin{equation*}
    \llan(\theta-\theta_0) - \llan(\thetat - \theta_0)
    = -t\frac{\bigpsi}{s_N} - \frac{Nt^2\|\opinfo(\barpsi)\|_{\ltwo(\X)}^2}{2s_N^2} 
    + \frac{Nt\langle\opinfo(\barpsi),\opinfo(\theta-\theta_0)\rangle_{\ltwo(\X)}}{2s_N}.
\end{equation*}
Plugging in the definition of $ts_N^{-1}Z_N$, we observe the integrand inside the numerator of \eqref{eq:laplace_transform} can be rewritten as
\begin{eqnarray*}
    &\exp\bigg\{
    ts_N^{-1}b_N + 
    \frac{Nt^2}{2s_N^2}\|\opinfo\barpsi\|^2_{\ltwo(\Z)} +
    ts_N^{-1}\langle
    \psi -
    N\opinfo^*\opinfo\barpsi,
    \theta - \theta_0\rangle_{\ltwo(\Z)} -\\
    &\llan(\thetat-\theta_0) -\stoch(\theta,\theta_0)
    \bigg\}.
\end{eqnarray*}

\subsubsection{Change of Measure by an Absolutely Continuous Shift}

A basic consequence from Lemma \ref{lem:perturbation_lemma}, along with the fact that $\barpsi \in D(\Lambda)$, is the equality
\begin{equation*}
    \psi - N\opinfo^*\opinfo\barpsi = \tau_N^{-2}\Lambda\barpsi,
\end{equation*}
Write by $\Pi_N'$ the pushforward of the prior $\nprior$ under the shift map $h\mapsto h + ts_N^{-1}\barpsi$.
Since, as observed in Lemma \ref{lem:perturbation_lemma}, the shift lies in $D(\Lambda^{1/2})$,
by the Cameron-Martin theorem, we obtain
\begin{equation} \label{eq:cameron_martin_app}
    e^{t s_N^{-1}\langle \theta , \tau_N^{-2}
    \Lambda(\barpsi) \rangle_{\ltwo(\X)}} = 
    \frac{d\nprior'}{d\nprior}(\theta) \times
    \exp\left(t^2 \frac{\|\tau_N^{-1}\Lambda^{1/2}(\barpsi)\|_{\ltwo(\X)}^2}{2s_N^2}\right).
\end{equation}
Therefore, by change of measure and recalling definitions of $b_N$ and $s_N$, the numerator may be rewritten as
\begin{equation} \label{eq:change_of_measure}
    e^{t^2/2}\int_{\Theta_N^{(t)}} e^{-\stoch(\theta^{(-t)},\theta_0)}e^{-\llan(\theta,\theta_0)}~d\nprior(\theta),
\end{equation}
where $\Theta_N^{(t)} = \Theta_N - ts_N^{-1}\barpsi$.
Using definition  \eqref{eq:posterior_lan}, the Laplace transform \eqref{eq:laplace_transform} is simplified to
\begin{equation} \label{eq:laplace_simplified}
    e^{t^2/2}\times \frac{\int_{\Theta_N^{(t)}} e^{-\stoch(\theta^{(-t)},\theta_0)}~d\nlan(\theta|\data)}{\int_{\Theta_N} e^{-\stoch(\theta,\theta_0)}~d\nlan(\theta|\data)}.
\end{equation}
We isolate the second multiplicative term that is a ratio of two integrals with respect to $\nlan(\cdot|\data)$. It may rewritten as
\begin{equation*}
    1 + \frac{\int_{\Theta_N^{(t)}} e^{-\stoch(\theta^{(-t)},\theta_0)}~d\nlan(\theta|\data) - 
    \int_{\Theta_N} e^{-\stoch(\theta,\theta_0)}~d\nlan(\theta|\data)}{\int_{\Theta_N} e^{-\stoch(\theta,\theta_0)}~d\nlan(\theta|\data)},
\end{equation*}
where the second summand, by taking a uniform bound, is absolutely bounded above by
\begin{equation*}
    e^{2\sup_{\theta\in\Theta_N}|\stoch(\theta,\theta_0)|}
    \left|1 - \frac{\nlan\{\Theta_N-ts_N^{-1}\barpsi|\data\}}{\nlan\{\Theta_N|\data\}}\right|.
\end{equation*}
By conditions \eqref{eq:stoch_unif_control} and \eqref{eq:lan_ratio_control}, this term is $o_{\ptheta}(1)$. 

\subsection{Proof of Lemma \ref{lem:posterior_mean_conv}}

The proof is again broken down into several parts. 
Parts of our proof closely follow the arguments of the proof for \cite[Theorem~3.8]{monard2021statistical}.
The basic idea is to strengthen the weak convergence of the posterior (in probability) to a Gaussian limit in Lemma \ref{lem:gaussian_limit} by showing the sequence of posterior random variables is also uniformly integrable.
The scaled posterior mean should then converge to the $s_N^{-1}\bigpsi$, which itself is a random variable.
Some additional arguments will be needed to show that $s_N$ and $t_N$ are orderwise equivalent, and to derive the limiting distribution of $t_N^{-1}\bigpsi$.

\subsubsection{Controlling the Sequence Ratio}

We first state a lemma that complements the previous upper bound statement in Lemma \ref{lem:st_order_equivalence} by showing $t_N/s_N$ is also bounded \emph{from below} away from zero.

\begin{lemma}
    \label{lem:st_bound_below}
    If $\psi\notin{\rm ker}(\opinfo)$, then $\liminf_{N\to\infty}t_N/s_N > 0$.
\end{lemma}

\begin{proof}
    We proceed similarly as in the proof of \cite[Theorem~5.4]{knapik2011bayesian}.
    Write $s_N^2 = \int_0^\infty f_N^2(\opcomb^*\opcomb)\\~dF_{\lpsi}^2(t)$ for $f_N(t) = \frac{1}{\sqrt{1+N\tau_N^2t}}$,
    and write $s_N^2 = A_N + B_N$ with
    \begin{equation*}
    A_N  =  \int_{\frac{1}{N\tau_N^2}}^\infty f_N(t)~dF_{\lpsi}^2(t),\;
    B_N  =  \int_0^{\frac{1}{N\tau_N^2}}f_N^2(t)~dF_{\lpsi}^2(t).
    \end{equation*}
    Similarly, write $t_N^2 = \int_0^\infty g_N^2(\opcomb^*\opcomb)~dF_{\lpsi}^2(t)$ for $g_N(t) = \frac{\sqrt{N}\tau_N\sqrt{t}}{\sqrt{1+N\tau_N^2t}}$,
    and $t_N^2 = C_N + D_N$ with 
    \begin{equation*}
    C_N  =  \int_{\frac{1}{N\tau_N^2}}^\infty g_N(t)~dF_{\lpsi}^2(t),\;
    D_N  =  \int_0^{\frac{1}{N\tau_N^2}} g_N^2(t)~dF_{\lpsi}^2(t).
    \end{equation*}
    We claim that it suffices to prove $\limsup_{N\to\infty}B_N/A_N < \infty$.
    This is because 
    \begin{equation*}
        \frac{t_N^2}{s_N^2} = \frac{C_N + D_N}{A_N + B_N} \geq \frac{C_N}{A_N + B_N} = \frac{C_N/A_N}{1 + B_N/A_N}.
    \end{equation*}
    It is straightforward to check that for each $N$, $C_N/A_N\geq\frac{1}{2}$, so whenever $\limsup_{N\to\infty}B_N/A_N$ is finite, 
    we have $\liminf_{N\to\infty} t_N^2/s_N^2 > 0$.

    We now prove if $\psi\notin{\rm ker}(\opinfo)$, then $\limsup_{N\to\infty}B_N/A_N<\infty$.
    It suffices to prove that $\limsup_{N\to\infty} \\ B_N/A_N < \infty$ if there exists some $L\in(0,1)$ such that for all large enough $N$,
    \begin{equation*}
        A_N  =  \int_{\frac{1}{N\tau_N^2}}^\infty f_N^2(t)~dF_{\lpsi}^2(t) \geq L
        \int_0^\infty f_N^2(t)~dF_{\lpsi}^2(t).
    \end{equation*}
    Furthermore, the above display is true if there exists some $c>0$ such that $A_N\geq cB_N$; for then we have $A_N \geq \frac{c}{c+1}(A_N+B_N)$.
    We now prove this sufficient condition, and to do so, we suppose otherwise:
    for every $c>0$ and positive integer $N_0$, there exists an $N\geq N_0$ such that $A_N < cB_N$.
    For this to be true on a given $N$, since both $A_N$ and $B_N$ are positive quantities bounded from above by 1, it is necessary that $A_N = 0$.
    Write by $(\sigma_n^2,\nu_n)$ the eigenvalues and eigenvectors of $\opcomb^*\opcomb$.
    By definition of $A_N$, this can only be the case if for every $n>0$, either $\sigma_n^2 < \frac{1}{N\tau_N^2}$ or $\langle\lpsi,\nu_n\rangle=0$ (or both).
    Extract a subsequence of positive integers $N_k$ such that for every $N_k$ we have $A_{N_k} = 0$.
    Then $N_k$ diverges to $+\infty$ by the contradiction hypothesis, so we conclude that for a non-trivial operator $\opcomb^*\opcomb$, we must have $\langle\lpsi,\nu_n\rangle_{\ltwo(\Z)}=0$ for every $n$ such that $\sigma_n^2 > 0$.
    Since we have assumed $\psi\notin {\rm ker}(\opinfo)$, which is if and only if $\lpsi\notin{\rm ker}(\opcomb^*\opcomb)$ by injectivity of $\Lambda^{-1/2}$, we obtain a contradiction.
\end{proof}

\subsubsection{Triangular Central Limit of the Score Term}

Next, we derive a limiting law for the sequence of random variables $t_N^{-1}\bigpsi$ by the
Lindeberg-Feller central limit theorem \cite[Theorem~9.8.1]{resnick2019probability}. 
Recalling definitions \eqref{eq:seq_lf} and \eqref{eq:bias}, the sequence of random variables $\{t_N^{-1}\epsilon_i\opinfo(\barpsi)(X_i)\}_{i,N}$ are independent, and the subsequence $\{t_N^{-1}\epsilon_i\opinfo(\barpsi)(X_i)\}_{i=1}^{N}$ for each $N$ comprises i.i.d. zero-mean random variables such that 
\begin{equation*}
    \sigma_{N,i}^2 := \E[t_N^{-2}\epsilon_i^2\opinfo(\barpsi)(X_i)^2]
    = \frac{1}{N}.
\end{equation*}
We now verify the Lindeberg's condition for the application of the central limit theorem as in \cite[(9.23)]{resnick2019probability}.
From the definition of perturbation \eqref{eq:asymp_scale}, we have $\lpsi \in B_{\Hs^{2}}(M)$ for some $M>0$. 
From Assumption 4, we can choose some dense subspace $\Hs^p\subset\ltwo(\Z)$ over which $\Pi_0(\Hs^p) = 1$. Then,
since the Cameron-Martin space $\Hs^1$ of $\Pi_0$ is compactly embedded in $\Hs^p$  (cf. \cite[Theorem~6.11]{stuart2010inverse}), we may conclude that $\|\opinfo(\barpsi)(x)\|_{\infty} \leq M'$ for some 
$M' =M'(\theta_0) > 0$. Now, by Cauchy-Schwarz inequality and using the fact that $\E[\epsilon_i^4] = 3$, for every $r > 0$ we have
\begin{eqnarray*}
   & t_N^{-2} \sum_{i=1}^{N}     
    \E\left[\epsilon_i^2\opinfo(\barpsi)(X_i)^2
    \mathbf{1}\left(
    |t_N^{-1}\epsilon_i\opinfo(\barpsi)(X_i)| > r
    \right) \right]\\
   & \leq \sqrt{3}
    {M'}^2
    \max_{1\leq i\leq N}
    \sqrt{Pr(t_N^{-1}
    \epsilon_i\opinfo(\barpsi)(X_i) > r)}.
\end{eqnarray*}
The latter probability is identical for all $i$ and is bounded above by $r/N$ due to Chebyshev's inequality.
This concludes the verification of the conditions for the central limit, and we conclude that $t_N^{-1}\bigpsi$ weakly converges to $Normal(0,1)$.

\subsubsection{Uniform Tightness Argument}

We proceed to show that $s_N^{-1}\{\Psi(\theta-\theta_0) - \bigpsi\}$ is uniformly tight, using the arguments from the proof of \cite[Theorem~3.8]{monard2021statistical}.
Before doing so, we claim that indeed this is enough to also show that
$t_N^{-1}\{\Psi(\bartheta - \theta_0) - \bigpsi\}$ is uniformly tight.
This is because by Lemma \ref{lem:st_bound_below} and then a conditional version of Jensen's inequality, we have
\begin{equation*}
    t_N^{-2}\{\Psi(\bartheta-\theta_0)\}^2 \lesssim s_N^{-2}\E[\Psi(\theta-\theta_0)\mid\data]^2\leq s_N^{-2}\E[\{\Psi(\theta-\theta_0)\}^2\mid\data].
\end{equation*}

To show the uniform tightness of $s_N^{-1}\{\Psi(\theta-\theta_0) - \bigpsi\}$, we first note that
\begin{equation*}
    s_N^{-2}\E[\{
    \Psi(\theta-\theta_0)-\bigpsi\}^2 \mid \data] \lesssim 
    2s_N^{-2}
    \E[\Psi^2(\theta - \theta_0) \mid \data] + 2s_N^{-2}\bigpsi^2,
\end{equation*}
where the latter term is $O_{\ptheta}(1)$, since $t_N^{-1}\bigpsi$ is, by our previous triangular central limit theorem.
The former term involves an expectation, which we may write as $\E[\Psi^2(\theta-\theta_0)\mathbf{1}_{\Theta_N}\mid\data] + \E[\Psi^2(\theta-\theta_0)\mathbf{1}_{\Theta_N^c}\mid\data]$,
for a sequence of sets $\Theta_N$ from Assumption 
\ref{assmp:posterior_contraction}. Since $t^2 \leq e^t$ for all $t\geq 0$, we have
(recall the notation $Z_N := \Psi(\theta - \theta_0) - \bigpsi + b_N$):
\begin{equation}
\label{eq:restricted_fnl}
    s_N^{-2}\E[\Psi^2(\theta-\theta_0)\mathbf{1}_{\Theta_N}\mid \data] \leq \E[e^{s_N^{-1}(Z_N +\bigpsi - b_N)}\mathbf{1}_{\Theta_N} \mid\data]
    \vee
    \E[e^{-s_N^{-1}(Z_N + \bigpsi - b_N)}\mathbf{1}_{\Theta_N}\mid\data].
\end{equation}
Since $s_N^{-1}b_N = o_{\ptheta}(1)$ and $t_N \leq s_N$ for every $N$ (Lemma \ref{lem:st_order_equivalence}), 
we also have $s_N^{-1}\bigpsi = O_{\ptheta}(1)$.
Moreover, $\E[e^{s_N^{-1}Z_N}\mathbf{1}_{\Theta_N}\mid\data]$ and $\E[e^{-s_N^{-1}Z_N}\mathbf{1}_{\Theta_N}\mid\data]$
coincide with the Laplace transform of the posterior restricted to $\Theta_N$ in \eqref{eq:laplace_transform}, with choices $t = 1$ and $-1$, respectively. 
The conclusion of the proof of Lemma \ref{lem:gaussian_limit} showed that both these two expectations converge to $e^{1/2}$ in probability, so the left hand side of \eqref{eq:restricted_fnl} is indeed $O_{\ptheta}(1)$.

It remains to handle $s_N^{-2}\Psi(\theta-\theta_0)\mathbf{1}_{\Theta_N^c}$ for the proof of uniform tightness. 
First, we observe that $s_N\gtrsim 1/\sqrt{N}$ in all three cases given in Lemma \ref{lem:order_s}, so that by Cauchy-Schwarz inequality,
\begin{equation*}
    s_N^{-2}\E[\Psi(\theta-\theta_0)^2\mathbf{1}_{\Theta_N^c} \mid \data] \lesssim 
    N^2\sqrt{\E[\Psi(\theta-\theta_0)^4 \mid \data]}\sqrt{\nprior(\Theta_N^c \mid \data)}.
\end{equation*}
We bound from below the denominator of the posterior integral $\E[\Psi(\theta - \theta_0)^4 \mid \data]$
using \cite[Lemma~1.3.3]{nickl2023bayesian}.
Based on Assumption \ref{assmp:unif_bd_continuity}, fix some $M>0$ and envelope constant $U = U(M) > 0$ and define the set $\mathcal{A}_N := \{ \theta : \|G(\theta)-G(\theta_0)\|_{\ltwo(\X)} \leq 
    \epsilon_N,\; \|G(\theta)\|_{\infty} \leq U\}$.
Then under Assumption \ref{assmp:prior_mass}, we can apply the Lemma and obtain an in-probability bound
\begin{equation*}
    \ptheta\left( \int_{\mathcal{A}_N} e^{-\loss(\theta) + \loss(\theta_0)}~d\nprior(\theta) \geq e^{-KN\epsilon_N^2}\right) \to 1\quad\forall K\geq 2.
\end{equation*}
We fix $K = A+2$ for a constant $A$ from Assumption \ref{assmp:prior_mass}. Now, using furthermore Assumption \ref{assmp:posterior_contraction}, we obtain the following:
\begin{eqnarray*}
    &\ptheta\left(N^2\E[\Psi(\theta-\theta_0)^4 \mid \data]
    \nprior(\Theta_N^c \mid \data) > \frac{1}{\log^2 N}\right)\\
    &\leq \ptheta\left(\E[\Psi(\theta-\theta_0)^4 \mid \data]
    e^{-bN\epsilon_N^2} > \frac{1}{N^2\log^2 N}\right) + o_N(1) \\
    &\leq \ptheta\left( e^{(A+2) N\epsilon_N^2}
    \int \langle\psi,\theta-\theta_0\rangle_{\ltwo(\Z)}^4
    e^{-\loss(\theta) + \loss(\theta_0)}~d\nprior(\theta)
    > \frac{e^{bN\epsilon_N^2}}{N^2\log^2 N}\right) + o_N(1) \\
    &\leq N^2e^{-(b-A-2)N\epsilon_N^2}\log^2 N
    \iint \langle\psi,\theta - \theta_0\rangle^4_{\ltwo(\Z)} 
    e^{-\loss(\theta)+\loss(\theta_0)}~d\nprior(\theta \mid \data) d\ptheta(\data)
    + o_N(1).
\end{eqnarray*}
Now $b-A-2 > 0$ by Assumption \ref{assmp:posterior_contraction}, 
and by Fubini's theorem and Cauchy-Schwarz inequality, the latter integral is bounded above by a constant multiple of
\begin{equation*}
    \|\psi\|_{\ltwo(\Z)}^4
    \int\|\theta-\theta_0\|_{\ltwo(\Z)}^4~d\nprior(\theta),
\end{equation*}
and the latter integral is finite by Fernique's theorem (e.g., \cite[Theorem~6.9]{stuart2010inverse}).
The above argument shows that $s_N^{-2}\E[\Psi(\theta-\theta_0)^2\mathbf{1}_{\Theta_N^c}\mid \data] = o_{\ptheta}(1)$.

\subsubsection{Conclusion}

The rest of the proof now follows the same argument as the steps of the proof given by 
\cite{monard2021statistical} in p.3277, but we include them for completeness. We show through a contradiction argument that 
\begin{equation}
    \label{eq:mean_score_coupling}
    s_N^{-1}\{ \Psi(\bartheta - \theta_0) - \bigpsi\} \stackrel{\ptheta}{\to} 0.
\end{equation}
If the above is true, then the desired two weak convergence claims of Lemma \ref{lem:posterior_mean_conv} follow:
first \eqref{eq:pm_centering} from Lemma \ref{lem:gaussian_limit} and Slutsky's lemma,
and second \eqref{eq:limit_mean} from multiplying the above display by $s_N/t_N$ and then using the result from a previous section, that $t_N^{-1}\bigpsi$ weakly converges to $Normal(0,1)$.

For contradiction, suppose that the \eqref{eq:mean_score_coupling} is not true.
Write by $(\Omega,\mathcal{S},\Pr) := ((\R\times\X)^{\infty},\mathcal{S})$ the measurable space with Borel $\sigma$-algebra induced by the sequence of random variables $(X_i,\;\epsilon_i)_{i=1}^{\infty}$.
Since we assume the above limit does not hold, there exists a measurable event $E\in\mathcal{S}$ with $\Pr(E) = \tau > 0$ and $\zeta > 0$ such that along a subsequence, still denoted $N$:
\begin{equation*}
    |s_N^{-1}\{ \langle\E[\theta -\theta_0 \mid \data(\omega)],\psi\rangle_{\ltwo(\Z)} - 
    \bigpsi(\omega) \}|\geq \zeta > 0\quad\forall\omega\in E.
\end{equation*}
Recall that by Lemma \ref{lem:gaussian_limit} and by the fact that we assumed $|b_N|/s_N\ll 1$, 
\begin{equation*}
    d_{weak}(Law(s_N^{-1}\{
    \langle\psi,\theta-\theta_0\rangle_{\ltwo(Z)} - \bigpsi\}\mid \data),\;
    Normal(0,1)) \stackrel{\ptheta}{\to} 0.
\end{equation*}
The above statement remains true in an almost sure sense along some subsequence of $\data$.
Let $E'\subset\Omega$ denote this event, with $\Pr(E') = 1$ and such that the weak convergence holds for every $\omega\in E'$.
Write by $X_N = s_N^{-1}\{\langle\psi,\theta-\theta_0)-\bigpsi\}$ and $X\sim Normal(0,1)$.
For every $\omega\in E'$, by Skorohod imbedding \cite[Theorem~8.3.2]{resnick2019probability}, 
there exists a probability space
supporting real-valued random variables $\{\tilde{X}_N\}_{N=1}^{\infty},\; \tilde{X}$ which strengthen weak convergence to almost sure convergence:
\begin{equation*}
    Law(\tilde{X}_N) = 
    Law(X_N),\;
    Law(\tilde{X}) = Law(X),\;
    \tilde{X}_N(\omega) \to
    \tilde{X}(\omega)~\forall\omega\in E'.
\end{equation*}
From the previous section, we know that $\E[X_N^2\mid\data]$ is uniformly tight.
Thus there exists an event $\tilde{E}\subset E'$  with $\Pr(\tilde{E}) > 1 - \tau$
such that for every $\omega\in \tilde{E}$, we have
$\E[\tilde{X}_N^2\mid\data(\omega)] = O_N(1)$. 
The choice of $\tau$ here can be made arbitrarily close to 0. 
This implies that the sequence $\tilde{X}_N$ is uniformly integrable on the event $\tilde{E}$, and therefore 
\begin{equation*}
    \E[|X_N - X|\mid \data(\omega)]
    = \E|\tilde{X}_N(\omega) - \tilde{X}(\omega)| \to 0
\end{equation*}
for every $\omega\in \tilde{E}$. 
Therefore, the random variable $\E[X_N\mid\data(\omega)]$ converges to 0 in this event, and thus
\begin{equation*}
    t_N^{-1}\{
    \langle\E[\theta\mid\data(\omega)]-\theta_0,\psi\rangle_{\ltwo(\Z)}-\bigpsi(\omega)\}\to 0~\forall\omega\in \tilde{E}.
\end{equation*}
However, since $\Pr(\tilde{E}) > 1 - \tau$ and since we have assumed the existence of an event $E$ of probability $\tau$ over which the above convergence cannot hold,
we must have $\Pr(\Omega) \geq \Pr(\tilde{E}) + \Pr(E) > 1-\tau + \tau = 1$, 
a contradiction. 
Therefore, we conclude that the convergence statement \eqref{eq:mean_score_coupling} must be true.

\section{Proofs for Section 4}
\label{ssec:proof4}

\begin{proof}[Proof of Lemma \ref{lem:higher_cm_embedding}]
    The covariance operator $\mathcal{C}$ of $\Pi_0' = \tilde{\Pi}_0\circ M_\chi$, with $\tilde{\Pi}_0$ the measure corresponding to a Gaussian process with stationary kernel \eqref{eq:whittle_matern}, is given by the convolution $\mathcal{C}g(s) = \int_{K'} \chi(s)\mathcal{K}(t-s)\chi(t)g(t)~dt$ for $g\in\ltwo(\Z)$.
    Restrict the domain of two elements $f,g\in\ltwo(\Z)$ to $K'$ (and set them to be 0 outside $K'$).
    Since $\chi$ is chosen to only vanish at $\partial K'$, which has Lebesgue measure zero, $\mathcal{C}(f-g) = 0$ implies $f=g$ almost everywhere in $\Rd$, so ${\rm ker}(\mathcal{C})$ is trivial in $\ltwo(K')$.
    Next, for $x\in D(\Lambda)$, the multiplication operator $M_\chi$ is boundedly invertible, and we obtain using Fourier transforms
    \begin{equation*}
        (\Lambda x)(t) = \chi^{-1}(t)\frac{1}{(2\pi)^d}\int_{K'}\int_{\Rd}
        e^{\sqrt{-1}\langle s-t,\xi\rangle_{\Rd}}(1+\|\xi\|_{\Rd})^{\alpha}x(s)~d\xi~ds,\;
        \Lambda x\in \ltwo(\Z);
    \end{equation*}
    for $x\in D(\Lambda^{1/2})$, this is still valid if interpreted appropriately, 
    with $\Lambda x\in D(\Lambda^{-1/2}$),
    and we obtain $\|\Lambda^{1/2}x\|_{\ltwo(\Z)}\geq\|x\|_{\ltwo(\Z)}$ as claimed.
    Finally, the above shows that for every positive integer $k$,
    \begin{equation*}
        \|\Lambda^{k/2} x\|_{\ltwo(\Z)} \geq
        \left\{\inf_{x\in K'}\chi(x)\right\}^{k}\|M^{(k)}\mathcal{F}(x)\|_{\ltwo(\Z)}\geq 
        \|M^{(k)}\mathcal{F}( x)\|_{\ltwo(\Z)},
    \end{equation*}
    where $\mathcal{F}$ is the Fourier transform operator and $M^{(k)}\hat{u}(\xi) := (1+\|\xi\|_{\Rd}^{\alpha})^{k/2}\hat{u}(\xi)$ for every $\hat{u}\in\ltwo(K')$.
    Since $(1+\|\xi\|_{\Rd}^{\alpha})^k\geq 1+\|\xi\|_{\Rd}^{k\alpha}$ for every positive integer $k$, by the Fourier transform characterization of Sobolev norms \cite[Chapter~4,~Eq.(1.2)]{taylor2013partial},
    we conclude $\|\Lambda^{k/2}x\|_{\ltwo(\Z)}\gtrsim\|x\|_{H^{k\alpha}(\Z)}$ for every such $k$, proving the desired continuous embedding result for positive integer $s$.
    For general $s\in\R$, we note that \cite{neubauer1988sobolev} shows that while there is no operator generating Hilbert scales that coincide with Sobolev spaces on a domain, there does exist an operator generating Hilbert scales which agree with Sobolev spaces up to an arbitrary order $p>0$ (cf. Corollary 2.3 and remark on p.562; see \cite[Example~6.3]{gugushvili2020bayesian} for a constructive example).
    For an arbitrary positive integer $k\geq1$, consider Hilbert scales $(H^*_s)_{s\geq 0}$ generated by an operator $L_*$ which coincide with Sobolev spaces $H^s(\Z)$ for every $s\leq k+1$. 
    What we have just shown is that $D(\Lambda^{k/2})\subset D(L_*^{k/2})$, and that viewed as two Hilbert spaces, the subset relation is a continuous embedding.
    Since both $\Lambda^{1/2}$ and $L_*^{1/2}$ are closed, it follows that $L_*^{k/2}$ is $\Lambda^{k/2}$-bounded \cite[Chapter~IV,~Remark~1.5]{kato1966perturbation}: i.e., there exists some $C>0$ such that
    $\|L_*^{k/2}x\|_{\ltwo(\Z)} \leq C^2(\|\Lambda^{k/2}x\|_{\ltwo(\Z)} + \|x\|_{\ltwo(\Z)})$ for every $x\in D(\Lambda^{k/2})$. Since $t\mapsto t^{r/2}$ is an operator monotone function for every $r\in(0,1)$, we also have
    $\|L_*^{kr/2}x\|_{\ltwo(\Z)} \leq C^2(\|\Lambda^{kr/2}x\|_{\ltwo(\Z)} + \|x\|_{\ltwo(\Z)})$ for every such $r$,
    which relation proves that $D(\Lambda^{kr/2})$ is continuously embedded into $D(L_*^{k/2})$.
    Since this argument can be reiterated for each integer $k$, we conclude that in general $\Hs^s$ is continuously embedded in $H^{s\alpha}(\Z)$.
\end{proof}

\subsection{Proof of Theorem \ref{thm:darcy_coverage}}

We verify the asserted claim, that the stated assumptions of Theorem \ref{thm:darcy_coverage} are sufficient for Assumptions \ref{assmp:truth_prior}--\ref{assmp:deriv_stability},
conditions \eqref{eq:stoch_unif_control}--\eqref{eq:lan_ratio_control}, and that $|b_N|/s_N\ll1$.
For invoking results in Appendix \ref{appendix:vi_lemmas}, we note that $\psi$ is by assumption not trivial when restricted to $K'$. By Lemma \ref{lem:higher_cm_embedding}, we know that $\Lambda$ is densely defined in $\ltwo(K')$. so it suffices to view $\mathcal{C},\Lambda$ as being defined on $\ltwo(K')$, and use the fact that on a compact manifold $K'$, we have equivalence of norms between $H^s(K')$ and $H^s(\Z)$ for arbitrary $s\geq 0$.

\noindent {\bf 1. Assumption \ref{assmp:truth_prior}}.
We have already assumed that the prior is Gaussian, that $(Y_i,X_i)\stackrel{iid}{\sim} P_{\theta_0}$ for some $\theta_0\in\ltwo(\Z)$, and that
$\theta_0\in\Hs^\beta$.

\noindent {\bf 2. Assumption \ref{assmp:unif_bd_continuity}}.
Recall that $\alpha > 1 + d$. Choose an $s\in (1+d/2,\alpha-d/2)$ so that $H^s(\Z)$ is a Hilbert support of $\Pi_0$.
The assumption can now be verified based on Lemmas \ref{lem_supp2} and \ref{lem_supp3}.

\noindent {\bf 3. Assumption \ref{assmp:derivative}}. This assumption is verified by explicitly deriving the requisite linearization operator $\opinfo$ in Section \ref{ssec:darcy_lemmas} below.
In brief, there exists a continuous linear extension of $\mathbb{I}_\theta : H^s(\Z)\to\ltwo(\Z)$, with $s>1+d/2$, given by $-L_\theta^{-1}[\nabla\cdot\{(\phi'\circ\theta)h\nabla u_\theta\}]$, with the remainder term $\|R_\theta(h)\|_{\ltwo(\Z)} = \|G(\theta+h)-G(\theta)-\mathbb{I}_\theta(h)\|_{\ltwo(\Z)} =  O(\|h\|_\infty^2)$.
Lemmas \ref{lem_supp4} and \ref{lem_supp5} derive a Lipschitz continuity estimate of $\opinfo$ when considered as a mapping from $H^{-1}(\Z)\to\ltwo(\Z)$ and a mapping from continuously differentiable functions $C_1(\Z)\to C(\Z)$, respectively.

\noindent {\bf 4. Assumption \ref{assmp:prior_mass}}.
The choice of  $\tau_N$ clearly satisfies the condition $\frac{1}{\sqrt{N}}\ll \tau_N\lesssim 1$. 
The prior mass condition was verified in \cite[Theorem~4]{giordano2020consistency} with the choice of a sequence $\epsilon_N = N^{-\frac{\alpha+1}{2\alpha+2+d}}$.

\noindent {\bf 5. Assumption \ref{assmp:posterior_contraction}}.
\cite[Theorem~5]{giordano2020consistency} showed this is true for the above choice of $\epsilon_N$ and the pre-defined sets $\Theta_N$ \eqref{eq:lsupset}
valid for all $p \in (1 + d/2,\alpha - d/2)$ when $\alpha > 1 + d$.
Clearly, we have $\sup_{\theta\in\Theta_N}|\Psi(\theta-\theta_0)| \ll 1$ for any continuous linear functional $\Psi$ on $\ltwo(\Z)$.

\noindent {\bf 6. Assumption \ref{assmp:deriv_stability}}.
This assumption is verified through the following lemma. Note that our prior, by construction of \cite[Example~25]{giordano2020consistency}, satisfies the continuous embedding of $\Hs^1$ into $H^\alpha_0(\Z)$, with $\alpha > 1$.

\begin{lemma} \label{lem:darcy_regularity}
    Suppose $\alpha > 1$ and $h\in \Hs^s$, with $\Hs^1$ continuously embedded in $H^\alpha_0(\Z)$, and each $\Hs^s$ in $H^{s\alpha}(\Z)$ for $s>1$.
    Under Assumption \ref{assmp:darcy_flow},
    Assumption \ref{assmp:deriv_stability} is satisfied for the linear operator $\opinfo$ at $\theta_0\in B_{\Hs^p}(M)$ for some $p>1+d/2$, with $\delta(s) = \frac{s\alpha-1}{s\alpha+1}$.
\end{lemma}

\begin{proof}
By \cite[Theorem~5]{nickl2022some}, if $\theta_0\in C^\infty(\Z)$, then we have
\begin{equation*}
    C\|\opinfo(h)\|_{H^2(\Z)} \geq \|h\|_{\ltwo}\quad\forall h\in H_0^1(\Z)
\end{equation*}
with constant $C = C(k_1,k_2,\theta_0,\Z)$, where $k_1,k_2$ are constants in Assumption \ref{assmp:darcy_flow}, 
and the derivative of the forward operator $G$ is explicitly derived in Section \ref{ssec:darcy_lemmas}.
The sufficient condition that $\theta_0\in C^\infty(\Z)$ can be weakened, following the proof of \cite[Proposition~2.1.7.]{nickl2023bayesian}, where a similar is derived for $\theta_0\in H^p(\Z),\; p > 1 + d/2$.
Let $h\in B_{H_0^p}(M)$ for $p>1$; then, by interpolation between Sobolev scales,
\begin{equation*}
    \|\opinfo(h)\|_{H^2(\Z)} \lesssim \|\opinfo(h)\|^{\frac{p - 1}{p + 1}}_{\ltwo(\Z)}
    \|\opinfo(h)\|^{\frac{2}{p + 1}}_{H^{p + 1}(\Z)},
\end{equation*}
with the latter $H^{p+1}(\Z)$-norm bounded uniformly in $h\in B_{H^p}(M)$, by Lemma \ref{lem_supp1}.
Choosing $p = s\alpha$ and using Lemma \ref{lem:higher_cm_embedding}, we conclude that Assumption \ref{assmp:deriv_stability} holds with the given $\delta(s)$-exponent when $h\in\Hs^{s\alpha}$.
\end{proof}

\noindent {\bf 7. Verifying $|b_N|/s_N\ll1 $ by Theorem \ref{thm:order_estimates}.}
First, that $\opinfo$ for $\theta_0\in\Hs^s,\;s\geq 1$ is injective on $\ltwo(\Z)$ is shown by \cite[Proposition~5]{nickl2022some}, so $\psi\neq 0$ does not belong to the kernel of $\opinfo$.
Some algebra based on the estimate in Lemma \ref{lem:darcy_regularity} above shows that the maintained condition on $\beta$ in \eqref{eq:complex_condn3} of Theorem \ref{thm:darcy_coverage} implies Theorem \ref{thm:order_estimates} (ii) under the choice $\tau_N \asymp (\sqrt{N}\epsilon_N)^{-1}$;
thus, we have $|b_N|/s_N\ll 1$.

\noindent {\bf 8. Showing \eqref{eq:stoch_unif_control}.}
We now focus on showing that two conditions, given in conditions \eqref{eq:stoch_unif_control} and \eqref{eq:lan_ratio_control} of Lemma \ref{lem:gaussian_limit}, are satisfied under the maintained assumptions.
We recall by definition
\begin{eqnarray*}
   &\stoch(\theta,\theta_0) :=
    -\sum_{i=1}^{N}\epsilon_i\{
    G(\theta) - G(\theta_0) - \opinfo(\theta-\theta_0)\}(X_i)
    +\\
   &\quad \frac{1}{2}\sum_{i=1}^{N}\left[
    \{G(\theta)(X_i) - G(\theta_0)(X_i) \}^2 - \|\opinfo(\theta-\theta_0)\|_{\ltwo(\Z)}^2\right].
\end{eqnarray*}
We now claim the following are true:
\begin{eqnarray}
    \E\sup_{\theta\in\Theta_N} \left|\sum_{i=1}^{N}
    \epsilon_i\{
    G(\theta) - G(\theta_0) - \opinfo(\theta-\theta_0)\}(X_i)\right| &= o_N(N\epsilon_N^2);
    \label{eq:empirical_bound1} \\
    \E\sup_{\theta\in\Theta_N} \left|
    \sum_{i=1}^{N} \{G(\theta) - G(\theta_0)\}^2(X_i) - 
    N\|G(\theta) - G(\theta_0)\|_{\ltwo(\Z)}^2
    \right|
    &= o_N(N\epsilon_N^2);
    \label{eq:empirical_bound2}
\end{eqnarray}
and, moreover, pointwise for every $t\in \R$:
\begin{equation}
    \E \sup_{\theta\in \Theta_N^{(t)}}
    \left|\sum_{i=1}^{N}
    \{\opinfo(\theta - \theta_0)(X_i)\}^2 -
    N\| \opinfo(\theta - \theta_0) \|_{\ltwo(\Z)}^2
    \right| = o_N(N\epsilon_N^2);
    \label{eq:empirical_bound3}
\end{equation}
where $\Theta_N^{(t)} = \Theta_N - ts_N^{-1}\barpsi$ for each $t$.

Given these claims, which we will prove in Lemma \ref{lem:empirical_process_bounds}, 
we now observe that 
\begin{equation*}
\sup_{\theta\in\Theta_N} 
|\stoch(\theta,\theta_0)| = 
N\left|\|G(\theta) - G(\theta_0)\|_{\ltwo(\Z)}^2 - \|\opinfo(\theta - \theta_0)\|_{\ltwo(\Z)}^2 \right| + 
o_{\ptheta}(N\epsilon_N^2).
\end{equation*}
By Cauchy-Schwarz inequality and since $\opinfo$ satisfies \eqref{eq:quadratic_approx}, we have
\begin{equation*}
    N|\|G(\theta) - G(\theta_0)\|_{\ltwo(\Z)}^2 - \|\opinfo(\theta - \theta_0)\|_{\ltwo(\Z)}^2|
    \leq 2N(\|\theta-\theta_0\|_{\infty}^4 + \|\theta-\theta_0\|_{\infty}^2
    \|\theta-\theta_0\|_{C^1}).
\end{equation*}
We claim 
\begin{equation}
    \label{eq:remainder_bound}
    \sup_{\theta\in\Theta_N} (\|\theta - \theta_0\|_{\infty}^4 + N\|\theta-\theta_0\|_{\infty}^2\|\theta-\theta_0\|_{C^1}) = o(N\epsilon_N^2),
\end{equation}
which is proven in Lemma \ref{lem:empirical_process_bounds}.

\begin{lemma}
    \label{lem:empirical_process_bounds}
    We have
    \eqref{eq:empirical_bound1}, \eqref{eq:empirical_bound2}, \eqref{eq:empirical_bound3} and
    \eqref{eq:remainder_bound}.
\end{lemma}

For proof, we use similar arguments as used to verify \cite[conditions~(20),~(21)]{monard2021statistical} and apply \cite[Theorem~3.5.4]{gine2021mathematical}.
We define two sequences
\begin{equation*}
    \barepsu := \sup_{\theta \in \Theta_N}\|\theta - \theta_0\|_{\infty},\;
    \barepsc := \sup_{\theta\in\Theta_N}
    \|\theta - \theta_0\|_{C^1},
\end{equation*}
and a Dudley-type integral
\begin{equation}
    \label{eq:dudley_int}
    J_{N,s}(u,t) := \int^u_0\sqrt{\log 2\mathcal{N}(\Theta_N,\; \norm_{H^s},\; t\delta)}~d\delta.
\end{equation}
The sequences $\barepsu$ and $\barepsc$ play similar roles as the terms ``$\sigma_N$'' and ``$\overline{\delta}_N$'' appearing in \cite[Condition~3.6]{monard2021statistical}.

\begin{proof}[Proof of \eqref{eq:remainder_bound}]
Since $\Theta_N\subset B_{H^p}(M)$ for $p < \alpha - d/2$, based on the definition of \eqref{eq:lsupset}, we have 
\begin{eqnarray*}
    \barepsu & \lesssim & \|G(\theta) - G(\theta_0)\|_{C^2(\Z)} \lesssim 
    \epsilon_N^{1-\frac{s+1}{p+1}},\;
    s > 1 + \frac{d}{2};\\
    \barepsc & \lesssim & \|\theta-\theta_0\|_{H^s(\Z)}
    \lesssim \|\theta - \theta_0\|_{\ltwo(\Z)}^{1-\frac{s}{p}}
    \lesssim \epsilon_N^{\frac{p-1}{p+1}
    \left(1 - \frac{s}{p}\right)},
\end{eqnarray*}
with uniform constants over $\Theta_N$ \eqref{eq:lsupset},
from interpolation, regularity estimate \eqref{eq:uniform_elliptic} and conditional stability estimates \eqref{eq:l2_inv_lipschitz} and \eqref{eq:lsup_inv_lipschitz}. 
Since $\alpha > 3 + \frac{3}{2}d$, we have $\frac{1}{2} > \frac{2+d/2}{\alpha + 1 - d/2}$, 
which shows that $\barepsu^4 \ll \epsilon_N^2$.
Also, under condition (25) of Theorem 4.3.1,
we can choose $s$ and $p$ so that
$2(1-\frac{s+1}{p+1})+\frac{p-1}{p+1}(1-\frac{s}{p}) > 2$,
which implies that $\barepsu^2\barepsc \ll \epsilon_N^2$,
proving claim \eqref{eq:remainder_bound}.
\end{proof}

\begin{proof}[Proof of \eqref{eq:empirical_bound1}]
For the first term from \eqref{eq:empirical_bound1}:
\begin{equation*}
\sup_{\theta\in\Theta_N} \left|\sum_{i=1}^{N}
\epsilon_i\{
G(\theta) - G(\theta_0) - \opinfo(\theta-\theta_0)\}(X_i)\right|,
\end{equation*}
we may view it as an empirical process 
with random variables $|\epsilon_i R_{\theta_0}(\theta - \theta_0)(X_i)|$ dominated by $F_N(\epsilon_i,X_i)=|\epsilon_i|\barepsc$ by estimates \eqref{eq:forward_lsup_lipschitz}, \eqref{eq:deriv_lsup_lipschitz} and $\E[\max_{1\leq i\leq N}|\epsilon_i|]\lesssim \sqrt{\log N}$.
For the application of \cite[Theorem~3.5.4]{gine2021mathematical}, one needs to bound from above the ``Koltchinskii-Pollard entropy integral'' appearing in display (3.169) therein:
\begin{equation*}
    \int_0^{\barepsu^2/\barepsc}\sup_Q
    \sqrt{\log 2\mathcal{N}(\Theta_N,\; \norm_{L^2(\Z,Q)},\; \delta\|F_N\|_{\ltwo(\Z,Q)})}~d\delta.
\end{equation*}
The notation $\ltwo(\Z,Q)$ indicates the $\ltwo$-norm is measured with respect to a measure $Q$ placed on $\Z$.
The supremum in $Q$ is taken over all discrete probability measures with finitely many atoms and rational probabilities.
This integral can be related to $J_N$ \eqref{eq:dudley_int} by the relation, for some fixed $c_0$:
\begin{equation*}
    \|\epsilon_i\{
    R(\theta_1,\theta_0)(X_i) - 
    R(\theta_2,\theta_0)(X_i)\}\|_{\ltwo(\Z,Q)}
    \leq c_0\frac{\|F_N\|_{\ltwo(\Z,Q)}}{\barepsc}
    \|\theta - \theta_0\|_{C^1}.
\end{equation*}
Since $H^p(\Z)$-norm is bounded uniformly over $\Theta_N$ for $p < \alpha - d/2$, 
we use a continuous embedding of $C^1(\Z)$ into $H^s(\Z),\; s > 1 + d/2$ and conclude for the entropy integral appearing in the above display can be bounded above by a constant multiple of the Dudley-type integral from \eqref{eq:dudley_int}:
\begin{equation*}
    J_{N,s}(\barepsu^2/\barepsc,\barepsc) =  \frac{1}{\barepsc}J_{N,s}(\barepsu^2,1),\; s > 1 + \frac{d}{2}.
\end{equation*}
Based on the definition \eqref{eq:dudley_int}, the problem is reduced to an upper bound on the metric entropy number $\log 2\mathcal{N}(\Theta_N,\; \norm_{H^s},\; t\delta)$.
We use known relations between approximation numbers for linear operators and the metric entropy (cf. \cite{edmunds1996function}, or a more condensed account from \cite[Appendix~B]{gugushvili2020bayesian}).
By definition, the $k$-th \emph{entropy number} $e_k$ of a linear operator $T: A\to B$ between Banach spaces $A,B$ is the infimum of all $\epsilon > 0$ such that $T(B_A(1)) \subset \cup_{i=1}^{2^k-1}\{v_i + \epsilon B_A(1)\}$ for some collection of $v_i$s in $B$.
If $T$ is compact, its relation to the metric entropy is given by (see \cite[Remark~1,~Section~1.3]{edmunds1996function})
\begin{equation*}
\log_2\mathcal{N}(T(B_A(1)),
\norm_B,\; e_k) = k-1.
\end{equation*}
Choose $A=H^p(\Z)$ and $B=H^s(\Z)$.
Based on the Besov space embeddings in \cite[Section~3.3.1]{edmunds1996function} and the first, upper bound theorem \cite[Section~3.3.2]{edmunds1996function}, we obtain the following metric entropy bound:
\begin{equation}
    \label{eq:entropy_bound}
    \log 2\mathcal{N}(\Theta_N,\;\norm_{H^s},\; t\delta) \leq 
    \log 2\mathcal{N}(B_{H^p(\Z)}(M),\;
    \norm_{H^s},\; t\delta) \leq
    A(d,s,p)
    \left(\frac{M}{t\delta}\right)^{\frac{d}{p-s}},\; p < \alpha - \frac{d}{2}.
\end{equation}
\eqref{eq:dudley_int} is thus convergent if $\alpha > 1 + \frac{3}{2}d$, which is implied by our assumption that $\alpha > 3 + \frac{3}{2}d$.

\cite[Theorem~3.5.4]{gine2021mathematical} now shows the left hand side of \eqref{eq:empirical_bound1} is bounded from above by 
\begin{equation*}
\sqrt{N} J_{N,s}(\barepsu^2,1)\vee 
\frac{\barepsc\sqrt{\log N} J_{N,s}^2(\barepsu^2,1)}{\barepsu^4}
\end{equation*}
By \eqref{eq:entropy_bound}, the first term is bounded above by a constant multiple of $\sqrt{N}\barepsu^{2-\frac{d}{2\zeta}},\; \zeta < \alpha - 1 - d$;
below, it will be shown that 
$\sqrt{N}\barepsc^{2-\frac{d}{2\zeta}}$ is $o(N\epsilon_N^2)$ under the maintained assumptions, so that the first term is also $o_N(N\epsilon_N^2)$, as $\barepsu\leq\barepsc$ for every $N$.
Next, the second term is bounded by a constant multiple of $\sqrt{\log N}\barepsc\barepsu^{-\frac{2d}{\zeta}},\; \zeta < \alpha - 1 - d$.
By using upper bounds on $\barepsu,\barepsu$ involving $\epsilon_N$, it can be checked that condition \eqref{eq:complex_condn2} implies this upper bound is also $o_N(N\epsilon_N^2)$.
\end{proof}

\begin{proof}[Proof of \eqref{eq:empirical_bound2}]
Next, for the second term \eqref{eq:empirical_bound2}:
\begin{equation*}
    \sup_{\theta\in\Theta_N} \left|
    \sum_{i=1}^{N} \{G(\theta) - G(\theta_0)\}^2(X_i) - 
    N\|G(\theta) - G(\theta_0)\|_{\ltwo(\Z)}^2
    \right|,
\end{equation*}
we may view it as an empirical process 
with random variables $\{G(\theta) - G(\theta_0)\}^2(X_i)$ dominated by $\barepsc^2$.
From the Lipschitz estimate for $G$ \eqref{eq:lsup_inv_lipschitz}, 
it can be shown that the following provides an upper bound on the left hand side of \eqref{eq:empirical_bound2}, 
due to a consequence of the same theorem (cf. \cite[Remark~3.5.5]{gine2021mathematical}):
\begin{equation}
    \label{eq:discretize_condn}
    \sqrt{N}\barepsc^2 
    J_{N,s}\left(1,\barepsc^2\right),
\end{equation}
which, by \eqref{eq:entropy_bound}, is bounded above by a constant multiple of $\sqrt{N}\barepsc^{2-\frac{d}{2\zeta}},\; \zeta < \alpha - 1 - d$.
It can be checked then that condition \eqref{eq:complex_condn1} implies this upper bound is $o_N(N\epsilon_N^2)$.
\end{proof}

\begin{proof}[Proof of \eqref{eq:empirical_bound3}]
Finally, for the third term \eqref{eq:empirical_bound3}:
\begin{equation*}
    \sup_{\theta\in \Theta_N^{(t)}}
    \left|\sum_{i=1}^{N}
    \{\opinfo(\theta - \theta_0)(X_i)\}^2 -
    N\| \opinfo(\theta - \theta_0) \|_{\ltwo(\Z)}^2
    \right|,
\end{equation*}
We can again use the Lipschitz estimate for the linear operator $\opinfo$ from \eqref{eq:deriv_lsup_lipschitz}.
We deduce that by Lemma \ref{lem:lsupnorm_bound} with the case of $r=1$,
plugging in the lower bound for $s_N$ from Lemma 5.2,
we have $|t|s_N^{-1}\|\barpsi\|_{C^1} \lesssim \barepsc$ for every $t\in\R$ if the following relation holds for $\tau_N=(\sqrt{N}\epsilon_N)^{-1}$:
\begin{equation*}
    \frac{1}{2}\geq \frac{\alpha+1}{2\alpha+2+d}\left\{
    \frac{\alpha-1-d/2}{\alpha+1-d/2}\frac{\alpha-1-d}{\alpha-d/2} +
    2\left( 1 - 
    \frac{(2+q)\alpha-1}{(2+q)\alpha+1}\frac{2\alpha-1-d/2}{2\alpha}
    \right)\right\}.
\end{equation*}
We have assumed $\alpha >3 + \frac{3}{2}d$.
By plugging in this lower bound, it can be checked the above inequality always holds for $q = 0$ and every $d \geq 1$, and thus holds for every $q \geq 0$ since the display is decreasing in $q$.
Therefore, we conclude that indeed $|t|s_N^{-1}\|\barpsi\|_{C^1} \lesssim \barepsc$, so that similar calculations as done for \eqref{eq:empirical_bound2} can be now done, and it again then suffices to show that \eqref{eq:discretize_condn} is true.
\end{proof}

\noindent {\bf 9. Showing \eqref{eq:lan_ratio_control}.}
We define a new ``linearized posterior'' $\dnlan = \dnlan(\cdot \mid \data)$ as
\begin{eqnarray*}
    \dnlan(A\mid \data) & = & \frac{\int_A e^{-\dllan(\theta - \theta_0)}~d\nprior(\theta)}
    {\int e^{-\dllan(\theta - \theta_0)}~d\nprior(\theta)},\\
    \dllan(\theta - \theta_0) & = & -\sum_{i=1}^{N}
    \left[ 
    \epsilon_i\opinfo(\theta-\theta_0)(X_i) - \frac{\{\opinfo(\theta-\theta_0)(X_i)\}^2}{2}\right].
\end{eqnarray*}
Then, given \eqref{eq:empirical_bound3},  we have
\begin{equation*}
    \left| \frac{\nlan\{\Theta_N - ts_N^{-1}\barpsi\mid \data\}}
    {\nlan\{\Theta_N \mid \data\}} - 1\right|
    \leq e^{o(N\epsilon_N^2)}
    \left| \frac{\dnlan\{ \Theta_N - ts_N^{-1}\barpsi| \data \}}{\dnlan\{\Theta_N | \data \}} - 1\right|.
\end{equation*}
Therefore, to verify \eqref{eq:lan_ratio_control},
it suffices to show that there exists some $c' > 0$ for which
\begin{equation*}
    \ptheta\left(
    \left| \frac{\dnlan\{ \Theta_N - ts_N^{-1}\barpsi| \data \}}{\dnlan\{\Theta_N | \data \}} - 1\right| \geq e^{-c'N\epsilon_N^2}\right) \to 0.
\end{equation*}
Intuitively, this new distribution $\dnlan(d\theta|\data)$ is a Bayes posterior naturally arising from a ``linearized'' statistical experiment, in which we observe $N$ independent observations $Y_i = \opinfo(\theta_0)(X_i) + \epsilon_i$ with the underlying measure space associated with a $\sigma$-algebra generated by $(X_i,\epsilon_i)_{i=1}^{N}$.

We first argue that we can show this posterior has the same contraction rate upper bound in $\ltwo(\Z)$-norm as $\Pi(\cdot\mid\data)$, i.e., with some $b' > 0$:
\begin{equation}
    \label{eq:lan_contraction}
    \ptheta(\dnlan\{\Theta_N\mid\data\} \geq e^{-b'N\epsilon_N^2}) \to 0.
\end{equation}
By the Lipschitz continuity estimate on $\opinfo : \ltwo(\Z) \to \ltwo(\Z)$ given by \eqref{eq:forward_lipschitz_deriv}, 
we know that this operator satisfies the same local Lipschitz estimate as $G$ \eqref{eq:forward_lipschitz} (with possibly different constants),
and we also know that $\|\opinfo(\theta)\|_{\infty} \leq \|\theta\|_{C_1}$ is bounded whenever $\theta\in B_{\Hs^s}(M)$ for some suitable $s$. As remarked in, up to multiplicative constants that may be chosen uniform over $\theta\in B_{\Hs}(M)$, such estimates are essentially the same as obtained for the nonlinear operator $G$.
Applying \cite[Theorems~2.2.2~and~2.3.1]{nickl2023bayesian} to $\dnlan$, then, we may conclude the posterior contraction \eqref{eq:lan_contraction} is indeed true, and that the sequence $\epsilon_N\to 0$ can be chosen exactly the same as the original posterior $\nprior(\cdot\mid\data)$ in Assumption \ref{assmp:derivative}.

Next, using Lemma \ref{lem:lsupnorm_bound} and plugging in the lower bounds for $s_N/\tau_N$ from Lemma \ref{lem:order_s},
we obtain the estimate, for each $t\in\R$ and with our choice of $\tau_N$,
\begin{equation*}
    |t|s_N^{-1}\|\barpsi\|_{\ltwo(\Z)}
    \leq |t|\frac{1}{\sqrt{N}\epsilon_N}
    \epsilon_N^{(2\kappau+1)\delta(2+q) - 2\kappau},
\end{equation*}
where $\delta(2+q) = \frac{(2+q)\alpha - 1}{(2+q)\alpha +1}$ by Lemma \ref{lem:darcy_regularity}.
The obtainable estimate of $\kappau$ is different depending on the three different regimes in Lemma \ref{lem:order_s}; however,
it is bounded above by 1/2, and the exponent $(2\kappau+1)\delta(2+q) - 2\kappau$ is decreasing in $\kappau\in [0,1]$ when $\delta(2+q) < 1$. Therefore, an upper bound follows:
\begin{equation*}
    |t|s_N^{-1}\|\barpsi\|_{\ltwo(\Z)}\lesssim \frac{1}{\sqrt{N}}\epsilon_N^{2\{\delta(2+q) - 1\}}.
\end{equation*}
Now we plug in the definition of $\epsilon_N = N^{-\frac{\alpha+1}{2\alpha+2+d}}$.
By the conditional stability estimate \eqref{eq:l2_inv_lipschitz}, we see that $\|\theta-\theta_0\|_{\ltwo(\Z)}\lesssim \epsilon_N^\lambda$ for $\lambda < \frac{\alpha - d/2 - 1}{\alpha - d/2 + 1}$, uniformly over the set $\Theta_N$ given in \eqref{eq:diffusion_eq}. The above display is $o_N(\epsilon_N^\lambda)$ for every $\alpha\geq 1 + d/2$, $d\geq 1$ and $q\geq 0$,
since we obtain
\begin{equation*}
    \frac{1}{2} > \frac{\alpha+1}{2\alpha+2+d}\left\{\frac{\alpha-d/2-1}{\alpha-d/2+1} + \frac{4}{(2+q)\alpha+1}\right\}.
\end{equation*}
Similarly, again using Lemma \ref{lem:lsupnorm_bound},
we see that there exists $s > 1 + d/2$ for which we can obtain $s_N^{-1}\|\barpsi\|_{H^s(\Z)} = o_N(1)$.
Therefore, we conclude that \eqref{eq:lan_contraction} remains true if $\Theta_N$ is replaced with the shifted set $\Theta_N - ts_N^{-1}\barpsi$ for each $t$.

Combining the above estimates and \eqref{eq:lan_contraction}, we conclude that for each $t\in\R$, there exist $c_1,c_2 > 0$ such that the event
\begin{equation*}
    \left| \frac{\dnlan\{ \Theta_N - ts_N^{-1}\barpsi| \data \}}{\dnlan\{\Theta_N | \data \}} - 1\right|
    \leq \frac{e^{-c_2N\epsilon_N^2}}{1-e^{-c_1N\epsilon_N^2}}
\end{equation*}
has $\ptheta$-probability convergent to 0.
Hence, condition \eqref{eq:lan_ratio_control} is verified.

\begin{lemma}
    \label{lem:lsupnorm_bound}
    For every $t > 0$ and $0 \leq r < \alpha - d$:
    \begin{eqnarray*}
        ts_N^{-1}\|\barpsi\|_{\ltwo(\Z)}
         & \lesssim & t\tau_N\left(\frac{1}{\sqrt{N}\tau_N}\right)^{\delta(2+q)}
        \left(\frac{s_N}{\tau_N}\right)^{\delta(2+q)-1}
        ;\\
        ts_N^{-1}\|\barpsi\|_{C^r(\Z)} 
        & \lesssim & t\tau_N
        \left(\frac{1}{\sqrt{N}\tau_N}\right)^{\delta(2+q)(1-\frac{p}{2})}
        \left(\frac{s_N}{\tau_N}\right)^{\delta(2+q)(1-\frac{p}{2}) - 1},\;
        p > \frac{r}{\alpha} + \frac{d}{2\alpha}.
    \end{eqnarray*}
\end{lemma}

\begin{proof}
    For the first statement, we proceed similarly as in proving \eqref{eq:bias_ineq1} and \eqref{eq:bias_ineq2} in the proof of Lemma \ref{lem:order_sb}, to conclude
    \begin{equation*}
        \|\barpsi\|_{\ltwo(\Z)} =
        \|f^2(\opcomb^*\opcomb)\lpsi\|_{\Hs^1}
        \leq \|\opcomb f^2(\opcomb^*\opcomb)
        \lpsi\|_{\ltwo(\Z)}^{\delta(2+q)}
    \end{equation*}
    and therefore
    \begin{equation*}
        ts_N^{-1}\|\barpsi\|_{\ltwo(\Z)} \lesssim t\tau_N\left(\frac{1}{\sqrt{N}\tau_N}\right)^{\delta(2+q)}
        \left(\frac{s_N}{\tau_N}\right)^{\delta(2+q)-1}.
    \end{equation*}
    For the second statement,
    since the embedding from $H^\alpha(\Z)$ into $H^s(\Z)$ is Hilbert-Schmidt when $s < \alpha - d/2$,
    we deduce that the prior places full measure on $\Hs^s(\Z)$ for $s < 1 - \frac{d}{2\alpha}$. 
    Since we also know that $H^{s\alpha}(\Z)$ is continuously embedded into $C^{r\alpha}(\Z)$ for $r < s - \frac{d}{2\alpha}$, 
    we conclude that
    \begin{equation*}
        \|\barpsi\|_{C^r} \lesssim \|\barpsi\|_{\Hs^{p}} = \tau_N^2
        \|f^2(\opcomb^*\opcomb)\lpsi\|_{\Hs^{p-1}},
    \end{equation*}
    for any $p > \frac{r}{\alpha} + \frac{d}{2\alpha}$, and $f(t) = \frac{1}{\sqrt{1 + N\tau_N^2t^2}}$,
    with $C^0(\Z) = (C(\Z), \norm_{\infty})$.
    From hereon, proceeding similarly as in the proof of \eqref{eq:cauchy_schwarz}, \eqref{eq:bias_ineq1} and \eqref{eq:bias_ineq2},
    we obtain
    \begin{equation*}
        ts_N^{-1}\|\barpsi\|_{C^r(\Z)} \lesssim
        t\tau_N
        \left(\frac{1}{\sqrt{N}\tau_N}\right)^{\delta(2+q)(1-\frac{p}{2})}
        \left(\frac{s_N}{\tau_N}\right)^{\delta(2+q)(1-\frac{p}{2}) - 1}.
    \end{equation*}
\end{proof}

\section{Technical Results for Darcy Flow}
\label{ssec:darcy_lemmas}

As in the main text, we consider the elliptic PDE model given by \eqref{eq:diffusion_eq}, 
with its conductivity parameter $a$ belonging to the set $\Gamma_{\alpha,K_{min}}$ \eqref{eq:constraint_set}, and
known source $f$ and boundary condition $u_D$ belonging to $C^\infty(\Z)$.
We consider a smooth reparameterization $\phi : H_0^\alpha(\Z) \to \Gamma_{\alpha,K_{min}}$, in the sense of Definition \ref{def:link_fn}. 
The forward map $G(\theta) = u_{\phi\circ\theta}$ uniquely defines the solution with $a=\phi\circ\theta$.
Below, we summarize some known results about the PDE that we use in our proofs.

\begin{lemma} \label{lem_supp1}
    Let $\theta\in H_0^s(\Z),\;s>1+d/2$. For each $M>0$, there exists a constant \\
    $C \equiv C(s,d,\Z,K_{\min},f,u_D,M) > 0$ such that
    \begin{equation} \label{eq:uniform_elliptic}
        \sup_{\|\theta\|_{H^s(\Z)} \leq M}\| G(\theta) \|_{H^{s+1}(\Z)} \leq C < \infty.
    \end{equation}
\end{lemma}
See \cite[Proposition~A.5.2]{nickl2023bayesian} for a proof and the explicit expression of the constant.

\begin{lemma} \label{lem_supp2}
For $a = \phi(\theta) \in C^1(\Z)$ and $a \geq K_{min}$ on the domain,    
\begin{equation}
    \label{eq:bound_elliptic}
    \|G(\theta)\|_{\infty} \leq c(\Z,K_{min})
    \|f\|_{\infty} + \|u_D\|_{\infty}.
\end{equation}
\end{lemma}
The proof can be derived by a Feynman-Kac representation of the solution $G(\theta)$; 
see, e.g., the remarks leading up to display (A.27) in \cite{nickl2023bayesian}.

\begin{lemma} \label{lem_supp3}
    For $\theta_1,\theta_2\in B_{H_0^s(\Z)}(M),\;s>1+d/2$, 
    \begin{equation}
        \label{eq:forward_lipschitz}
        \|G(\theta_1) - G(\theta_2)\|_{\ltwo(\Z)} \lesssim \|\theta_1 - \theta_2\|_{(H_0^1)^*}
    \end{equation}
    with a multiplicative constant uniform in $\theta\in B_{H_0^s(\Z)}(M)$.
    Furthermore, 
    \begin{equation}
    \label{eq:forward_lsup_lipschitz}
    \|G(\theta_1) - G(\theta_2)\|_{\ltwo(\Z)} \leq c(\Z,K_{min})\|\theta_1-\theta_2\|_{C^1(\Z)}.
    \end{equation}
\end{lemma}

\begin{proof}
    We can utilize the fact that two solutions $u_{\theta}$ and $u_{\theta_0}$, corresponding respectively to elliptic PDE coefficients $\phi(\theta)$ and $\phi(\theta_0)$,
    solve an inhomogeneous PDE
    \begin{equation}
    \label{eq:inhomogeneous_pde}
    \begin{cases}
    (L_{\theta}u_{\theta} -L_{\theta_0}u_{\theta_0}) = 
    L_{\theta}(u_\theta - u_{\theta_0})\quad\text{on }\Z; \\
    0 = u_\theta - u_{\theta_0}\quad\text{on }\partial\Z.
    \end{cases}
    \end{equation}
    The second estimate then follows using the Feynman-Kac representation as in \ref{lem_supp2}.
    For the proof of the first claim, see \cite[Section~2.1.1.3]{nickl2023bayesian}.
\end{proof}

For linearizing $G$, arguing as in the proof of \cite[Theorem~3.3.2]{nickl2023bayesian}, we can show that the following bounded operator for any $\theta \in H^s(\Z)$,
mapping $h$ from $H^{s}(\Z),\; s > 1 + \frac{d}{2}$ to $\ltwo(\Z)$:
\begin{equation*}
    \mathbb{I}_{\theta}(h) = -L_{\theta}^{-1}[\nabla\cdot \{(\phi'\circ\theta)h\nabla u_{\theta}\}],
\end{equation*}
can be continuously extended to $\ltwo(\Z)$ and satisfies a quadratic approximation in $\ltwo(\Z)$-norm, as shown in \cite[Theorem~3.3.2]{nickl2023bayesian}:
\begin{equation}
    \label{eq:quadratic_approx}
    \|R_{\theta}(h)\|_{\ltwo(\Z)} = O(\|h\|_{\infty}^2).
\end{equation}
The adjoint of this operator $\mathbb{I}^*_{\theta} : \ltwo(\Z) \to \ltwo(\Z)$ for $\theta\in H^s(\Z),\; s > 1 + \frac{d}{2}$ is given by
\begin{equation*}
    \mathbb{I}_{\theta}^*(g) = (\phi'\circ\theta)\nabla u_{\theta}\cdot\nabla 
    L_{\theta}^{-1}[g].
\end{equation*}
In particular, the choice of link function $\phi(\theta) = e^{\theta} + K_{min}$ yields $\phi'(\theta) = e^{\theta}$, 
even though this link function does \emph{not} satisfy Definition 4.1.

\begin{lemma}\label{lem_supp4}
    Let $\theta_1,\theta_2\in B_{H_0^s}(M),\;s>1+d/2$, and $\theta_1=\theta_2$ on $\partial\Z$. 
    Then
    \begin{equation}
    \label{eq:forward_lipschitz_deriv}
    \| \opinfo(\theta_1 - \theta_2) \|_{\ltwo(\Z)}
    \lesssim \|\nabla\cdot 
    (\phi'(\theta_0)
    (\theta_1 - \theta_2)\nabla u_{\theta_0})\|_{(H_0^2)^*}
    \lesssim \|\theta_1 - \theta_2\|_{(H_0^1)^*}
\end{equation}
with a multiplicative constant uniform in $\theta\in B_{H_s^0}(M)$.
\end{lemma}
\begin{proof}
    By Sobolev embedding, $\theta_1,\;\theta_2\in B_{C^1}(M)$ and they agree on the boundary $\partial\X$ by the assumption. We then also know that $(H^1)^*$-norm and the $(H_0^1)^*$-norm are equivalent for $\theta_1 - \theta_2$.
    The Lipschitz estimate then follows a similar argument as in the initial estimates from the proof of \cite[Theorem~9]{nickl2020convergence}: 
    see the arguments in pp. 390--391, displays (4.10), (4.11) and the Lipschitz estimate right below it.
\end{proof}

\begin{lemma}\label{lem_supp5}
    In the setting of \ref{lem_supp4}, we also have
    \begin{equation}
    \label{eq:deriv_lsup_lipschitz}
    \|\mathbb{I}_\theta(g-h)\|_{\infty} \lesssim
    \|g-h\|_{C^1}.
\end{equation}
\end{lemma}

\begin{proof}
    Since $\phi$ realizes a bijection from $C^1(\Z)$ to an image contained in $C^1(\Z)$, and $u_{\theta}$ has $C^{2+\xi}(\Z)$ regularity for some $\xi > 0$, 
    we deduce that $\phi'(\theta),\; \nabla\phi'(\theta),\;\nabla u_{\theta},\;\Delta u_{\theta}$ are all bounded functions on $\X$.
    The Lipschitz estimate then follows using the definition of $\mathbb{I}_\theta$ derived above, with a multiplicative constant depending on the norms of
    $\phi'(\theta),\; \nabla\phi'(\theta),\;\nabla u_{\theta},\;\Delta u_{\theta}$
    measured in $\norm_\infty$, 
    which may be made uniform in $\theta_1,\theta_2\in B_{C^1}(M)$.
\end{proof}

Comparison of \eqref{eq:forward_lipschitz_deriv}
and \eqref{eq:deriv_lsup_lipschitz} with \eqref{eq:forward_lipschitz} and \eqref{eq:forward_lsup_lipschitz}
shows that the nonlinear operator $G$ and its linearization $\opinfo$ essentially have the same degree of smoothness
when restricted to balls in $C^1(\Z)$.

\begin{lemma}\label{lem_supp6}
    Let $\theta,\theta_0\in B_{H^s(\Z)}(M),\;s>1+d/2$ with $\theta=\theta_0$ on $\partial\Z$. Suppose Assumption \ref{assmp:darcy_flow} is met.
    Then
    \begin{equation}
    \label{eq:l2_inv_lipschitz}
    \|\theta-\theta_0\|_{\ltwo(\Z)} \leq C(M)\|G(\theta) - G(\theta_0)\|_{H^2(\Z)};
    \end{equation}
    furthermore,
    \begin{equation}
    \label{eq:lsup_inv_lipschitz}
    \|\theta - \theta_0\|_{\infty} \leq 
    C(M)\|u_{\theta} - u_{\theta_0}\|_{C^2(\X)}.
\end{equation}
\end{lemma}

\begin{proof}
    See \cite[Proposition~2.1.7]{nickl2023bayesian} for the proof of the first assertion, which relies on the fact that $G(\theta)-G(\theta_0)$ solves the inhomogeneous PDE \eqref{eq:inhomogeneous_pde}. 
    For the second assertion, under we can apply \cite[Theorem~1]{richter1981inverse}
    to this inhomogeneous PDE since Assumption \ref{assmp:darcy_flow} is met.
\end{proof}

\section{Numerical Experiments} \label{ssec:simulation}

\begin{figure}[!tb]
    \centering
    \begin{tabular}{ccc}
        \captionsetup[subfigure]{oneside,margin={-1cm,0cm}}
        \begin{subfigure}{0.45\linewidth}
            \centering
            \hspace*{-1cm}
            \includegraphics[trim={0.8cm 2.3cm 1.5cm 1.8cm},clip,
            width=\linewidth]{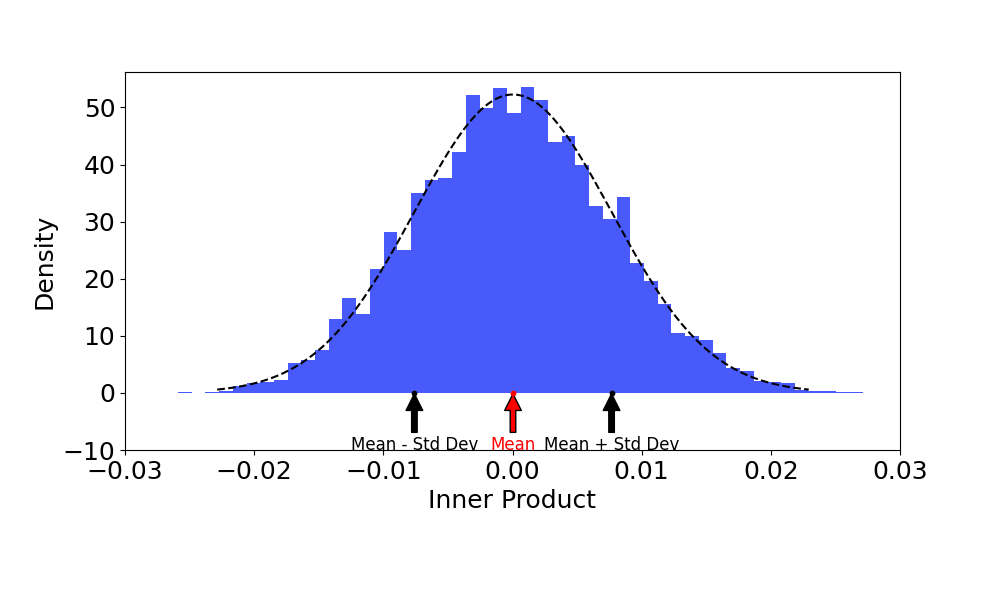}
            \caption{$\sqrt{N}= 20$}
        \end{subfigure} &
        \captionsetup[subfigure]{oneside,margin={-1cm,0cm}}
        \begin{subfigure}{0.45\linewidth}
            \centering
            \hspace*{-1cm}
            \includegraphics[
            trim={0.8cm 2.3cm 1.5cm 1.8cm},clip,
            width=\linewidth]{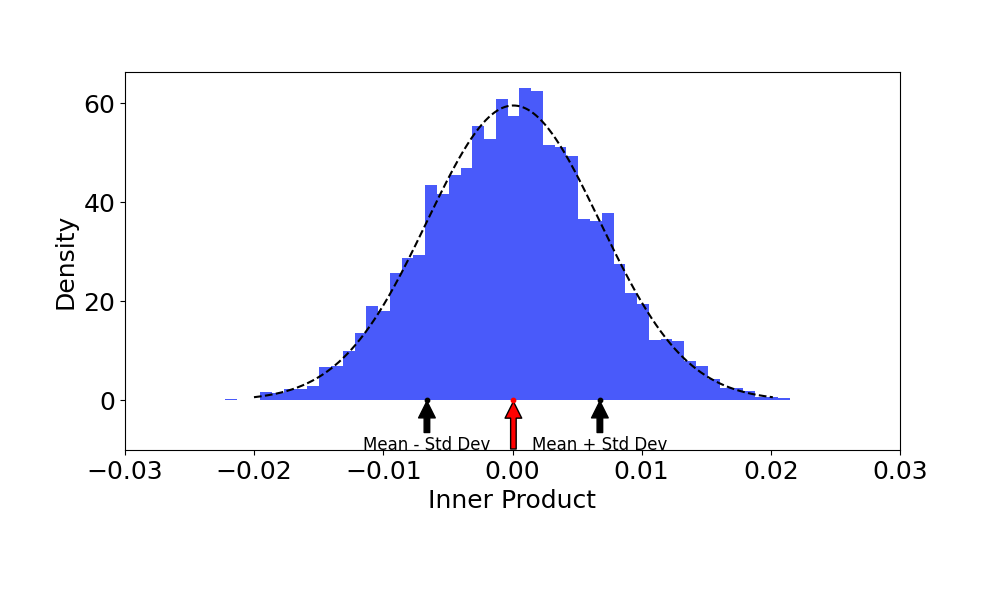}
            \caption{$\sqrt{N} = 40$}
        \end{subfigure} \\
        \captionsetup[subfigure]{oneside,margin={-1cm,0cm}}
        \begin{subfigure}{0.45\linewidth}
            \centering
            \hspace*{-1cm}            
            \includegraphics[
            trim={0.8cm 2.3cm 1.5cm 1.8cm},
            clip,width=\linewidth]{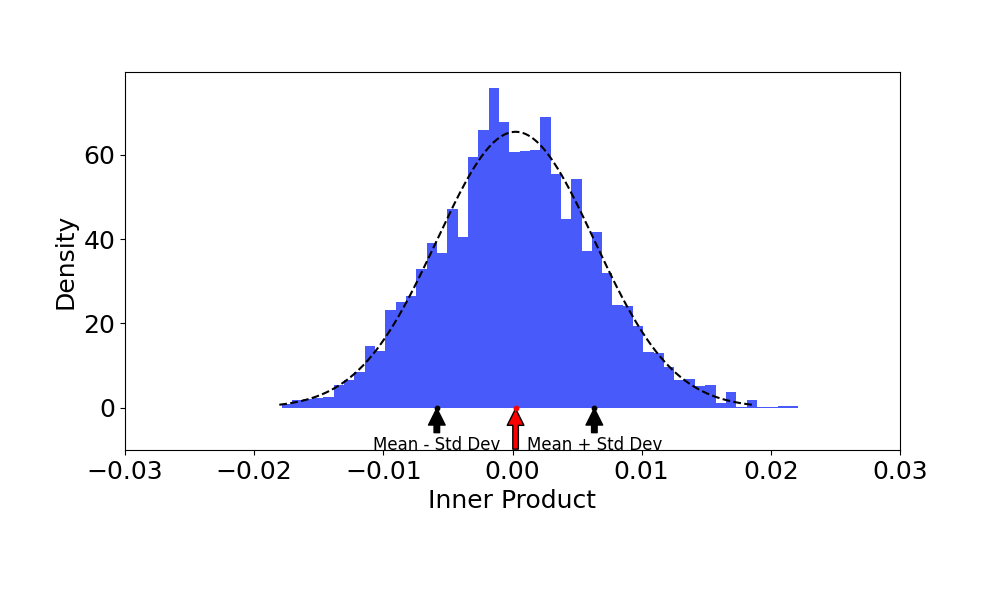}
            \caption{$\sqrt{N} = 80$}
        \end{subfigure}   
    \end{tabular}
    \caption{
        Histograms of posterior draws $\langle \psi, \theta_s \rangle,\; s = 1,2,\ldots,S$ with $S = 10,000$ (number of MCMC draws), corresponding to different number of observations $N$.
        Posterior mean (red arrow) and a credible interval of 1 posterior standard deviation width (black arrows) are shown. The truth is at the origin ($\langle\theta_0,\psi\rangle = 0$).
    }
\label{fig:figure}
\end{figure}

We illustrate our theory by considering an example similar to the one first studied by \cite{nickl2022some}, based on a two-dimensional Darcy flow. 
Our emphasis is not on a practical demonstration of implementing the MCMC algorithm; rather, it is to consider an artificial but theoretically non-trivial example with a large but finite number $N$, to which the BvM results of \cite{monard2021statistical} and \cite{nickl2023bayesian} are inapplicable.
Consider the equation
\begin{equation} \label{eq:darcy_flow}
    -\nabla \cdot (a(\theta_0) \nabla u(x, y)) = f(x,y).
\end{equation}
Let $\theta_0 = 0$ and $a_0 = e^{\theta_0} = 1$. The resulting equation with Dirichlet boundary condition is an inhomogeneous Poisson equation
\begin{equation*}
    \begin{cases}
        -\Delta u(x, y) = f(x, y),\; &(x,y)\in (0,1)\times(0,1); \\
        u(x, y) = 0,\; &(x,y)\in [0,1]\times[0,1]\setminus(0,1)\times(0,1).
    \end{cases}
\end{equation*}
We choose $f = -2x(1-x) - 2y(1-y)$. One can check that the unique solution is $u(x, y) = -xy(1-x)(1-y)$. Consider the following bump function defining a functional of inferential interest:
\begin{equation}\label{bump}
\psi(x, y) = \begin{cases}
    \exp\left\{-
    \left(\frac{1}{(8x-3)(5-8x)} + 
    \frac{1}{(5y-1)(3-5y)}\right)\right\},\; 
    &(x,y)\in \left(\frac{3}{8} , \frac{5}{8}\right)\times
    \left(\frac{1}{5} , \frac{3}{5}\right); \\
    0 &\text{otherwise}.
\end{cases}
\end{equation}
The fact that $\theta_0$, an unknown parameter, is actually zero, and that $\psi\in C_0^\infty(\Z)$, seems to imply that the problem is trivial and of little practical interest. 
However, an integral curve argument similar to the one found in \cite{nickl2022some} shows that because $\psi\notin R(\opinfo^*)$, a necessary condition for BvM is \emph{not} met. 
Therefore, despite the apparent triviality of the estimation problem, previous literature does not provide  a positive statement about the posterior coverage of $\langle\psi,\theta_0\rangle_{\ltwo(\Z)}$. 
On the other hand, since $\theta_0,\psi\in C^\infty_0(\Z)$ and $\theta_0 = 0$, we clearly have $|b_N|/s_N = 0\ll 1$ for $b_N$ \eqref{eq:bias} and $s_N$ \eqref{eq:asymp_scale},
and the rest of the sufficient conditions on $\alpha,\beta$ in Theorem \ref{thm:darcy_coverage} are also trivially met. 
We thus expect that for reasonable Gaussian prior specifications, the Bayes credible interval covers $\langle\psi,\theta_0\rangle$ with high probability for large enough $N$.

For the numerical study, we first generated discrete observations $G(\theta_0)(X_i)$ over a $\sqrt{N}\times\sqrt{N}$ square grid over $[0,1]\times[0,1]$, with $N\in\{20,40,80\}$. For each grid configuration, we corrupted each $G(\theta_0)(X_i)$ with a Gaussian additive noise with standard deviation $\sigma = 5.0$ to simulate each $Y_i$.
The posterior samples given a Gaussian likelihood centered around $G(\theta)$ were drawn using the preconditioned Crank-Nicolson (pCN) Metropolis-Hastings algorithm introduced by \cite{cotter2013mcmc} for $S = 10000$ iterations targeting the posterior distribution of $\theta|\data$.
The tuning parameter $\beta$, controlling the dependence between the previous iterate of the Markov chain and the new proposal, was fixed at 0.99. Independent Gaussian proposals were used, with variance tuned to aim at an acceptance ratio between 50 and 60\%. 

In Figure \ref{fig:figure}, we show the histograms of the tracked quantities $\langle\theta_s, \psi\rangle,\; s = 1,2,\ldots,S$, with $\theta_s$ denoting the $s$-th sample drawn from MCMC.
Three plots of posterior sample histograms for the target functional $\langle\psi,\theta\rangle$ are shown, each corresponding to an increasing number of observed data points $N$.
These results demonstrate that the posterior samples appear approximately Gaussian, that the credible intervals cover the origin (the truth), and that there is a reduction in variance with growing $N$. 

\end{document}